\documentclass[11pt]{amsart}
\usepackage{rotating}
\usepackage{dynkin-diagrams}
\usepackage{etoolbox}
\usepackage[flushleft]{threeparttable}
\usepackage{hhline}
\usepackage{comment}
\usepackage{tikz}
\usepackage[page]{appendix}
\usepackage[utf8]{inputenc}
\usepackage{filecontents}
\usepackage{graphicx} 
\usepackage{mathtools}
\usepackage{etoolbox}
\usepackage{amsmath}
\usepackage{amsthm}
\usepackage{amsfonts}
\usepackage{amssymb}
\usepackage{mathrsfs}
\usepackage{xypic}
\usepackage{color}
\usepackage{xcolor}
\usepackage{titlesec}
\usepackage{titletoc}
\usepackage{etoc}
\usepackage{minitoc}
\usepackage{tocloft}
\usepackage{tabularx}

\newtheorem{thm}{Theorem}
\newtheorem{conj}{Conjecture}
\newtheorem*{quest*}{Question}

\newtheorem*{thm*}{Theorem}
\newtheorem{dfn}{Definition}

\usepackage[style=alphabetic,firstinits,backend=bibtex, uniquename=init,
url=true,hyperref]{biblatex}
\DeclareFieldFormat{labelalpha}{\thefield{entrykey}}
\DeclareFieldFormat{extraalpha}{}
\DeclareFieldFormat{postnote}{#1}
\DeclareFieldFormat{multipostnote}{#1}
\usepackage[colorlinks]{hyperref}
\definecolor{cclr}{rgb}{25,25,112}
\hypersetup{
citecolor=[rgb]{0.15, 0.15, 0.68},
linkcolor=gray,   
urlcolor=magenta}

\usetikzlibrary{decorations.markings}
\usetikzlibrary{decorations.shapes}
\usetikzlibrary{shapes.geometric}
\usetikzlibrary{decorations.fractals}
\usetikzlibrary{decorations.footprints}

\tikzset {
    , bedknobs one/.style = {
        decorate
        , fill = red!50
        , decoration = {
            shape backgrounds
            , shape = #1
            , shape size = 2mm
        }
    }
    , bedknobs two/.style = {
        decorate
        , decoration = {
            #1
        }
        , fill = blue!50
    }
}

\nocite{*}
\bibliography{\jobname}

\newcommand{\DL}{\operatorname{DL}}
\newcommand{\JH}{\operatorname{JH}}

\newcommand{\Aut}{\operatorname{Aut}}
\newcommand{\Int}{\operatorname{Int}}
\newcommand{\red}{\operatorname{red}}

\newcommand{\crys}{\operatorname{crys}}

\newcommand{\sems}{\operatorname{ss}}
\newcommand{\std}{\operatorname{std}}

\newcommand{\Res}{\operatorname{Res}}

\newcommand{\dirlim}[1]{\underset{#1}{\underrightarrow{\operatorname{lim}}}}
\newcommand{\invlim}[1]{\underset{#1}{\underleftarrow{\operatorname{lim}}}}
\usepackage{xstring}

\titleformat{\section}[runin]{\normalfont\bfseries}{\thesection.}{3pt}{}
\titleformat{\subsection}[runin]{\normalfont\bfseries}{\thesubsection.}{3pt}{}
\titleformat{\subsubsection}[runin]{\normalfont\bfseries}{\thesubsubsection.}{3pt}{}
\renewcommand{\thesection}{\arabic{section}}
\titleformat{\section}{\normalfont\large\bfseries}{\thesection.~~}{1em}{}
\DeclareMathAlphabet{\mathpzc}{OT1}{pzc}{m}{it}

\dottedcontents{section}[2.5em]{}{2.9em}{1pc}
\newcommand{\mychapter}[1]{\section{#1}}
\newcommand{\mysection}[1]{\subsection{#1}}
\newcommand{\mysubsection}[1]{\subsubsection{#1}}

\newcommand{\myappsection}[1]{\section{#1}}
\newcommand{\myappsubsection}[1]{\subsubsection{#1}}

\begin{document}

\newcommand{\chicyc}{\chi_{\operatorname{cyc}}}
\newcommand{\Brnr}{\operatorname{Br}_{\operatorname{nr}}}
\newcommand{\Bronr}{\operatorname{Br}_{1,\operatorname{nr}}}
\newcommand{\Br}{\operatorname{Br}}
\newcommand{\LT}{\operatorname{LT}}
\newcommand{\Frob}{\operatorname{Frob}}
\newcommand{\Fr}{\operatorname{Fr}}
\newcommand{\Qp}{\mathbb{Q}_p}
\newcommand{\Zp}{\mathbb{Z}_p}
\newcommand{\TODO}{{\color{red} TODO}}
\newcommand{\Gm}{\mathbb{G}_m}
\newcommand{\Scin}{\operatorname{Scin}}
\newcommand{\Fil}{\operatorname{Fil}}
\newcommand{\Ind}{\operatorname{Ind}}
\newcommand{\Sym}{\operatorname{Sym}}
\newcommand{\sym}{\operatorname{sym}}
\newcommand{\semis}[1]{#1\operatorname{-ss}}
\newcommand{\Alt}{\operatorname{Alt}}
\newcommand{\DGK}{\operatorname{DG}_K}
\newcommand{\Img}{\operatorname{Im}}
\newcommand{\Ker}{\operatorname{Ker}}
\newcommand{\grp}{\operatorname{grp}}
\newcommand{\cont}{\operatorname{cont}}
\newcommand{\ord}{\operatorname{ord}}
\newcommand{\Hom}{\operatorname{Hom}}
\newcommand{\Ext}{\operatorname{Ext}}
\newcommand{\Ad}{\operatorname{Ad}}
\newcommand{\Id}{\operatorname{id}}
\newcommand{\Lie}{\operatorname{Lie}}
\newcommand{\Lift}{\operatorname{Lift}}
\newcommand{\ad}{\operatorname{ad}}
\newcommand{\rk}{\operatorname{rk}}
\newcommand{\la}{\operatorname{l. alg.}}
\newcommand{\Det}{\operatorname{Det}}
\newcommand{\LHS}{\operatorname{LHS}}
\newcommand{\RHS}{\operatorname{RHS}}
\newcommand{\bFp}{{\bar{\mathbb{F}}_p}}
\newcommand{\bZp}{\bar{\mathbb{Z}}_p}
\newcommand{\bQp}{\bar{\mathbb{Q}}_p}
\newcommand{\uF}{{\breve{F}}}
\newcommand{\uS}{{\underline{S}}}
\newcommand{\uchi}{{\underline{\chi}}}
\newcommand{\uT}{{\underline{T}}}
\newcommand{\uGm}{{\underline{\Gm}}}
\newcommand{\uB}{{\underline{B}}}
\newcommand{\uG}{{\underline{G}}}
\newcommand{\kNm}[1]{{\widehat{\kappa_{#1}^\times}^{\operatorname{Nm}}}}
\newcommand{\ONm}[1]{{\widehat{\cO_{#1}^\times}^{\operatorname{Nm}}}}
\newcommand{\bF}{{\mathbb{F}}}
\newcommand{\cC}{{\mathcal{C}}}
\newcommand{\cZ}{{\mathcal{Z}}}
\newcommand{\bfF}{{\mathbf{F}}}
\newcommand{\bZ}{{\mathbb{Z}}}
\newcommand{\bA}{{\mathbb{A}}}
\newcommand{\bQ}{{\mathbb{Q}}}
\newcommand{\bR}{{\mathbb{R}}}
\newcommand{\bC}{{\mathbb{C}}}
\newcommand{\ucong}[1]{\overset{#1}{\cong}}
\newcommand{\scrC}{{\mathscr{C}}}
\newcommand{\Fq}{\bar{\mathbb{F}}_q}
\newcommand{\Sems}{\operatorname{ss}}
\newcommand{\MinR}{\operatorname{mr}}
\newcommand{\GSp}{\operatorname{GSp}}
\newcommand{\Sp}{\operatorname{Sp}}
\newcommand{\GL}{\operatorname{GL}}
\newcommand{\Sha}{\operatorname{Sha}}
\newcommand{\SC}[1]{\mathfrak{SC}_{#1}}
\newcommand{\TI}[1]{\mathfrak{TI}_{#1}}
\newcommand{\SCu}[1]{\mathfrak{SC}_{#1,\operatorname{ur}}}
\newcommand{\QZp}{\bQ_{p'}/\bZ}
\newcommand{\LS}{\operatorname{LS}}
\newcommand{\Tors}{\operatorname{tors}}
\newcommand{\Tf}{\operatorname{tf}}
\newcommand{\BC}{\operatorname{BC}}
\newcommand{\Nm}{\operatorname{Nm}}
\newcommand{\Tft}[1]{\mathscr{#1}^{\operatorname{ft}}}
\newcommand{\To}[1]{\mathscr{#1}^{0}}
\newcommand{\bTft}[1]{\overline{\mathscr{#1}}}
\newcommand{\Diag}{\operatorname{Diag}}
\newcommand{\Aniso}{\operatorname{an}}
\newcommand{\Cl}{\operatorname{Cl}}
\newcommand{\PGL}{\operatorname{PGL}}
\newcommand{\Dyn}{\operatorname{Dyn}}
\newcommand{\Out}{\operatorname{Out}}
\newcommand{\cris}{\operatorname{cris}}
\newcommand{\modp}{\bF}
\newcommand{\cts}{\operatorname{cts}}
\newcommand{\Grp}{\operatorname{Grp}}
\newcommand{\End}{\operatorname{End}}
\newcommand{\SL}{\operatorname{SL}}
\newcommand{\SO}{\operatorname{SO}}
\newcommand{\GO}{\operatorname{GO}}
\newcommand{\Disc}{\operatorname{Disc}}
\newcommand{\cO}{\mathcal{O}}
\newcommand{\cA}{\mathcal{A}}
\newcommand{\clR}{\mathcal{R}}
\newcommand{\cX}{\mathcal{X}}
\newcommand{\fH}{\mathfrak{H}}
\newcommand{\sC}{\mathfrak{C}}
\newcommand{\cG}{\mathcal{G}}
\newcommand{\cS}{\mathcal{S}}
\newcommand{\cT}{\mathcal{T}}
\newcommand{\Fix}{\operatorname{Fix}}
\newcommand{\cB}{\mathcal{B}}
\newcommand{\Lieg}{\mathfrak{g}}
\newcommand{\Lieh}{\mathfrak{h}}
\newcommand{\Mat}{\operatorname{Mat}}
\newcommand{\Gal}{\operatorname{Gal}}
\newcommand{\Art}{\operatorname{Art}}
\newcommand{\BHT}{\mathbb{B}_{\operatorname{HT}}}
\newcommand{\DXT}{{\hatI}^{X(T)}}
\newcommand{\HTGC}{\prod_{\sigma:K_m\hookrightarrow \bC}X_*(G_{\bC})}
\newcommand{\HTGCa}{\prod_{\tilde\tau\in S}X_*(G_{\bC})}
\newcommand{\HT}{{\operatorname{HT}}}
\newcommand{\hT}{{\mathcal{HT}}}
\newcommand{\dR}{{\operatorname{dR}}}
\newcommand{\hatbI}{\widehat{\bar I}_K}
\newcommand{\hatI}{\widehat{I}_K}
\newcommand{\symt}{\operatorname{sym}^2}
\newcommand{\altt}{\operatorname{alt}^2}
\newcommand{\socle}{\operatorname{soc}}
\newcommand{\wh}{\widehat}
\newcommand{\wt}{\widetilde}
\newcommand{\lsup}[2]{{^{#1}\!#2}}
\newcommand{\lsupp}[2]{{^{#1}\!#2}_{\bFp}}
\newcommand{\lsupw}[2]{{^{\mathfrak{#1}}\!#2}}
\newcommand\nlcup{%
  \mathrel{\ooalign{\hss$\cup$\hss\cr%
  \kern0.7ex\raise0.6ex\hbox{\scalebox{0.4}{$\diamond$}}}}}
\newcommand{\matt}[9]{
\left(
\begin{matrix}
#1 & #2 & #3 \\
#4 & #5 & #6 \\
#7 & #8 & #9
\end{matrix}
\right)
}
\newcommand{\bigequiv}[2]{
\begin{tabular}[t]{ccc}
$\left\{\text{
\begin{tabular}{p{.4\linewidth}}
#1
\end{tabular}}
\right\}$
& $\xrightarrow{\cong}$ &
$\left\{\text{
\begin{tabular}{p{.4\linewidth}}
#2 
\end{tabular}}
\right\}$
\end{tabular}
}
\newcommand{\bigequivlabel}[3]{
\begin{tabular}[t]{ccc}
$\left\{\text{
\begin{tabular}{p{.4\linewidth}}
#1
\end{tabular}}
\right\}$
& $\xrightarrow{#3}$ &
$\left\{\text{
\begin{tabular}{p{.4\linewidth}}
#2 
\end{tabular}}
\right\}$
\end{tabular}
}

\newcommand{\bigmap}[3]{
\begin{tabular}[t]{ccc}
$\left\{\text{
\begin{tabular}{p{.4\linewidth}}
#1
\end{tabular}}
\right\}$
& $#2$ &
 {#3}
\end{tabular}
}

\newcommand{\bigbigmap}[8]{
\begin{tabular}[t]{ccc}
#1
& $#5$ &
#2
\\
\rotatebox{270}{$#7$}
&&
\rotatebox{270}{$#8$}
\\
#3
& $#6$ &
#4
\end{tabular}
}

\newcommand{\cellwrap}[2]{
\left\{\text{
\begin{tabular}{p{#1\linewidth}}
#2
\end{tabular}}
\right\}
}
\newcommand{\cellwrapsupset}[3]{
$\left\{\text{
\begin{tabular}{p{.07\linewidth}|p{.25\linewidth}}
#1 & #2
\end{tabular}}
\right\}_{\text{#3}}$
}

\newcommand{\cellwrapsup}[2]{
$\left\{\text{
\begin{tabular}{p{.34\linewidth}}
#1
\end{tabular}}
\right\}_{\text{#2}}$
}

\author{Lin, Zhongyipan}
\date{\today}
\title{A Deligne-Lusztig type correspondence for tame $p$-adic groups}

\begin{abstract}
We establish a ``matrix simultaneous diagonalization theorem'' for disconnected reductive groups which relaxes both the semisimplicity condition and the commutativity condition.
As an application, we prove the following basic results concerning mod $p$ Langlands parameters for quasi-split tame groups $G$ over a $p$-adic field $F$:
\begin{itemize}
\item All semisimple $L$-parameters $\operatorname{Gal}_F\to {^ L\!G}({\bar{\mathbb{F}}_p})$ factor through the $L$-group of a maximal $F$-torus of $G$;
\item All semisimple mod $p$ $L$-parameters admit a de Rham lift of regular $p$-adic Hodge type;
\item A version of tame inertial local Langlands correspondnece; and
\item A group-theoretic description of irreducible components of the reduced Emerton-Gee stacks away from Steinberg parts.
\end{itemize}
We also propose generalizations of the explicit recipe for Serre weights (after Herzig) and the geometric Breuil-M\'ezard for tame groups.
\end{abstract}

\maketitle

\tableofcontents

\newpage

\mychapter{Introduction}

\etocsettocdepth{2}
\localtableofcontents

\vspace{3mm}

In this paper, we present (without proving) novel extensions of both the explicit Serre weight conjecture and the Breuil-M\'ezard conjecture.
Prior to stating these conjectures, we establish foundational results concerning mod $p$ Langlands parameters for
tame groups, such as the classification of elliptic Langlands parameters,
the construction of their de Rham lifts,
and a version of tame inertial local Langlands correspondence.
The conjectures we formulate interpolate the respective conjectures for ramified general linear groups
$\Res_{F/\Qp}\GL_n$ (as discussed in \cite{EG23} and \cite{LLHLM22})
and unramified reductive groups (as discussed in \cite{GHS} and \cite{FLH}).

To set the stage, let us
denote by $F$ a finite extension of $\Qp$
with residue field $\kappa_F$.
Write $\Gal_F:=\Gal(F^s/F)$ for the absolute Galois group of $F$,
and $W_F$ for the Weil group of $F$.
Let $G$ be a quasi-split group over $F$
which splits over a tame extension $L$ of $\Qp$
(so $F$ is also a tame extension of $\Qp$).
Write $\lsup LG=\wh G\rtimes \Gal(L/F)$
for the $L$-group of $G$.

\mysection{Motivation}
In \cite{Se87}, Serre formulated a precise conjecture predicting the minimal weight
of mod $p$ Galois representations arising from modular eigenforms. He also posed a question regarding the potential connection between the weight recipe and a ``mod $p$ Langlands philosophy,'' as well as whether the weight recipe generalizes to encompass general reductive groups (Question 3.4, loc. cit.).

Subsequently, extensive research has been conducted on the Serre weight conjectures, with references available in the introduction of \cite{GHS}. 
The modern version of Serre weight conjecture seeks to classify congruences of Hecke eigensystems
in the cohomology of locally symmetric space associated to a reductive group $G$, with coefficients
in local systems induced by different weights of $G$.

Since the formulation of the Serre weight conjecture for unramified groups in \cite{GHS}, substantial progress has been made in verifying cases of the conjecture beyond $\GL_n$, particularly for groups such as $\GSp_4$ (\cite{Lee23}) and the unramified quasisplit unitary group $U_2$ (\cite{KM22}). These accomplishments motivate further investigation into mod $p$ Langlands parameters and their associated moduli stacks.

In this paper, as well as in the companion paper \cite{L23C}, we undertake a systematic classification of mod $p$ Langlands parameters and explore their lifts in characteristic $0$.

\mysection{Mod $p$ Langlands parameters}
An $L$-parameter $\bar\rho:\Gal_F\to\lsup LG(\bFp)$
is said to be {\it parabolic}
if it factors through a parabolic subgroup
$\lsup LP(\bFp)$,
and is said to be {\it elliptic}
if otherwise.
Here $\lsup LP:=\wh P\rtimes \Gal(L/F)\subset \lsup LG$
and $\wh P$ is a $\Gal(L/F)$-stable parabolic subgroup of $\wh G$.
An elliptic $L$-parameter
is {\it semisimple}
in the sense that its image
is a completely reducible subgroup of $\lsup LG(\bFp)$.

The first main theorem we prove is the following
characterization of elliptic mod $p$ $L$-parameters.

\begin{thm} (Theorem \ref{thm:LS})
\label{thm:B}
If $\bar\rho:\Gal_F\to\lsup LG(\bFp)$
is semisimple, then there exists a maximally unramified
maximal $F$-torus $S$ of $G$
such that $\rho$ admits a factorization
$$
\Gal_F \xrightarrow{\bar\rho_S} \lsup LS(\bFp)\xrightarrow{\lsup Lj}\lsup LG(\bFp)
$$
where $\lsup Lj$ is the mod $p$ Langlands-Shelstad $L$-embedding.
Moreover, if $\rho$ is elliptic, then $S$ is an elliptic torus.
\end{thm}

The immediate consequence of the mod $p$ Langlands-Shelstad
factorization theorem is the existence of de Rham lifts
$\rho:\Gal_F\to \lsup LG(\bZp)$.
Indeed, to construct a de Rham lift of $\bar\rho$,
it suffices to construct a lift of $\bar\rho_S$,
which can be done via the $p$-adic Local Langlands Correspondence
for algebraic tori (see \cite{Ch20}).
In Section \ref{subsec:Hodge-LLC},
we discuss the $p$-adic Hodge-theoretic refinement of the LLC
for algebraic tori and prove the following.

\begin{thm} (Theorem \ref{thm:dR-lift})
\label{thm:C}
If $\bar\rho:\Gal_F\to\lsup LG(\bFp)$
is semisimple, then there exists a potentially crystalline
lift of $\bar\rho$ of regular Hodge type.
\end{thm}

If $G$ is a ramified group, then 
there is no semistable or crystalline $L$-parameters;
therefore potentially crystalline lifts are
the best we can hope for.
Theorem \ref{thm:B} plays a pivotal role in \cite{L23}, where it is used as a crucial input for establishing the Noetherian formal algebraicity of the Emerton-Gee stack $\cX_{\lsup LG}$
which is the foundation of many recent developments such as
\cite{LLHLM23} and \cite{FLH}.

We delve into the theory of parabolic/non-semisimple Langlands parameters in the companion paper \cite{L23C}.
For classical groups, a parabolic mod $p$ $L$-parameter is an iterated Heisenberg-type extension
of elliptic $L$-paramaters.
Based on Theorem \ref{thm:C}, we reduce the general existence of de Rham lifts
to a question about the dimension of certain closed substacks of the reduced Emerton-Gee stacks,
and answer this question affirmatively for unitary groups using
the geometry of Grassmannian manifolds in loc. cit..
So, for tamely ramified unitary groups, Theorem \ref{thm:C} holds for any $L$-parameter $\bar\rho$,
not just semisimple ones (see Theorem 5, loc.cit.).

\mysection{Parahoric Serre weights and the qualitative Breuil-M\'ezard conjecture}~
In this subsection, we briefly discuss implications of the existence of de Rham lifts.

\vspace{3mm}

\noindent
{\bf The Emerton-Gee stacks}
The Emerton-Gee stacks are a version of moduli stacks of $L$-parameters
that enables us to apply the powerful machinery of geometric representation theory
to the study of Galois deformation rings. 
They are first constructed for $\GL_n$
in \cite{EG23}, and then generalized to
general tame groups in \cite{L23}.

\vspace{3mm}

\noindent
{\bf The $\GL_n$-case}
In \cite{EG23}, the irreducible components of the reduced Emerton-Gee stack
$\cX_{\GL_n, \red}$ are shown to be in bijection
with isomorphism classes of irreducible $\bFp$-representations of the finite group $\Res_{F/\Qp}\GL(\bF_p)$.
An irreducible $\bFp$-representation of the finite group $\Res_{F/\Qp}\GL(\bF_p)$
can be inflated to an irreducible $\bFp$-representation of the compact group $\Res_{F/\Qp}\GL(\bZ_p)$,
which is a superspecial parahoric subgroup of $\GL(F)$.
The subgroup $\Res_{F/\Qp}\GL(\bZ_p)$ is also a maximally bounded subgroup of $\GL(F)$.

\vspace{3mm}
\noindent
{\bf The tame group case}
In general, a parahoric subgroup of $G(F)$ is not necessarily a maximally bounded subgroup
of $G(F)$.
We fix a superspecial parahoric $\cG^\circ$ of $G(F)$ and let $\cG\supset \cG^\circ$ be the maximally bounded subgroup containing $\cG^\circ$.
Write $\uG$ for the reductive quotient of $\cG^\circ$,
which is a connected reductive group over the finite field $\kappa_F$.
We call isomorphism classes of irreducible $\bFp$-representations of $\cG^\circ$ {\it parahoric Serre weights}.
A {\it Serre weight} is defined to be an isomorphism class of pairs $(\sigma, \wt \sigma)$
where $\wt \sigma$ is an irreducible $\bFp$-representations of $\cG$
and $\sigma\subset \wt \sigma$ is an irreducible $\bFp$-representations of $\cG^\circ$.
By abuse of notation, we write $\wt \sigma$ for the pair $(\sigma, \wt \sigma)$.

Conjecturally, for each Serre weight $\wt\sigma$, we can attach to it a finite union of irreducible components of $\cX_{\lsup LG, \red}$ which we denote by
$\cC(\wt \sigma)$, such that
$$\cX_{\lsup LG, \red} = \cup_{\wt\sigma} \cC(\wt\sigma)$$
and $\cC(\wt\sigma)\cap \cC(\wt \sigma')$ is nowhere dense in $\cX_{\lsup LG, \red}$
if $\wt\sigma\not\cong \wt \sigma'$.
Moreover, if $\uG$ has simply-connected derived subgroup,
we cojecture that each $\cC(\wt\sigma)$ is an irreducible component of $\cX_{\lsup LG, \red}$.
This is the natural generalization of the so-called {\it qualitative Breuil-M\'ezard conjecture},
see \cite[Section 8.1]{EG23}.
The content of the qualitative Breuil-M\'ezard conjecture is that
$\cX_{\lsup LG, \red}$ is equidimensional of the same dimension as the special fiber potentially semistable/crystalline
deformation rings and its irreducible components admit a group-theoretic parameterization.
(In the literature, the qualitative Breuil-M\'ezard conjecture sometimes refers to a stronger statement, see \cite[Theorem 1.4.5]{LLHLM23}, which further claims that a certain transition matrix is upper triangular. See Subsection \ref{subsec:BM} below for details.) 
If $\sigma$ is a parahoric Serre weight, 
we can define a closed substack
$\cC(\sigma)=\cup_{(\sigma,\wt \sigma)} \cC(\wt \sigma)$.
We will postpone the construction of $\sigma\mapsto \cC(\sigma)$ for {\it regular} parahoric Serre weights
to Subsection \ref{subsec:BM} because we haven't introduced the necessary notations yet.

For classical groups, the qualitative Breuil-M\'ezard conjecture
follows from the existence of de Rham lifts of regular Hodge type, see \cite{L23B}.
In loc. cit., we prove the qualitive Breuil-M\'ezard conjecture for unitary groups;
in particular, for even unitary groups, the irreducible components of $\cX_{\lsup LU, \red}$
are in bijection with Serre weights (which, in this particular case, coincide with parahoric Serre weights).
We also note that the association $\sigma\mapsto \cC(\sigma)$
is constructed for all parahoric Serre weights in loc. cit., not just the regular ones.

\mysection{Simultaneous diagonalization of matrices}~
Before we proceed to formulate the conjectures, we digress and explain the proof of Theorem \ref{thm:B}.

\vspace{3mm}
\noindent
{\bf Characteristic $p$ versus characteristic $0$ coefficients}
The characteristic $0$ coefficient version of Theorem \ref{thm:B}
essentially follows from the work of Borel-Serre
on solvable subgroups of compact Lie groups (\cite{BS53}).
Some variants of the argument can be found in the literature;
see, for example, \cite{Kal19b}.
However, these results make use
of the assumption that the image of semisimple/elliptic $\bar\rho$
consists of semisimple elements,
which is not true in characteristic $p$.
For example, consider $\Ind_{\bQ_{4}}^{\bQ_2}1:\Gal_{\bQ_2}\to \Gal_2(\bF_2)$,
which is semisimple but sends a Frobenius element to $
\begin{bmatrix}
\bar 0 & \bar 1 \\
\bar 1 & \bar 0
\end{bmatrix}
$
(a unipotent element of $\Gal_2(\bF_2)$).

It is also not helpful to consider embeddings $\lsup LG\hookrightarrow\GL_d$,
because the property of being a completely reducible subgroup
is {\it not well-behaved undering embeddings} for disconnected groups.
For example, consider
\begin{align*}
\bZ/3 & \to \GL_2(\bF_3)\\
\bar 1 &\mapsto 
\begin{bmatrix}
\bar 1 & \bar 1\\
\bar 0 & \bar 1
\end{bmatrix}
\end{align*}
$\bZ/3\subset \bZ/3$ is a completely reducible subgroup (because it is the full group),
but its image in $\GL_2(\bF_3)$ is not completely reducible.
As a consequence, even proving $\bar\rho$ is tamely ramified
is not straightforward.

We managed to reduce Theorem \ref{thm:B} to the linear algebra problem (see below).
The author does not know of an elementary proof of Theorem \ref{thm:B}
for low rank ramified unitary groups.

\vspace{3mm}
\noindent
{\bf A linear algebra problem}
A familiar fact from linear algebra is that
two diagonalizable matrices $X$ and $Y$ are simultaneously
diagonalizable if and only if they commute with each other, that is,
$XY=YX$.

We can reinterpret the simultaneous diagonalization theorem
from the perspective of algebraic groups.
First of all, the center of $\GL_d$ does not play a role
in a diagonalization problem, so we can assume $X, Y$
as elements of $\SL_d$.
Matrices $X$ and $Y$ correspond to inner automorphisms
$\Int(X),~\Int(Y)$
of the algebraic group $\SL_d$.
Now, $X$ is a diagonalizable matrix if and only
if there exists a Borel $B\subset \SL_d$
and a maximal torus $T\subset B$
such that the Borel pair $(B, T)$
is $\Int(X)$-stable.
Two matrices $X$ and $Y$ are simultaneously diagonalizable
if and only if there exists a Borel pair $(B, T)$
of $\SL_d$ fixed by both $\Int(X)$ and $\Int(Y)$.
The simultaneous diagonalization theorem
follows from the following 
standard linear algebraic group fact:
two commuting semisimple elements are contained
in a maximal torus (\cite[Corollary 8.6]{St68}).

In this paper, we need a generalized simultaneous diagonalization problem.
Instead of inner automorphisms, we allow for outer automorphisms;
and instead of the commuting relation, we only impose the metacyclic relation;
finally, we relax the semisimplicity assumption.
We prove the following:

\begin{thm}
\label{thm:A}
(Corollary \ref{cor:metacyclic})
Let $H$ be a connected reductive group over an algebraically closed field
and let $\tau, \sigma$ be two automorphisms of $H$.
Assume
\begin{itemize}
\item[(MC)] $\sigma\tau\sigma^{-1}=\tau^q$ for some integer $q$,
\item[(WS)] $\tau$ is a semisimple automorphism and $\langle \sigma, \tau\rangle$ is a pseudo-completely reducible subgroup
of $\Aut(H)$.
\end{itemize}
Then there exists a Borel pair $(B, T)$ of $H$ such that
$T$ is both $\sigma$- and $\tau$-stable
while $B$ is $\tau$-stable.
\end{thm}
Here (MC) stands for metacyclicity and (WS) stands for weak semisimplicity.
We remark that
by forming the semi-direct product
$H':=H\rtimes \langle\sigma, \tau\rangle$,
we can instead 
insist that both $\sigma$ and $\tau$
are inner automorphisms of $H'$
while allowing $H'$ to be a disconnected
reductive group.

\mysection{A tame inertial local Langlands correspondence}

The second application of the factorization theorem (Theorem \ref{thm:B})
is tame inertial Local Langlands Correspondence.
Recall that $\uG$ denotes the reductive quotient
of a superspecial parahoric of $G$.
An {\it inertial Deligne-Lusztig datum}
is a pair $(\uS, \uchi)$
where $\uS$ is a maximal torus of $\uG$
and $\uchi$ is a character $\uS(\kappa_F)\to \bFp^\times$.
Computation shows that
geometric conjugacy classes of inertial Deligne-Lusztig data
and of tame inertial types $I_F\to \lsup LG(\bFp)$
admit the same combinatorial parametrization,
and thus there exists a set-theoretic bijection
between these two.
However, such a bijection depends a priori on certain choices
and its arithmetic significance is not clear.
Using the Langlands-Shelstad factorization,
we give a theoretic explanation of
this bijection
and show it is compatible with (indeed determined by)
the LLC for algebraic tori.

\begin{thm} (Corollary \ref{cor:DL})
\label{thm:D}
There exists a natural bijection
$$
\{\text{Inertial Deligne-Lusztig data $(\uS, \uchi)$}\}_{/\uG(\bFp)}
\xrightarrow{\cong}
\{\text{Tame inertial types $I_F\to \lsup LG(\bFp)$}\}_{/\wh G(\bFp)}
$$
\end{thm}

The theorem above is a natural extension
of the Deligne-Lusztg duality for finite groups of Lie type
to quasi-split tame $p$-adic groups,
and we call it {\it the Deligne-Lusztig correspondence}.
It is the starting point of 
the generalization of
the Gee-Herzig-Savitt recipe for Serre weights
for unramified groups to tame groups;
we will elaborate on this topic in the next subsection.

\mysection{Explicit Serre weight conjectures, after Herzig}

In this subsection, we present a generalization of Herzig's explicit Serre weight recipe.

By the non-abelian Shapiro's lemma (\cite[Proposition 1]{L23}),
working with the $F$-group $G$ is equivalent to working with the $\bQ_p$-group
$\Res_{F/\Qp}G$.
Assume $G$ is a quasi-split tame group over $\Qp$ in the rest of the introduction for simplicity
of notation, without loss of generality.
Recall that $\cG^\circ$ is a superspecial parahoric of $G$, and
$\uG$ is the reductive quotient of $\cG^\circ$.

For technical simplicity,
we assume both $G$ and $\uG$ admit
a local twisting element
and $\uG$ has a simply-connected derived subgroup.
Since $G$ is the generic fiber of
$\cG^\circ$,
we denote by $\eta_{\Qp}$
a local twisting element for $G$;
and since $\uG$ is the reductive part of
the special fiber of $\cG^\circ$,
we denote by $\eta_{\bF_p}$
a local twisting element for $\uG$.
By restricting to the maximal unramified subtorus
of a maximal torus of $G$,
we can regard $\eta_{\Qp}$
as an element of the character lattice of $\uG$.

For each tame inertial $L$-parameter
$\tau:I_{\Qp}\to \lsup LG(\bFp)$,
we define
$$
W^?(\tau):= \clR(\JH(\bar V(\DL^{-1}(\tau)) \otimes W(w_0(\eta_{\bF_p}-\eta_{\Qp})))).
$$
Here $\clR$ is Herzig's involution operator
and $\bar V$ is the reduction mod $p$ of
the Deligne-Lusztig induction functor
(with a sign modification).
See Subsection \ref{subsec:Serre}
for other unfamiliar notations and clarifications.

\vspace{3mm}

\noindent {\bf Remark}
\begin{itemize}
\item
Our definition is inspired by the work of \cite{LLHLM22}
for $\Res_{E/\Qp}\GL_d$.

\item 
When $G$ is unramified,
$\eta_{\Qp}-\eta_{\bF_p}=0$
and thus $W^?(\tau)$
recovers the definition found in
\cite[Section 9]{GHS}.

\item
$\DL^{-1}(\tau)$ is only well-defined up to geometric conjugacy while the input for the Deligne-Lusztig induction
needs to be well-defined up to rational conjugacy.
As a consequence, the composition
$\bar V\circ \DL^{-1}$ is ambiguous
(it is a multi-valued map).
\end{itemize}

\vspace{3mm}

\noindent {\bf $L$-packets}
We want to elaborate on the last bullet point of the remarks above.
To resolve the ambiguity for $\bar V\circ \DL^{-1}$,
\cite{GHS} considers only the {\it maximally split}
rational conjugacy class of $\DL^{-1}$,
which is unique after imposing certain technical assumptions on $G$,
indeed,
for sufficiently generic tame types $\tau$,
maximal splitness of $\DL^{-1}(\tau)$ is automatic (\cite[Proposition 6.20, Lemma 6.24]{Her09}), and therefore so is the uniqueness of the rational conjugacy class (under certain technical assumptions).

From the perspective of the Local Langlands Correspondence,
an $L$-parameter should correspond to an $L$-packet
of admissible representations,
rather than a single admissible representation.
When $G$ is an unramified $p$-adic group,
the maximally split rational conjugacy class of $\DL^{-1}(\tau)$
is expected to correspond to the {\it generic constituent}
of the $L$-packet.
For more general groups, there can be multiple generic constituents
in $L$-packets.
In generic situations, we expect the various rational conjugacy classes
of $\DL^{-1}(\tau)$ to correspond to
the various generic constituents in the corresponding $L$-packet.

The reader may want to interpret $W^?(\tau)$
not as a set, but rather
as a packet of sets
where $\DL^{-1}(\tau)$
ranges over all rational conjugacy classes.

\mysection{The geometric Breuil-M\'ezard conjecture: the potentially crystalline case}~
\label{subsec:BM}

We refer to \cite[Section 2.4, 2.5, 2.6]{L23B}
for the definition of the potentially crystalline stacks $\cX^{\crys,\lambda, \tau}_{\lsup LG}$
of Hodge type $\lambda$ and tame inertial type $\tau$.
Roughly speaking, a Hodge type is a conjugacy class of cocharacters $\lambda$ of $\wh G$.
Since cocharacters of $\wh G$ are identified with characters of $G$,
we regard $\lambda$ as an element of the character lattice of $G$.

Denote by $V(\lambda-\eta_{\bQ_p})$ the (restriction to $\cG^\circ$ of the)
irreducible algebraic representation of highest weight $\lambda-\eta_{\bQ_p}$.
For a finite-dimensional representation $R$ of $\cG^\circ$ over a finite extension of $\bQ_p$,
write $\bar R$ for the semisimplification of the reduction mod $p$ of a $\cG^\circ$-invariant lattice of $R$.

The natural generalization of \cite[Theorem 1.4.5(1)]{LLHLM23} is the following conjecture.

\begin{conj}
If $\lambda$ is regular dominant and $\tau$ is a sufficiently generic tame inertial type,
then
$$\cX^{\crys, \lambda, \tau}_{\lsup LG, \red}=\underset{\sigma\in \JH(\overline{\bar V(\DL^{-1}(\tau))\otimes V(\lambda-\eta_{\bQ_p})}}{\bigcup}~\cC_\sigma$$
\end{conj}

Write $[\cX^{\crys, \lambda, \tau}_{\lsup LG,\bF_p}]$ for the cycle class of 
$\cX^{\crys, \lambda, \tau}_{\lsup LG,\bF_p}$ in the Chow group of $\cX_{\lsup LG, \red}$.
The conjecture above has the following refinement generalizing \cite[Conjecture 8.2.2]{EG23}.

\begin{conj}
For each parahoric Serre weight $\sigma$, there exists an effective top-dimensional cycle
$\cZ_\sigma$ on $\cX_{\lsup LG, \red}$ such that
$$
[\cX^{\crys, \lambda, \tau}_{\lsup LG,\bF_p}] = \sum_{\sigma} [\overline{\bar V(\DL^{-1}(\tau))\otimes V(\lambda-\eta_{\bQ_p})}:\sigma] \cZ_\sigma
$$
for all regular $\lambda$ and tame inertial types $\tau$.

\end{conj}

We don't have evidence for the conjecture except that it agrees with the ramified general linear case
$\Res_{F/\Qp}\GL_n$ (\cite{EG23}) and the split group case (\cite{FLH}).
However, since the method of \cite{FLH} is purely local and group-theoretical,
and it seems plausible that the arguments of \cite{FLH} generalize to general tame groups
once the corresponding players are correctly defined.

Finally, we explain the construction of $\sigma\mapsto \cC(\sigma)$.
There exist natural bijections between the following objects:
\begin{itemize}
\item[(1)]
equivalence classes of regular parahoric Serre weights,
\item[(2)]
geometric conjugacy classes of regular based inertial Deligne-Lusztig data of niveau $1$,
\item[(3)]
$\wh G(\bFp)$-conjugacy classes
of regular inertial $L$-parameters $I_F\to \lsup LB(\bFp)$
where $\lsup LB=\wh B\rtimes \Gal(L/F)$ is a Borel of $\lsup LG$.
\end{itemize}
The bijection (1) $\Leftrightarrow$ (2)
follows from Proposition \ref{prop:parahoric-based};
and the bijection (2) $\Leftrightarrow$ (3) follows from
Theorem \ref{thm:D}.
So, we identify (1), (2) and (3) implicitly in the rest of this subsection.
A regular parahoric Serre weight $\sigma$ of niveau $1$ admits Herzig's presentation
$(1, \mu_\sigma)$ (see \ref{subsec:herzig-pre}).

\begin{dfn}
$\cC(\sigma)$ is defined to be the closure
of $\bFp$-points of $\cX_{\lsup LG, \red}$
corresponding to $L$-parameters $\bar\rho:\Gal_{\Qp}\to \lsup LG(\bFp)$
that factor through a unique Borel $\lsup LB$
such that $\bar\rho|_{I_{\Qp}}:I_{\Qp}\to \lsup LB(\bFp)$
has Herzig's presentation $(1, -w_0(\mu_{\sigma}) -\eta_{\bF_p})$.
\end{dfn}

The reader can verify that our definition is consistent with \cite[Definition 5.5.1]{EG23}
for $\Res_{F/\Qp}\GL_n$.

\mysection{Future directions}

Before we finish the introduction,
we raise the following natural questions that are not addressed in this paper.

\vspace{3mm}

\noindent {\bf Question A}
How about wildly ramified $p$-adic groups?

\vspace{3mm}
If $p\ge 5$, then any connected reductive group over a $p$-adic field $F$ is isogenous to a product of groups of the form $\Res_{K/F}G$
where $G$ is tame over $K$.
The Weil restricted case follows from the tame case by
non-abelian Shapiro's lemma (see, for example, \cite[Section 7.2]{L23}).
So it remains to consider the $p=2, 3$-case.
It is probably approachable by a case-by-case analysis,
and we expect the most complicated case to be related to the triality.

\vspace{3mm}

\noindent {\bf Question B}
How about the mod $l$ ($l\ne p$) analogue?

\vspace{3mm}
To study mod $l$ Langlands parameters,
we need a further generalization of Theorem \ref{thm:A},
which allows $\tau$ to be a non-semisimple automorphism.
We expect generalizations of Theorem B-D to hold
under some mild assumptions on the prime $p$.

\vspace{3mm}

\noindent {\bf Question C}
How about crystalline lifts of semisimple mod $p$ Langlands parameters?
Can we formulate a conjecture relating crystalline lifts
and Serre weights?

\vspace{3mm}
When $G$ is a split group, crystalline lifts of semisimple $L$-parameters are constructed in \cite{L22}.
When $G$ is ramified, crystalline $L$-parameters do not exist.
So it is an interesting question whether crystalline lifts exist
when $G$ is unramified and non-split.

The reason that crystalline $L$-parameters do not exist for ramified groups is probably simply because the na\"ive definition of crystallinity is incorrect.
Similar issues arise when people study the Satake correspondence for 
ramified groups,
in which case unramified $L$-parameters do not exist;
nonetheless, we can still define ``spherical'' $L$-parameters
for ramified groups (see \cite{Zhu15}).
Similarly, we say a $p$-adic $L$-parameter is ``crystalline''
if it is potentially crystalline and its corresponding Weil-Deligne representation
is ``spherical'' in the sense of \cite[Definition 6.3]{Zhu15}.
It seems an interesting question to explore analogues of the crystalline lift aspects of the Serre weight conjectures.

\vspace{3mm}

\mysubsection{Acknowledgements}
The author thanks Bao Le Hung for many inspiring conversations.
The author is grateful to Tasho Kaletha for helpful correspondence
and sharing the manuscript of his book.
The author thanks Florian Herzig for his interest in this project and for his correspondence.
Many ideas in this paper originate from
the author's thesis research and the author
thanks his PhD supervisor David Savitt for introducing the subject to
the author.
Finally, the reader will find the influence of the work of Emerton-Gee, Gee-Herzig-Savitt and Le-Le Hung-Levin-Morra
to the author.

\mysubsection{Notation and conventions}
Write $\breve{F}$ for the maximal unramified extension of $F$
inside a fixed separable closure $F^s$.
For each finite extension $E$ of $F$, denote by
$\kappa_E$ the residue field of $E$,
and denote by $\cO_E$ the ring of integers of $E$.

We denote $-\otimes_{\bZ}-$ by $-\otimes -$,
and denote $\Hom_{\operatorname{\Grp}}(-,-)$ by $\Hom(-,-)$.

Write $\QZp$ for the prime-to-$p$ divisible subgroup of $\bQ/\bZ$.
Note that $\QZp$ is isomorphic to $\bFp^\times$
as an abelian group.

\mychapter{Metacyclic actions on reductive groups}

\etocsettocdepth{2}
\localtableofcontents

\mysection{Extending pseudo-parabolics of disconnected groups}

\mysubsection{Dynamic methods}
\label{dyn-fil}
We recall the definitions in \cite[2.3]{L22}.
Let $H$ be an algebraic group over a ring $k$.
Let $f: \Gm \to H$ be a \textit{scheme morphism}.
Define the following functor on the category of $k$-algebras
$$P_H(f)(A)=\{g\in H(A)|\lim_{t\to 0}f(t)g f(t)^{-1}\text{~exists.}\}$$
where $A$ is a general $k$-algebra. 
We call $f$ a \textit{fake cocharacter}.
Here ``a limit exists'' means
the scheme morphism  $\Gm \to H$, defined by $t\mapsto f(t)gf(t)^{-1}$,
extends to a scheme morphism $\bA^1\to H$.
Note that $P_H(f)$ is not representable in general.
We define similarly $U_H(f)$ by setting $A\mapsto \{g\in H(A)|\lim_{t\to 0}f(t)g f(t)^{-1}=1\}$, and 
$Z_H(f)$ by setting $A\mapsto \{g\in H(A)|f g =g f\}$.

\mysubsection{Definition}
Let $H$ be an algebraic group over an algebraically closed field.
A fake cocharacter $f:\Gm \to H$ is said to be {\it relevant}
if the functor $P_H(f)$ is representable
    by a linear algebraic group.
In particular, all maximal tori of $P_H(f)$
    are conjugate to each other, by an element of $P_H(f)$.

    \vspace{3mm}

The reason we introduce the notion of fake cocharacters is because of the following
powerful tool:

\mysubsection{Lemma}
\label{lem:dynamic}
Let $M \hookrightarrow G$ be two possibly disconnected algebraic groups over an algebraically closed field $k$.
Write $M^\circ$ for the neutral component of $M$.
Let $\lambda,\mu:\Gm\to M^\circ$
be two relevant cocharacters of $M^\circ$.
Assume $P_{M^\circ}(\lambda) = P_{M^\circ}(\mu)=:P$.

\noindent(1) There exists a relevant cocharacter $f:\Gm\to M^\circ$ such that
$P_{G}(f) = P_G(\mu\lambda)$ as a functor.
In particular, the fake cocharacter
$\mu\lambda:\Gm\to M^\circ$ is relevant.

We have $P_{M}(\lambda) \cap P_M(\mu) \subset P_M(\mu\lambda)=P_M(f)$.

Moreover, if $M^\circ$ is a reductive group, then
$P_{M^\circ}(\lambda) = P_{M^\circ}(\mu\lambda)$.

\noindent(2)
The limit
$$
\lim _{t\to 0} \lambda(t)\mu(t)\lambda(t)^{-1}\mu(t)^{-1}
$$
exists in the sense of subsection \ref{dyn-fil}, and lies in $P$.

\noindent(3)
Let $u$ be an element of $P$.
The limit
$$
\lim _{t\to 0} \lambda(t)u \mu(t) u^{-1}\lambda(t)^{-1}\mu(t)^{-1}
$$
exists in the sense of subsection \ref{dyn-fil} and lies in $P$.

\noindent(4)
Now assume $\lambda$ is a product of cocharacters $\lambda_1$, \ldots, $\lambda_s$
such that $P_{M^\circ}(\lambda_i)=P$ for all $i$.
(1), (2) and (3)
remain true (for example, $P_G(\mu\lambda_1\dots\lambda_s)=P_G(f)$ for some cocharacter $f:\Gm\to M^\circ$).

\begin{proof}
Since all maximal tori in $P$ are conjugate to each other and the image of a cocharacter is contained in a maximal torus,
there exists an element $x\in P$
such that
$(x\lambda x^{-1})\mu = \mu (x\lambda x^{-1})$.
Write $\xi$ for $x \lambda x^{-1}$.

(1)
We have 
\begin{align*}
\lim_{t\to 0} \mu(t) \lambda(t) g \lambda(t)^{-1}\mu(t)^{-1} & = 
\lim_{t\to 0} \mu(t) x^{-1}\xi(t) x g  x^{-1} \xi(t)^{-1} x \mu(t)^{-1} \\
&= \lim_{t\to 0} (\mu(t) x^{-1} \mu(t)^{-1})\cdot(\mu(t)\xi(t) x g  x^{-1} \xi(t)^{-1} \mu(t)^{-1})\cdot (\mu(t) x \mu(t)^{-1}) \\
&= \lim_{t\to 0} \mu(t) x^{-1} \mu(t)^{-1} \cdot \lim_{t\to 0} \mu(t)\xi(t) x g  x^{-1} \xi(t)^{-1} \mu(t)^{-1} \cdot \lim_{t\to 0} \mu(t) x \mu(t)^{-1}
\end{align*}
Since $x\in P$, $\lim_{t\to 0} \mu(t) x \mu(t)^{-1}$ exists.
So we have $P_G(\mu\lambda) = x^{-1}P_G(\mu\xi)x = P_G(x^{-1}\mu\xi x)$.
Note that $f:=x^{-1}\mu\xi x$ is a cocharacter.

It is obvious that $P_{M}(\lambda) \cap P_M(\mu) \subset P_M(\mu\lambda)=P_M(f)$.

Next we consider the ``moreover'' part. Since $\mu\xi =\xi \mu$, we can regard $\mu$ and $\xi$ as elements in a cocharacter lattice $X_*(M^\circ,T)$
where $T$ is a maximal torus containing both $\mu$ and $\xi$.
Since $P_{M^\circ}(\mu) = P_{M^\circ}(\lambda) = P_{M^\circ}(\xi)$, $\mu$ and $\xi$ lie in the same Weyl chamber (the Borel case) or the same wall of Weyl chamber (the non-Borel case).
The cocharacter $\mu\xi$ is the sum of $\mu$ and $\xi$ in the cocharacter lattice $X_*(G,T)$,
and lies in the same (wall of) Weyl chamber.
So $P_{M^\circ}(\mu\xi) = P_{M^\circ}(\mu) = P_{M^\circ}(\lambda)$.
Since $x \in P$, we have
$P_G(\mu\lambda) = x^{-1}P_G(\mu\xi)x = P_G(\mu)=P_G(\lambda)$.

(2)
We have 
\begin{align*}
\lambda(t)\mu(t)\lambda(t)^{-1}\mu(t)^{-1}
&=
x^{-1}\xi(t) x \mu(t) x^{-1} \xi(t)^{-1} x \mu(t)^{-1} \\
&=
x^{-1}\cdot (\xi(t) x \xi(t)^{-1}) \cdot (\xi(t) \mu(t) x^{-1} \mu(t)^{-1}\xi(t)^{-1})
\cdot (\mu(t) x \mu(t)^{-1})
\end{align*}
By (1), $P = P_{M^\circ}(\lambda)\cap P_{M^\circ}(\mu)\subset P_{M^\circ}(\xi\mu)$, and thus the limits
$$\lim_{t\to 0}\xi(t)g \xi(t)^{-1},$$
$$\lim_{t\to 0}\xi(t)\mu(t)g^{-1} \mu(t)^{-1}\xi(t)^{-1}\text{, and}$$
$$\lim_{t\to 0}\mu(t)g \mu(t)^{-1}$$ all exist and lies in $P$.
As a consequence, 
$\lim_{t\to 0}\lambda(t)\mu(t)\lambda(t)^{-1}\mu(t)^{-1}\in P$.

(3)
We have
$$
\lambda(t)u \mu(t) u^{-1}\lambda(t)^{-1}\mu(t)^{-1}
 = 
(\lambda(t)u \mu(t) u^{-1}\lambda(t)^{-1} u\mu(t)^{-1} u^{-1}) (u \mu(t) u^{-1} \mu(t)^{-1}).
$$
So (3) follows from (2).

(4)
The method is the same but notations are more complicated.
The reader can consult \cite[Lemma 2.4]{L22} for a proof.
\end{proof}

\mysubsection{Theorem}(The Parabolic Extension Theorem)
\label{thm:par-ext}
Let $M$ be a possibly disconnected linear algebraic group over an algebraically closed field with reductive neutral component $M^\circ$.
Let $H$ be a disconnected linear algebraic group which is a group extension
$$
1\to M\to H \to \langle \bar \gamma \rangle \to 1,
$$
where $\langle \bar \gamma \rangle$ is a cyclic group $\bZ/d$,
regarded as a constant group scheme.

Let $f:\Gm\to M^\circ$ be a cocharacter
such that
\begin{itemize}
\item[(PE1)] $P_{M}(f)$ is stabilized by $\gamma$ for some $\gamma\in H$ which lifts $\bar \gamma$;
\item[(PE2)] The normalizer of $P_{M}(f)$ in $M$ is $P_{M}(f)$ itself.
\end{itemize}
Then there exists a cocharacter $f':\Gm\to M^\circ$ such that
\begin{itemize}
\item[(1)]
$P_{M^\circ}(f) = P_{M^\circ}(f')$,
\item[(2)]
$P_M(f)\subset P_M(f')$, and
\item[(3)]
$\gamma\in P_H(f')$.
\end{itemize}

\begin{proof}
We have $P_{M^\circ}(\gamma f \gamma^{-1})=\gamma P_{M^\circ}(f) \gamma^{-1} = P_{M^\circ}(f)$
by unravelling the definitions.
Write $F(-)$ for
$\gamma\cdot (-) \cdot \gamma^{-1}$
(conjugation-by-$\gamma$).
Define
$$
\mu:=F^{d-1}(f)F^{d-2}(f)\dots F(f)f.
$$
Say $\gamma^d = u \in M$.
By (PE1), $u$ is in the normalizer of $P_{M}(f)$ in $M$.
By (PE2), $u\in P_{M}(f)$.
Note that
$$
F(\mu) = (u f u^{-1})\mu f^{-1}.
$$
We have
\begin{align*}
\lim_{t\to 0}\mu(t) \gamma  \mu(t)^{-1}
&=
\lim_{t\to 0}\mu(t) \gamma  \mu(t)^{-1} \gamma^{-1} \gamma  \\
&=
\lim_{t\to 0}\mu(t)F(\mu)(t)^{-1} \gamma  \\
&=
\lim_{t\to 0}\mu(t) f(t) \mu(t)^{-1}u f(t)^{-1}u^{-1}  \gamma  \\
&=
\lim_{t\to 0}(\mu(t) f(t) \mu(t)^{-1}f(t)^{-1})(f(t)u f(t)^{-1})u^{-1}  \gamma\\
&=
\lim_{t\to 0}\mu(t) f(t) \mu(t)^{-1}f(t)^{-1}\cdot \lim_{t\to 0}f(t)u f(t)^{-1} \cdot u^{-1}  \gamma\\
\end{align*}
By Lemma \ref{lem:dynamic}, the limit $\lim_{t\to 0}\mu(t) f(t) \mu(t)^{-1}f(t)^{-1}$ exists
and lies in $P_{M^\circ}(f)$.
Since $u\in P_{M}(f)$, the limit 
$\lim_{t\to 0}f(t)u f(t)^{-1}$ exists
and lies in $P_M(f)$.
Therefore the limit $\lim_{t\to 0}\mu(t) \gamma  \mu(t)^{-1}$ exists and lies in $P_M(f)\gamma$.
So $\gamma \in P_H(\mu)$.
Since $P_M(F^i(f)) = P_M(f)$ for all $i\ge 0$,
we have
$P_M(f) = \cap_{i\ge 0}P_M(F^i(f))\subset P_M(\mu)$.
We've seen the fake cocharacter $\mu$ satisfies all the requirements. By the previous lemma, there exists a genuine cocharacter
$f':\Gm \to M^\circ$ such that
$P_{H}(\mu)=P_H(f')$.
\end{proof}

\mysubsection{Corollary}
\label{cor:par-ext}
Let $M$
be a linear algebraic group over an algebraically closed field such that
\begin{itemize}
\item 
the neutral component $M^\circ$ is reductive and
\item the component group $\pi:=M/M^\circ$
is a finite solvable group.
\end{itemize}

Let $\Pi \subset M$ be a solvable subgroup which maps surjectively onto $\pi$.
For each $\Pi$-stable parabolic $P\subset M^\circ$, there exists a cocharacter
$f:\Gm\to M^\circ$ such that
\begin{itemize}
\item $P = P_{M^\circ}(f)$, and
\item $\Pi\subset P_M(f)$.
\end{itemize}

\begin{proof}
Let $\Pi_0\subset \Pi$ be a normal subgroup such that $\Pi/\Pi_0$ is cyclic.
Write $\pi_0$ for the image of $\Pi_0$ in $\Pi$.
Choose $\gamma\in \Pi$ which maps
surjectively onto $\pi/\pi_0$.

Write $M^\circ{\Pi_0}$ for the subgroup of $M$
generated by $M^\circ$ and $\Pi_0$.

By induction, there exists a cocharacter $f_0:\Gm\to M^\circ\Pi_0$ such that
\begin{itemize}
\item $P_{M^\circ}(f_0) = P$, and
\item $\Pi_0 \subset P_{M^\circ{\Pi_0}}(f_0)$.
\end{itemize}
Since $\Pi_0$ maps surjectively onto $(M^\circ{\Pi_0})/M^\circ$, we have $P_{M^\circ{\Pi_0}}(f_0) = P\Pi_0$.
Thus $P_{M^\circ\Pi_0}(f_0)$ is stabilized by $\Pi$; and the assumption (PE1) in the previous theorem is verified.
Suppose $g h\in M^\circ \Pi_0$
lies in the normalizer of $P_{M^\circ\Pi}(f_0)$, and that $g\in M^\circ$ and $h\in \Pi_0$.
It follows that $g$ also stabilizes $P_{M^\circ\Pi}(f_0)\supset P$. So $g\in P$.
Thus the normalizer of $P_{M^\circ \Pi_0}(f_0)$ in $M^\circ \Pi_0$
is contained in $P\Pi_0= P_{M^\circ \Pi_0}(f_0)$; and we've verified the assumption (PE2).
By the theorem above, there exists
a cocharacter $f:\Gm\to M^\circ$ such that
\begin{itemize}
\item $P_{M^\circ}=P$,
\item $P_{M^\circ\Pi_0}(f_0)\subset P_{M^\circ\Pi_0}(f)$,
and 
\item $\gamma \in P_{M}(f)$.
\end{itemize}
The second and the third bullet point
jointly imply $\Pi \subset P_{M}(f)$.
\end{proof}

\mysection{Complete reducibility for disconnected groups}

We generalize Serre's notion of $G$-complete reducibility to disconnected groups.

\mysubsection{Definition}
Let $H$ be a linear algebraic group scheme.
A subgroup of $H$ is said to be \textit{a pseudo-parabolic}
if it is of the form $P_H(f)$
    for some cocharacter $f:\Gm\to H^\circ$.

A subgroup of $H$ is said to be \textit{a pseudo-Levi}
if it is of the form $Z_H(f)$
    for some cocharacter $f:\Gm\to H^\circ$.

%

%

\mysubsection{Remark}
The terminology is misleading.
If we only consider groups $H$ over a characteristic $0$ field,
then the notion of ``pseudo-parabolic''
is strictly stronger than 
the usual notion of parabolic (closed subgroups $P$ such that $G/P$ is a projective variety).

Here is a dumb example: let $H$ be $\GL_2 \rtimes (\{1, \sigma\}\times \{1, \tau\})$
where both $\sigma$ and $\tau$ acts on $\GL_2$ by $(-)\mapsto (-)^{-t}$;
any pseudo-parabolic that contains $\sigma$ must contain $\tau$ as well. So all pseudo-parabolics are conjugate
to either $B$ or $B\rtimes (\{1, \sigma\}\times \{1, \tau\})$
where $B$ is the subgroup of upper-triangular matrices.
However, the group $B \rtimes \{1, \sigma\}$
is a parabolic subgroup in the usual sense.

\mysubsection{Definition}
\label{def:pirr}
Let $H$ be an algebraic group over an algebraically closed field.
A \textit{pseudo-parabolic} $P$ of an algebraic group $H$
is said to be a {\it proper} pseudo-parabolic if $P\cap H^\circ\ne H^\circ$.

An abstract subgroup $\Gamma$ of $H$
is said to be \textit{pseudo-irreducible} in $H$ if $\Gamma$ is not contained in any proper pseudo-parabolic subgroup.

An abstract subgroup $\Gamma$ of $H$
is said to be \textit{pseudo-completely reducible} in $H$ if whenever $\Gamma$ is contained in a proper pseudo-parabolic,
it is also contained in a correponding pseudo-Levi.

\mysubsection{Remark}
The notion of pseudo-irreducibility is well-behaved in the following situation: let $M$ be a subgroup of $H$ and assume $\Gamma\subset M$; if $\Gamma$ is pseudo-irreducible in $H$ then $\Gamma$ is also pseudo-irreducible in $M$.

This is not immediate from the definitions.
Suppose $H = \GL_{2n} \rtimes \{1, \sigma\}$.
and let $M = \SO_{2n}\rtimes \{1, \sigma\}\subset H$. 
The above paragraph claims all $\sigma$-stable proper parabolics of
$\SO_{2n}$ extend to a parabolic of $\GL_{2n}$
which is $\sigma$-stable.
The Parabolic Extension Theorem \ref{thm:par-ext}
is where we comfirmed the above claim in full generality.

\mysubsection{Definition}
Let $M$ be a linear algebraic group over 
an algebraically closed field.
A pseudo-parabolic $P$ of $M$ is said to be {\it big} if $P$ maps surjectively onto
$M/M^\circ$.

\mysubsection{Lemma}
\label{lem:big-par}
Let $M$ be a disconnected linear algebraic group
over an algebraically closed field,
with reductive neutral component.
Assume $M/M^\circ$ is solvable.
A subgroup $P$ of $M$ is a big pseudo-parabolic if and only if
\begin{itemize}
\item $P\cap M^\circ$ is a parabolic of $M$, and
\item $P$ maps surjectively onto $M/M^\circ$.
\end{itemize}

\begin{proof}
By induction, we can assume $M/M^\circ$ is cyclic.
Choose an element $\gamma\in P$ which maps to a generator of
$M/M^\circ$.
By Corollary \ref{cor:par-ext},
there exists a cocharacter $f:\Gm\to M^\circ$ such that
\begin{itemize}
\item $P\cap M^\circ = P_{M^\circ}(f)$, and
\item $\gamma\in P_{M}(f)$.
\end{itemize}
So $P \subset P_{M}(f)$.
Since $P$ is big, we must have $P = P_M(f)$.
\end{proof}

\mysubsection{Lemma}
\label{lem:big-ps}
Let $M$ be a linear algebraic group over 
an algebraically closed field, with reductive neutral component.
Assume $M/M^\circ$ is solvable.
Two big pseudo-parabolics $P$, $Q$ of $M$ are the same
if and only if $P^\circ = Q^\circ$.

\begin{proof}
Suppose $P^\circ=Q^\circ$.
By induction, we can assume $M/M^\circ$ is cyclic.
Let $g\in P$ and $h\in Q$
be elements in the same component of $M$ which generate all components.
So $h^{-1}g\in N_{M}(P^\circ) \cap M^\circ = N_{M^\circ}(P^\circ) = P^\circ$.
Thus $P=Q$.
\end{proof}

\mysubsection{Remark}
Both of the lemmas above fail without the bigness assumption.
Let $\omega$ be a primitive cubic root of unity.
Let $H$ be $\GL_2 \rtimes (\{1, \sigma, \sigma^2\}\times \{1, \tau, \tau^2\})$
where both $\sigma$ and $\tau$ acts on $\GL_2$ by 
conjugation by $
\begin{bmatrix}
\omega & \\ & \omega^{-1}
\end{bmatrix}
$ and 
$\frac{1}{2}\begin{bmatrix}
-1& \omega^{-1}-\omega  \\  \omega^{-1}-\omega & -1
\end{bmatrix}=:x$, respectively.
The group $H$ contains no big pseudo-parabolic.
The two pseudo-parabolics
$B\rtimes \langle\sigma\rangle$ and
$B\rtimes \langle x^{-1}\tau\rangle$ are not conjugate to each other.

\vspace{3mm}

At the end of this section, we clarify the relation between
various notions of semisimplicity.

We recall Steinberg's definition (\cite[Section 9]{St68}) of quasi-semisimple automorphisms:

\mysubsection{Definition}
An automorphism $f:M\to M$ of connected linear algebraic groups
is said to be {\it quasi-semisimple} if there exists a Borel $B\subset M$ and a maximal torus $T\subset B$
such that $f(B)=B$ and $f(T)=T$.

A well-known theorem is

\mysubsection{Theorem}(\cite[Theorem 7.5]{St68})
Semisimple automorphisms are quasi-semisimple.

\mysubsection{Definition}
An automorphism $f:M\to M$ of connected linear algebraic groups
is said to be {\it pseudo-semisimple} if 
$\langle f \rangle$ is a pseudo-completely reducible
subgroup in $\Aut(M)$.

\mysubsection{Proposition}
(1) Semisimple automorphisms of reductive groups are pseudo-semisimple.

(2) Pseudo-semisimple automorphisms of reductive groups are quasi-semisimple.

(3) A quasi-semisimple automorphism is semisimple
if and only if its class in $\Aut(G)/\Int(G)$
has order prime to the characteristic (exponent) of the field.

\begin{proof}
(1) This is \cite[Proposition 2.2]{L22}.

(2)
Let $f:M\to M$ be an automorphism of reductive groups.
By replacing $M$ by $M/Z_M(M)$, we can assume $M$ is centerless
and $\Aut(M) = M \rtimes \Out(M)$ where $\Out(M)$ is a finite group.

Suppose $f:M\to M$ is pseudo-semisimple.
By \cite[Theorem 7.2]{St68}, there exists a Borel $B\subset M$
which is $f$-stable.
By Corollary \ref{cor:par-ext}, there exists
a pseudo-parabolic $P\subset \Aut(M)$
such that 
\begin{itemize}
\item $P\cap M=B$,
\item $f\in P$.
\end{itemize}
Since $\langle f \rangle$ is pseudo-completely reducible in 
$\Aut(M)$ and $\langle f \rangle$ lies in the pseudo-parabolic $P$, we can choose a cocharacter $\lambda:\Gm\to B$ suh that
the pseudo-Levi $Z_{\Aut(M)}(\lambda)$ contains $\langle f \rangle$.
The neutral component of $Z_{\Aut(M)}(\lambda)$ is a maximal torus of $M$, and is stable under $f$.


(3) See the remarks at the beginning of \cite[Section 9]{St68}.
\end{proof}



\mysection{Simultaneous diagonalization of metacyclic actions on reductive groups}~

In this section, we study when two automorphisms
of a reductive group fix a common maximal torus.

\mysubsection{Remark}
Instead of thinking about outer automorphisms of connected groups,
it is helpful to
form a semi-direct product and think about inner automorphisms
of disconnected linear algebraic groups
for two reasons:
\begin{itemize}
\item the powerful machinery of pseudo-parabolics are available, and
\item the framework of disconnected groups is more flexible
and allows us to maneuver to disconnected groups
that are not semi-direct products.
\end{itemize}
Since outer automorphism groups of semisimple connected groups
are finite, we don't lose much by allowing only considering
disconnected groups of finite type.

It is useful to pass to the adjoint quotient.

\mysubsection{Lemma}
\label{lem:ad-lift}
(1) Let $H$ be a possibly disconnected algebraic group.
Let $\pi:H'\to H$ be a homormorphism such that
$H^{\prime\circ}\to H^\circ$ is a central isogeny.
If $P_H(f)$ is a pseudo-parabolic subgroup of $H$,
then
there exists a pseudo-parabolic subgroup
$P_{H'}(f')$ of $H'$ such that
$\pi(P_{H'}(f'))=P_H(f)$.

(2) A subgroup $\Gamma$ of $H$ is pseudo-completely reducible in $H$
if and only if it is pseudo-completely reducible in the adjoint quotient of $H$.

\begin{proof}
(1) Let $T'$ be a maximal torus of $H^{\prime\circ}$.
Then $T=\pi(T')$ is a maximal torus of $H^\circ$.
We can assume $f$ is valued in $T$.
The homomorphism of cocharacter groups $X_*(T')\to X_*(T)$
has finite cokernel of cardinality $d$.
The cocharacter $f^d$ admits a lift
$f':\Gm \to T'$.
We have $\pi(P_{H'}(f'))=P_H(f^d)=P_H(f)$.

(2) It is an immediate consequence of (1).
\end{proof}

\mysubsection{Definition}
Let $\gamma\in \Aut(H)$ where $H$ is a connected algebraic group.
We say $\gamma$ acts {\it innerly} on $H$ if $\gamma$ is an inner automorphism, and we say $\gamma$ acts {\it semisimply}
on $H$ if $\gamma$ is a semisimple automorphism.

\mysubsection{Lemma}
\label{lem:central-ss}
Let $\Gamma$ be a pseudo-completely reducible subgroup of $H$,
and let $\gamma$ be a element of $\Gamma$
which acts on $H^\circ$ innerly by conjugation.
Assume
\begin{itemize}
\item $H^\circ$ is reductive,
\item $\Gamma$ is solvable, and
\item $\langle \gamma \rangle$ is a normal subgroup of $\Gamma$.
\end{itemize}
Then $\gamma$ acts on $H^\circ$ semisimply by conjugation.

\begin{proof}
By passing to the adjoint quotient, we can assume
$H^\circ$ is centerless.
Say $\gamma(g) = h g h^{-1}$ for all $g\in H^\circ$.
Since $H^\circ$ is centerless, $h$ is unique.
So for each $\sigma\in \Gamma$, $\sigma(h)=h^a$ for some $a$.
Write $h = h_sh_u$ where $h_s$ is the semisimple part
and $h_u$ is the unipotent part (the multiplicative Jordan decomposition).
For each $\sigma\in \Gamma$, by \cite[Theorem 2.4.8]{Sp98},
$\sigma(h_u)=\sigma(h)_u = (h^a)_u=h_u^a$ for some integer $a$.
So $\langle h_u\rangle$ is stable under $\Gamma$.
By \cite[Corollaire 3.9]{BT71}, there exists a parabolic subgroup
$P$ of $H^\circ$ such that
\begin{itemize}
\item $h_u\in R_u(P)$ (the unipotent radical of $P$),
\item $N_{H^\circ}(\langle h_u\rangle)\subset P$, and
\item $\Gamma$ fixes $P$.
\end{itemize}
By Corollary \ref{cor:par-ext}, there exists a cocharacter
$f:\Gm\to H^\circ$ such that
\begin{itemize}
\item $P = P_{H^\circ}(f)$, and
\item $\Gamma\subset P_H(f)$.
\end{itemize}
Since $\Gamma$ is pseudo-completely reducible in $H$,
we can choose $f$ so that
$\Gamma\subset Z_{H}(f)$.
Since $h=h_sh_u \in Z_H(f)$, we have $h_u\in Z_H(f)$ by \cite[Theorem 2.4.8]{Sp98}.
Therefore $h_u \in Z_H(f)\cap R_u(P) = \{1\}$. 
So $h=h_s$ is a semisimple element.
\end{proof}

\mysubsection{Theorem}
\label{thm:metacyclic}
Let $M$ be a possibly disconnected algebraic group
over an algebraically closed field
with reductive neutral component $M^\circ$.
Let $\langle \tau \rangle \rtimes \langle \sigma \rangle \subset M$
be a metacyclic group subgroup.
Assume $\tau$ acts on $M^\circ$ semisimply.
We have either
\begin{itemize}
\item[(I)] the neutral component of $(M^\circ)^{\tau}:=\{x\in M^\circ| \tau(x)=x\}$ is a torus, or
\item[(NI)] The subgroup $\langle \tau \rangle \rtimes \langle \sigma \rangle \subset M$
is not pseudo-irreducible in $M$.
\end{itemize}


\begin{proof}
By replacing $M$ by the subgroup generated by  $\langle \tau \rangle \rtimes \langle \sigma \rangle$ and
the derived subgroup of the neutral component $[M^\circ,M^\circ]$
we can and do assume $M^\circ$ is semi-simple.

Let $H = (M^\circ)^{\tau}$ be the fixed point subgroup.
Since $M^\circ$ is semi-simple and $\tau$ is semi-simple,
$H^\circ$ is a reductive group by \cite[Corollary 9.4]{St68}.

It is clear that $\sigma(H^\circ)=H^\circ$.
By \cite[Theorem 7.2]{St68}, there exists a Borel $B\subset H^\circ$
which is $\sigma$-stable.
Since $\tau$ acts trivially on $B$, the group
$B$ is $\langle \tau \rangle \rtimes \langle \sigma \rangle$-stable.
Write $H'$ for the subgroup of $M$ generated by
$H^\circ$ and $\langle \tau \rangle \rtimes \langle \sigma \rangle$.
By Corollary \ref{cor:par-ext}, there exists a cocharacter
$f:\Gm\to H^\circ$ such that
\begin{itemize}
\item $B = P_{H^\circ}(f)$, and
\item $\gamma, \tau\in P_{H'}(f)\subset P_M(f)$.
\end{itemize}
The group $P_M(f)$ is a pseudo-parabolic subgroup of $M$.
Since $P_M(f) \cap H^\circ = B$ and $H^\circ\subset M^\circ$,
we have either $P_M(f) \cap M^\circ \ne M^\circ$
or $P_M(f) \cap H^\circ = H^\circ$
(negating both yields an easy contradiction).
The first case implies 
 the subgroup $\langle \tau \rangle \rtimes \langle \sigma \rangle$
is not pseudo-irreducible in $M$.
The second case implies $H^\circ = B$.
Since $B$ is a Borel of $H^\circ$, it forces
$H^\circ$ to be a torus.
\end{proof}


\mysubsection{Lemma}
\label{lem:ks99}
Let $\gamma$ be a quasi-semisimple automorphism of a connected reductive group $H$ over an algebraically closed field $k$.
If $(H^{\gamma})^\circ\subset Z_H(H)$ (the center of $H$), then
$H$ is also a torus.

\begin{proof}
Choose a Borel pair $T\subset B$ of $H$ which is $\tau$-stable.
Let $U$ be the unipotent radical of $B$.
Assume $U$ is a non-trivial group.
Let $\beta$ be the highest positive root of $H$ with respect to
$B$ and $T$.
The root group $U_\beta$ is a characteristic subgroup of $U$,
and is $\tau$-stable.
Denote by $w_0$ the longest Weyl group element with respect to $T$ and $B$. It is clear that $\tau(w_0)=w_0$.
The root group $U_{-\beta} = w_0 U_{\beta}w_0^{-1}$
is also $\tau$-stable.
Write $H_\beta$ for the subgroup of $H$ generated by
$U_\beta$ and $U_{-\beta}$.
By the structure theory of reductive groups (\cite[Chapter 9]{Sp98}),
there exists a surjective homomorphism
$\psi: \SL_2\to H_\beta$.
So $H_\beta$ is either $\SL_2$ or $\PGL_2$,
and has no outer automorphisms.
In particular, the fixed-point subgroup $H_\beta^\tau$
contains a maximal torus $T_\beta$ of $H_\beta$.
By the assumption $(H^{\gamma})^\circ\subset Z_H(H)$,
we have $T_\beta \subset (Z_H(H) \cap H_\beta)^\circ \subset Z_{H_\beta}(H_\beta)^\circ = \{1\}$, which is a contradiction.
So $U$ must be a trivial group, and $B=T$.
\end{proof}

\mysubsection{Corollary}
\label{cor:metacyclic}
Let $M$ be a possibly disconnected algebraic group
whose neutral component is reductive.
Let $\langle \tau \rangle \rtimes \langle \sigma \rangle\subset M$
be a pseudo-completely reducible metacyclic subgroup in $M$.
If $\tau$ acts semisimply on $M^\circ$,
there exists a maximal torus $T$ of $M^\circ$ which
is $\langle \sigma, \tau \rangle$-stable,
and a Borel $B\supset T$ which is $\tau$-stable.

\begin{proof}
Assume without loss of generality that
$\langle \sigma, \tau \rangle$ is pseudo-irreducible in $M$.
By Theorem \ref{thm:metacyclic},
$C = ((M^\circ)^{\tau})^\circ$ is a torus.
Write $L = Z_{M^\circ}(C)$,
which is a Levi subgroup of $M^\circ$ stable under
both $\sigma$ and $\tau$-action.
We have $(L^\tau)^\circ = ((M^\circ)^\tau)^\circ \cap L
=C\cap L = C$.
By Lemma \ref{lem:ks99},
$L$ is a torus.
Since $\tau$ is a semisimple automorphism,
$\tau$ fixes a Borel pair $T\subset B$ of $M^\circ$.
Since $T$ is abelian, we have $T\subset L$.
Since $T$ is a maximal torus, we have $T=L$.
\end{proof}

\mysection{Semisimple conjugacy classes in disconnected groups}~

The key tool is \cite[Theorem 1.1.A]{KS99}. 
We include a proof here for lack of reference in characteristic $p$.

\mysubsection{Theorem}
\label{thm:ks99}
Let $M$ be a connected reductive group over a field.
Let $\theta:M\to M$ be a quasi-semisimple automorphism.
Write $M^{\theta\circ}$ for the neutral component of the fixed point group of $\theta$ in $M$.
We have
\begin{itemize}
\item[(1)] $M^{\theta\circ}$ is a reductive group,
\item[(2)] Let $T^{0}$ be a maximal torus of $M^{\theta\circ}$.
Then $Z_{M}(T^0)$ is a maximal torus of $M$ (which is clearly $\theta$-stable).
Moreover, for any Borel $B^0\supset T^0$ of $M^{\theta\circ}$,
there exists a $\theta$-stable Borel $B$ of $M^\circ$ containing 
$Z_{M}(T^0)$ and $B^0$.
\item[(3)] Conversely, let $(B, T)$ be a $\theta$-stable Borel pair of $M$. Then $(M^{\theta\circ}\cap B, M^{\theta\circ}\cap T)$ is a Borel pair of $M^{\theta\circ}$.
\item[(4)] If there is a $\theta$-stable pinning $(M, B, T, \{u_\alpha\}_{\alpha\in R(B,T)})$ of $M$,
then the map
$N_{M^{\theta\circ}}(T^{\theta\circ})/T^{\theta\circ}
\to (N_M(T)/T)^\theta$ is a bijection.
\end{itemize}

\begin{proof}
Part (1) and (3) can be found in \cite{St68}.
We first prove part (2).
Define $H:=Z_{M}(T^0)$, which is a reductive group since $M$ is reductive.
The fixed point group $H^{\theta\circ}$ is reductive by part (1),
and $H^{\theta\circ}\subset M^{\theta\circ}$
(so $T^0$ is a maximal torus of $H^{\theta\circ}$).
Since $H^{\theta\circ}$ is reductive and $T^0$ is a maximal torus thereof,
$Z_{H^{\theta\circ}}(T^0)=T^0$.
On the other hand, by definition, $Z_{H^{\theta\circ}}(T^0)=H^{\theta\circ}$. Thus $H^{\theta\circ}=T^0$.
Lemma \ref{lem:ks99} applied to the group $H$
shows $H$ is a torus.

Next we prove the ``moreover'' part.
There exists a $\theta$-stable Borel pair $(B, T)$.
By part (3), $(B^1, T^1):=(B, T)\cap M^{\theta\circ}$ is a Borel pair 
of $M^{\theta\circ}$.
Thus there exists $g\in M^{\theta\circ}$ such that
$g (B^1, T^1)g^{-1}=(B^0, T^0)$.
Since $g$ commutes with $\theta$, it is clear that
$g B g^{-1}$ is a $\theta$-stable Borel.

Now consider part (4).
Let $\alpha$ be a simple positive root of $R(B,T)$.
The root group $u_{\alpha|_{T^0}}:t\mapsto \underset{\beta|_{T^{\theta\circ}}=\alpha|_{T^{\theta\circ}}}{\prod}u_\beta(t)$ is $\theta$-stable 
of $T^{\theta\circ}$-weight $\alpha|_{T^{\theta\circ}}$.
Since $u_{\alpha|_{T^0}}$ lies in the unipotent radical of $B$,
we know it is a root group of $M^{\theta\circ}$,
and $\alpha|_{T^{\theta\circ}}\ne 1$ is a simple positive root of
$R(B^{\theta\circ}, T^{\theta\circ})$.
In particular, $\alpha\in R(B,T)$ is a positive root
if and only if $\alpha|_{T^{\theta\circ}}$ is a positive root.
Thus the $B$ is the unique Borel of $M$ containing $T$ and $B^{\theta\circ}$.
Let $w\in (N_M(T)/T)^\theta$. There exists $w'\in N_{M^{\theta\circ}}(T^{\theta\circ})/T^{\theta\circ}$
such that $w^{-1} w'\in N_M(B)/T$ since there is a bijection between Weyl chambers of $M^{\theta\circ}$
and $\theta$-stable Weyl chambers of $M$.
Thus $w^{-1} w'\in T/T$, and $w=w'$.
\end{proof}
We remark that part (4) of Theorem \ref{thm:ks99}
is also discussed in the paragraph after \cite[Theorem 1.1.A]{KS99};
and their argument works in characteristic $p$ without change.

\mysubsection{Definition}
In a disconnected algebraic group $M$,
two elements $x$ and $y$ are said to be {\it conjugate to each other}
if $x=gyg^{-1}$ for some $g\in M^{\circ}$.

\mysubsection{Lemma}
\label{lem:species-1}
Let $M$ be a possibly disconnected reductive group over an algebraically closed field.
There is a bijection
$$
\{\text{Conjugacy classes of tuples $(s, T^0)$}\}
\xrightarrow{(s, T^0)\mapsto s}
\{\text{Semisimple conjugacy classes in $M$}\}
$$
where a tuple $(s, T^0)$ consists of
a semisimple element $s\in M$
and a maximal torus $T^0$ of the fixed-point group $(M^{\circ})^{s\circ}$.

\begin{proof}
Let $(s_1, T^0_1)$ and $(s_2, T^0_2)$ be two pairs
such that $s_1=gs_2g^{-1}$ for some $g\in M^{\circ}$.
Then both $T^0_1$ and $g T^0_2 g^{-1}$
are maximal tori of the reductive group
$(M^{\circ})^{s_1\circ}$.
So we can choose $h\in (M^{\circ})^{s_1\circ}$
such that $h T^0_1 h^{-1} = g T^0_2 g^{-1}$.
Since $h s_1 h^{-1}=s_1$, we have
$h(s_1, T^0_1)h^{-1} = g(s_2, T^0_2)g^{-1}$.
\end{proof}

\mysubsection{Lemma}
\label{lem:species-2}
Let $\theta, \theta'\in M$ be semisimple elements
such that $\theta^{-1}\theta'\in M^\circ$.
Let $(B, T)$ be a $\theta$-stable Borel pair of $M^\circ$,
There exists a conjugate $g \theta' g^{-1}$, $g\in M^\circ$
such that $\theta^{-1}g \theta' g^{-1} \in T$.

\begin{proof}
Choose any Borel pair $(B^0, T^0)$
of $(M^{\circ})^{\theta'\circ}$.
By Theorem \ref{thm:ks99},
there exists a $\theta'$-stable Borel pair $(B', T')$
such that $(B^0, T^0)=(B',T')\cap M^{\circ\theta'\circ}$.
There exists $g\in M^\circ$
such that
$g(B', T')g^{-1}=(B,T)$.
For ease of notation, we assume $(B',T')=(B,T)$.
Now $(B, T)$ is simultaneously $\theta$-stable and $\theta'$-stable.
We have $\theta^{-1}\theta'\in N_{M^\circ}(T)\cap N_{M^\circ}(B)=T$.
\end{proof}


\mysubsection{Construction of $\xi_T$}
\label{cons:xiT}
Let $M$ be a disconnected reductive group over an algebraically closed field $k$. 
Fix a connected component $M^\circ \theta$ of $M$ that that $\theta$ is semisimple.

Fix a $\theta$-stable Borel pair $(B, T)$ of $M^\circ$.
Note that
$\theta$ acts on $X_*(T)$
and $\Omega:=N_{M^\circ}(T)/T$.
Write $X_*(T)_\theta$ for the $\theta$-coinvariants
$X_*(T)/(1-\theta)$, and write $\Omega^\theta$
for the subgroup of $\theta$-fixed points of $\Omega$.
Denote by $X_*(T)_{\theta,\Tf}$ the (maximal) torsion-free quotient
of the abelian group $X_*(T)_\theta$.

We construct a map

\bigmap
{Semisimple conjugacy classes of M that lie in the component $M^\circ\theta$}{\xrightarrow{\xi_T}}
{$X_*(T)_{\theta, \Tf}\otimes k^\times/ \Omega^{\theta}$}

\noindent
as follows.
Let $[\theta']$ be a semisimple conjugacy class of $M$
such that $\theta^{-1}[\theta']\subset M^\circ$.
Define $\xi_T([\theta'])$
to be the equivalence class of $\theta'\theta^{-1}\in T(\bFp)=X_*(T)\otimes \bFp^\times$.

\mysubsection{Proposition}
\label{prop:sscd}
If there is a $\theta$-stable pinning $(M^\circ, B, T,  \{u_\alpha\}_{\alpha\in R(B,T)})$,
then
the map $\xi_T$ in Paragraph \ref{cons:xiT}
is well-defined and is a bijection.

Moreover, each semisimple conjugacy class of $M$
that lies in the component $M^\circ \theta$
admits a representative in $T^{\theta\circ}\theta$.

\begin{proof}
Let $[\theta']$ be a semisimple conjugacy class
contained in $M^\circ \theta$.
By Lemma \ref{lem:species-2},
there exists a representative $\theta'$
such that $t=\theta'\theta^{-1}\in T$.
Now $T^{\theta'\circ}=T^{\theta (\theta^{-1}t\theta) \circ}= T^{\theta \circ}$.
Write $T^0$ for $T^{\theta \circ}$.
By Theorem \ref{thm:ks99},
$T^0$ is a maximal torus of
both $M^{\circ \theta\circ}$
and $M^{\circ \theta'\circ}$.
By Lemma \ref{lem:species-1},
if $\theta'_1$ is another representative of $[\theta']$
such that $t_1:=\theta'_1\theta^{-1}\in T$,
then there exists an element $g\in M^\circ$
such that $g(\theta'_1, T^0)g^{-1} = (\theta', T^0)$;
in other words, $g\in N_{M^\circ}(T^0)$.
By Theorem \ref{thm:ks99}, $Z_{M^\circ}(T^0)=T$,
and thus $N_{M^\circ}(T^0)\subset N_{M^\circ}(T)$.

So far, we have shown $\xi_T$ defines a bijection onto
$(X_*(T)\otimes k^\times)/N_{M^\circ}(T^\circ)(k)$,
with the caveat that the action of $N_{M^\circ}(T^\circ)(k)$
is $\theta$-twisted, that is,
\begin{equation}
\label{eq:q-twist-action}
w \cdot t = w t \theta w^{-1} \theta^{-1},
\end{equation}
$t\in T$, $w\in N_{M^\circ}(T^\circ)(k)$.
Next we show 
$$
(X_*(T)\otimes k^\times)/T(k) = X_*(T)_\theta \otimes k^\times.
$$
Let $s\in T(k)$, and $\theta'\in T(k) \theta$.
We have
$$
(s \theta' s^{-1}) \theta^{-1}
= s (\theta' \theta^{-1}) (\theta s^{-1} \theta^{-1})
= s (\theta s^{-1} \theta^{-1}) (\theta' \theta^{-1}).
$$
Thus
$(X_*(T)\otimes k^\times)/T(k) = \frac{X_*(T)\otimes k^\times}{(1-\theta)X_*(T)\otimes k^\times}$.
Since $-\otimes k^\times$ is right-exact, we have
$\frac{X_*(T)\otimes k^\times}{(1-\theta)X_*(T)\otimes k^\times}
=\frac{X_*(T)}{(1-\theta)X_*(T)}\otimes k^\times$.
Write $(-)_{\Tors}$ for the torsion part of an abelian group.
There exists a short exact sequence
$$
0\to X_*(T)_{\theta,\Tors}\to X_*(T)_{\theta}\to X_*(T)_{\theta, \Tf}\to 0
$$
where $X_*(T)_{\theta, \Tf}$ is the maximal torsion-free quotient of
$X_*(T)_\theta$.
Since $k^\times$ is a divisible abelian group, we have
$$
X_*(T)_\theta \otimes k^\times = X_*(T)_{\theta,\Tf}\otimes k^\times.
$$
By part (4) of Theorem \ref{thm:ks99},
we have
$$
N_{M^\circ}(T^\circ)/T = (N_{M^\circ}(T)/T)^\theta = N_{M^{\theta\circ}}(T^0)/T^0.
$$
Thus $\xi_T$ defines a bijection onto
$$
(X_*(T)_\theta \otimes k^\times)/(N_{M^{\theta\circ}}(T^0)/T^0).
$$
Moreover, the $\theta$-twisted action (\ref{eq:q-twist-action})
becomes untwisted;
indeed, for $t\in T$ and $\bar w\in N_{M^{\theta\circ}}(T^0)/T^0$,
we have
$w t \theta w^{-1} \theta^{-1} = w t w^{-1}$,
where $w\in N_{M^{\theta\circ}}(T^0)$ is a representative of $\bar w$.

{\bf Claim} $X_*(T^0)\otimes k^\times \to X_*(T)_{\theta, \Tf}\otimes k^\times$ is surjective.
\begin{proof}
By \cite[Proposition 13.2.4]{Sp98},
$T$ is an almost direct product $T^0T_a$
such that $T^0\cap T_a$ is finite.
Since $X_*(T)_{\theta,\Tf}$ is torsion-free, it defines
a quotient torus $T_0$ of $T$.
The composite $T_a \to T \to T_0$ must be trivial
(if otherwise $T_a^{\theta\circ}$ is non-trivial).
The composite $T^0\to T^0T_a/T_a \to T_0$
is thus surjective.
\end{proof}
Translating the claim to a statement about conjugacy classes,
we see there exists $s\in T$ such that
$s \theta' s^{-1} \in T^0\theta = \theta T^0$.
\end{proof}

\mysection{$\theta$-twisted semisimple conjugacy classes}~

Sometimes it is notationally easier to discuss
$\theta$-twisted conjugacy classes than
conjugacy classes in a disconnected group.
We keep notations and assumptions 
that are used in Proposition \ref{prop:sscd}.

In this subsection, $M$ is defined over an algebraically closed field
$k$ of characteristic $p$.
Let $q$ be a power of $p$.
Write $\Frob_q:x\mapsto x^q$ be the $q$-power map.

\mysubsection{Definition}
Two elements $g$, $g'\in M^\circ$ are said to be
$\theta$-conjugate if $g \theta$ and $g' \theta$
are $M^\circ$-conjugate.
Concretely, it means there exists $h\in M^\circ$
such that $g' = h g (\theta h^{-1} \theta^{-1})$.
Denote by $[g]_\theta$ the $\theta$-twisted conjugacy class
of $g$.

\mysubsection{Definition}
\label{def:QF}
A set-theoretic map $F:M^\circ\to M^\circ$
is said to be a {\it $\theta$-twisted Frobenius endomorphism}
if
\begin{itemize}
\item[(QF1)] $F(T)\subset T$;
\item[(QF2)] There exists an automorphism $\varphi\in\End_{\bZ}(X_*(T))$ such that $F|_{T^{\theta\circ}}=\varphi\otimes \Frob_q$
under the identification $T = X_*(T)\otimes k^\times$;
\item[(QF3)] If $y\in [x]_\theta$, then $F(y)\in [F(x)]_\theta$;
\item[(QF4)] $F(\theta^{q} x \theta^{-q})=\theta F(x) \theta^{-1}$
for all $x\in M^\circ$.
\end{itemize}

We say $[g]_\theta$ is $F$-stable
if $F(g')\in [g]_\theta$ for all $g'\in [g]_\theta$.

\mysubsection{Proposition}
\label{prop:QF}
If $F:M^\circ\to M^\circ$
is a {$\theta$-twisted Frobenius endomorphism},
then there exists a bijection

\bigmap
{$F$-stable $\theta$-twisted semisimple conjugacy classes in $M^\circ$}
{\cong}
{$(X_*(T)_{\theta,\Tf}\otimes k^\times/\Omega^{\theta})^{\varphi\otimes \Frob_q}$.}

\begin{proof}
By Proposition \ref{prop:sscd},
the bijection is automatic if
the $\varphi$-action descends to 
$X_*(T)_{\theta,\Tf}$,
which follows immediately from Definition \ref{def:QF}.
\end{proof}

\mychapter{Mod $p$ Langlands-Shelstad factorization}

\etocsettocdepth{2}
\localtableofcontents

\vspace{3mm}
\label{sec:LS}

Let $G$ be a connected quasi-split reductive group over $F$.
Fix an $F$-pinning $(B, T, \{X_\alpha\})$ of $G$.
The pinned group $(G, B, T, \{X_\alpha\})$
has a dual pinned group
$(\wh G, \wh B, \wh T, \{Y_\alpha\})$ defined over $\bZ$,
together with an isomorphism of based root data $\Psi_0(\wh G, \wh B, \wh T)\cong \Psi_0(G, B, T)^\vee$.

Let $L\subset F^s$ be a splitting field of $G$, that is, a Galois subfield of $F^s$  such that the Galois action on $\Psi_0(G_{F^s}, B_{F^s}, T_{F^s})$ factors through $\Gal(L/F)$.
The Galois action on the based root datum induces a Galois action on the corresponding pinned reductive group $\Gal(L/F)\to \Aut(\wh G, \wh B, \wh T, \{Y_\alpha\})$.
Write $\lsup L G$ for the semi-direct product $\wh G \rtimes \Gal(L/F)$.
See \cite[Section 1.1]{Zhu21} and \cite{BG14} for more details.

\mysection{Four types of $L$-groups}~

In the literature, four different kinds of $L$-groups
are used, and depending on the context,
they are not always interchangeable.
We list these $L$-groups as follows:

\vspace{3mm}

\begin{tabular}{l|l}
Galois form & $\wh G \rtimes \Gal(L/F)$\\
Absolute Galois form & $\wh G \rtimes \Gal_F$\\
Weil form & $\wh G \rtimes W_F$ \\
Relative Weil form & $\wh G \rtimes W_{L/F}$
\end{tabular}

\mysubsection{Definition}
Let $A$ be a ring.

A group homomorphism
$$f(A): \wh G(A) \rtimes \Gal(L/F)\to \wh H(A) \rtimes \Gal(L/F)$$
is said to be an {\it $L$-homomorphism}
if it sends $g\times \sigma$ to $h\times \sigma$
for all $\sigma\in \Gal(L/F)$.

A group homomorphism
$$f(A): \wh G(A) \rtimes \Gal_F\to \wh H(A) \rtimes \Gal_F$$
is said to be an {\it $L$-homomorphism}
if it sends $g\times \sigma$ to $h\times \sigma$
for all $\sigma\in \Gal_F$,
and there exists an open subgroup $\Gamma\subset \Gal_F$
such that $1\times \sigma$ is sent to $1\times \sigma$
for all $\sigma\in \Gamma$.

A group homomorphism
$$f(A): \wh G(A) \rtimes W_F\to \wh H(A) \rtimes W_F$$
is said to be an {\it $L$-homomorphism}
if it sends $g\times \sigma$ to $h\times \sigma$
for all $\sigma\in \Gal_F$,
and there exists an open subgroup $\Gamma\subset I_F$
such that $1\times \sigma$ is sent to $1\times \sigma$
for all $\sigma\in \Gamma$.

A group homomorphism
$$f(A): \wh G(A) \rtimes W_{L/F}\to \wh H(A) \rtimes W_{L/F}$$
is said to be an {\it $L$-homomorphism}
if it sends $g\times \sigma$ to $h\times \sigma$
for all $\sigma\in W_{L/F}$.

\mysubsection{Lemma}
\label{lem:Lgrp-1}
Let $A$ be a finite ring.

(1)
If 
$$f(A): \wh G(A) \rtimes \Gal_F\to \wh H(A) \rtimes \Gal_F$$
is an $L$-homomorphism,
then
the set $\sigma\in \Gal_F$
such that $f(A)(1\times \sigma)=1\times \sigma$
is an open normal subgroup.

(2)
If 
$$\wh G(A) \rtimes \Gal_F\to \wh H(A) \rtimes \Gal_F$$
is an $L$-homomorphism,
then it descends to
$$\wh G(A) \rtimes \Gal(L/F)\to \wh H(A) \rtimes \Gal(L/F)$$
for some finite extension $L/F$.

(3)
If
$$\wh G(A) \rtimes \Gal(L/F)\to \wh H(A) \rtimes \Gal(L/F)$$
is an $L$-homomorphism,
then it can be uniquely lifted
to an $L$-homomorphism
$$\wh G(A) \rtimes \Gal_F\to \wh H(A) \rtimes \Gal_F.$$

(4)
If 
$$\wh G(A) \rtimes W_F\to \wh H(A) \rtimes W_F$$
is an $L$-homomorphism,
then it descends to
$$\wh G(A) \rtimes \Gal(L/F)\to \wh H(A) \rtimes \Gal(L/F)$$
for some finite extension $L/F$.

(5)
If
$$\wh G(A) \rtimes \Gal(L/F)\to \wh H(A) \rtimes \Gal(L/F)$$
is an $L$-homomorphism,
then it can be uniquely lifted
to an $L$-homomorphism
$$\wh G(A) \rtimes W_F\to \wh H(A) \rtimes W_F.$$

\begin{proof}
(1)
Note that
$(1\times \sigma_1)(1\times \sigma_2)
=(1\times \sigma_1 \sigma_2)$.

(2)
In a profinite group, an open subgroup
is closed of finite index.
So part (2) follows from part (1).

(3)
Let $\sigma\in \Gal_F$ and write $\bar \sigma$ for its image in $\Gal(L/F)$.
If $g\times \bar\sigma\mapsto h\times \bar\sigma$,
set $g\times \sigma\mapsto h\times \sigma$.
It is easy to check it is a well-defined homomorphism.

(4)
Let $\theta\in W_F$ be a Frobenius element,
and say $1\times \theta\mapsto h\times \theta$.
Since $\wh H(A)$ is a finite group,
there exists an integer $n$ such that $(h\times\theta)^n=1\times \theta^n$
acts trivially on $\wh H(A)$.
So there exists an open subgroup $\Gamma$ of $W_F$ of finite index 
such that $1\times \sigma\mapsto 1\times \sigma$ for all $\sigma\in \Gamma$
and that $\Gamma$ acts trivially on $\wh H(A)$.

(5) It is similar to part (3).
\end{proof}

By Lemma \ref{lem:Lgrp-1},
the various forms of $L$-groups
make little difference when working with finite coefficients.

\mysubsection{Definition}
Let $A$ be a finite ring equipped with discrete topology, and
let $G$ and $H$ be two reductive groups over $F$.
Let $L$ be a sufficiently large field extension of $F$
that splits both $G$ and $H$.

An {$L$-morphism}
$f:\lsup LG_A\to \lsup LH_A$
is an algebraic group morphism over $A$
such that for all $A$-algebras $B$ with finitely many elements,
$f(B)$ is an $L$-homomorphism.

Write $\lsup L\{*\}$ for the $L$-group
of the trivial group.
An {\it $L$-parameter for $G$ with $A$-coefficients}
is an $L$-morphism
$\lsup L\{*\} \to \lsup LG_A$, defined up to $\wh G(A)$-conjugacy.
Equivalently, since $\lsup L\{*\}$
is a constant group scheme,
an $L$-parameter can also be defined as an $L$-homomorphism
$\Gal_F\to \lsup LG(A)$, defined up to $\wh G(A)$-conjugacy.

Let $R$ be a profinite ring.
A {\it profinite $L$-parameter for $G$ with $R$-coefficients}
is a compatible system of
$L$-parameters $\Gal_F \to \lsup LG(A)$,
where $A$ is a finite quotient of $R$.


\vspace{3mm}

From now on, assume $G$ is tamely ramified (that is, the splitting field $L$ can be chosen as a tame extension of $F$).

\mysection{Quasi-semisimplicity of semisimple mod $p$ $L$-parameters}
\label{sec:qs}

\mysubsection{Standard parabolic subgroups and Levi subgroups
of Langlands dual groups}
\label{def:std-par}
Write $\wh \Delta$ for $\Delta(\wh B, \wh T)$.
Any $\Gal(L/F)$-stable subset $S\subset \wh \Delta$ corresponds to a pinned subgroup $(M_S, M_S\cap \wh B, \wh T, \{Y_\alpha\}|_{S})\subset (\wh G, \wh B, \wh T, \{Y_\alpha\})$
which is stable under the $\Gal(L/F)$-action.
Write $P_S$ for a parabolic subgroup of $\wh G$ having $M_S$ as a Levi subgroup.
We call $M_S \rtimes \Gal(L/F)$ a {\it standard Levi subgroup} of $\lsup L G$,
and we call $P_S\rtimes \Gal(L/F)$ a {\it standard parabolic subgroup} of $\lsup L G$.

\mysubsection{Lemma}
\label{lem:std-par}
All big pseudo-parabolics of $\lsup LG$ are conjugate to
a standard parabolic.

\begin{proof}
Let $P$ be a big pseudo-parabolic of $\lsup LG$.
There exists a Borel subgroup $B$ of $P^\circ$.
Let $B_{\std}$ be the standard parabolic whose
neutral component is a Borel.
Then there exists an element $g\in \wh G$
such that $gBg^{-1}= B_{\std}^\circ$.
After replacing $B_{\std}$ by $g^{-1}B_{\std}g$,
we can assume $B\subset B_{\std}$.
For each $\bar \gamma\in \Gal(L/F)$,
let $\gamma\in P$ be a lift of $\gamma$ in $P$
and let $\gamma_{\std}\in B_{\std}$
be the ``standard'' lift of $\gamma$ in $B_{\std}$
(note that $B_{\std} = B \rtimes \Gal(L/F)$
admits a distinguished copy of $\Gal(L/F)$).
There exists $h\in P$ such that
$\gamma B \gamma^{-1}= h B h^{-1}$.
We have $h^{-1}\gamma \gamma_{\std}^{-1}(B)\gamma_{\std}\gamma^{-1}h = h^{-1}\gamma (B)\gamma^{-1}h = B$.
Since $B$ is self-normalizing in $\wh G$,
we have $h^{-1}\gamma \gamma_{\std}^{-1}\in B$.
Thus $\gamma_{\std}\in P$.

The same argument applies to all $\bar\gamma\in\Gal(L/F)$.
So $B_{\std}\subset P$.
Now it is clear that all positive simple roots
of $P$ with respect to the standard $(\wh B, \wh T)$-pair
are permuted by the $\Gal(L/F)$-action,
and $P$ is a standard parabolic.
\end{proof}

\mysubsection{Definition}
\label{def:par-sub}
For a ring $k$,
a \textit{parabolic subgroup} of $\lsup L G_k$
is a subgroup $\lsup LP$ conjugate to a standard parabolic
of $\lsup LG_k$.
Equivalently, a parabolic of $\lsup LG_k$
is a big pseudo-parabolic.

A {\it Levi subgroup} of a parabolic $\lsup LP$
is a subgroup which is conjugate to a standard Levi
of $\lsup LG$.

\mysubsection{Remark}
We will avoid using the term ``parabolic subgroup''
because it is not consistent with the usual definition
that $P$ is parabolic if $G/P$ is a projective variety.

Instead, we prefer the use ``big pseudo-parabolic''
to avoid confusions.

We also remark that the lemma above is true
only because (1) the group $\lsup LG$ is itself a semi-direct product, and (2) (for connected groups) parabolic subgroups are self-normalizing.
For a general (disconnected) reductive subgroup $H$ of $\lsup LG$,
we should not expect its big pseudo-parabolics to be a split extension of its neutral component,
or conjugate to a standard parabolic in general.

\mysubsection{Definition}
\label{def:sslpar}
Let $k$ be a field.
Let $\rho:\Gal_F\to \lsup L G(\bar k)$ be an $L$-parameter.
\begin{itemize}
\item 
$\rho$
is said to be {\it pseudo-irreducible}
if the image of $\rho$ is pseudo-irreducible (Definition \ref{def:pirr});
\item
$\rho$ is said to be {\it pseudo-semisimple} or {\it semisimple} if the image of $\rho$ is pseudo-completely reducible;
\item
$\rho$ is said to be \textit{quasi-semisimple}
if there exists a maximal torus $T$ of $\wh G$ such that
\begin{itemize}
\item $\rho(I_L)\subset T(k)$, and
\item $\Img \rho \subset N_{\lsup LG}(T)(k)$.
\end{itemize}
The term ``quasi-semisimple'' appears in \cite[Section 9]{St68}, which is loosely related to our situation.
\end{itemize}


\mysubsection{Lemma}
\label{lem:tame-ramification}
If $\rho:\Gal_F\to {\lsup LG}(\bFp)$ is a semisimple $L$-parameter,
then $\rho$ is tamely ramified.

\begin{proof}
Let $P_F\subset \Gal_F$ be the wild inertia.
The image $\rho(P_F) \subset \wh G(\bFp)$ is a $p$-group, and thus consists of unipotent elements.
By \cite[Corollaire 3.9]{BT71}, there exists a parabolic subgroup $P$ of $\wh G_{\bFp}$ with unipotent radical $R_u(P)$
such that
\begin{itemize}
\item $\rho(P_K)\subset R_u(P)(\bFp)$, 
\item $N_{\wh G}(\rho(P_K))\subset P(\bFp)$, and
\item all automorphism of $\wh G_{\bFp}$
which fix $\rho(P_K)$ also fix $P$.
\end{itemize}
Here $N_{\wh G}(\rho(P_K))$ is the normalizer of $\rho(P_K)$ in $\wh G$.
Since $P_K$ is a normal subgroup of $\Gal_L$, $\rho(\Gal_L)\subset N_{\wh G}(\rho(P_K))\subset P(\bFp)$.
Since $P_K$ is a normal subgroup of $\Gal_F$ and all automorphism of $\wh G_{\bFp}$
which fix $\rho(P_K)$ also fix $P$,
the subset $\Gamma:=P(\bFp)\rho(\Gal_F)\subset {\lsup LG}(\bFp)$
is a subgroup of $\lsup LG(\bFp)$.
We have $\Gamma \cap \wh G = P(\bFp) N_{\wh G}(\rho(P_K))=P(\bFp)$.
By Lemma \ref{lem:big-par},
$\Gamma$ is a big pseudo-parabolic of $\lsup LG_{\bFp}$.
Since $\rho$ is semisimple,
$\rho$ factors through a pseudo-Levi $M$
of $\Gamma$.
So
$
\rho(P_K)\subset R_u(P)(\bFp) \cap M = \{1\}
$.
\end{proof}

\mysubsection{Theorem}
\label{thm:ab-qss}
Semisimple $L$-parameters $\rho:\Gal_F \to {\lsup LG}(\bFp)$ are quasi-semisimple.

\begin{proof}
By Lemma \ref{lem:tame-ramification},
$\rho$ is tamely ramified.
Write $\mathfrak{F}$ for a Frobenius element of $\Gal_F$
and let $\mathfrak{T}$ be an element of $\Gal_F$
whose image in the tame quotient is a topological generator.
Write $\sigma:=\rho(\mathfrak{F})$
and $\tau:=\rho(\mathfrak{T})$.
We have $\sigma\tau\sigma^{-1}=\tau^{q}$.
Since $\tau$ has prime-to-$p$ order,
$\tau$ acts on $\wh G$ semisimply.
The theorem follows from Corollary \ref{cor:metacyclic}.
\end{proof}

\mysubsection{Definition}
Let $\rho$ be a tamely ramified $L$-parameter
$\Gal_F\to \lsup LG(\bFp)$.
Let $\lsup LP$ be minimal among all big pseudo-parabolic of $\lsup LG$
through which $\rho$ factors.
Let $\pi: \lsup LP\to \lsup LM$ be the quotient by the unipotent radical $U$ map,
and let $\iota:\lsup LM\hookrightarrow \lsup LP$ be
a splitting of the semi-direct product $1\to U \to \lsup LP\to \lsup LM\to 1$.
We say $\iota\circ \pi\circ \rho$ is a
{\it pseudo-semisimplification} of $\rho$,
and denote it by $\rho^{\Sems}$.

\mysubsection{Lemma}
\label{lem:Frob-ss-2}
For any tamely ramified $L$-parameter $\rho:\Gal_F\to \lsup LG(\bFp)$,
$\rho|_{I_F}$ is conjugate to $\rho^{\Sems}|_{I_F}$.

\begin{proof}
Let $\theta \in I_F$ be a topological generator of the tame
quotient of the inertia.
Let $\lsup LP$ be minimal among all big pseudo-parabolic of $\lsupp LG$
through which $\rho$ factors.
Since $\rho(\theta)$ is semisimple,
it fixes a Borel pair $(B_P, T_P)$ of $\lsup LP$
(\cite[Theorem 7.5]{St68}).
$(\pi(B_P),\pi(T_P))$ is a Borel pair of $\lsup LM$ fixed by $\pi(\rho(\theta))$,
and thus $(\iota\pi(B_P), \iota\pi(T_P))$ is a Borel pair of $\iota(\lsup LM)$ fixed by $\rho^{\Sems}(\theta)$.
Furthermore, $(\iota\pi(B_P)U, \iota\pi(T_P))$
is a Borel pair of $\lsup LP$.
There exists $g\in U$ such that
$g(B_P, T_P)g^{-1} = (\iota\pi(B_P)U, \iota\pi(T_P))$.
By replacing $\iota$ by a conjugate, we may assume $g=1$.
Write $t$ for $\rho(\theta)\rho^{\Sems}(\theta)^{-1}\in (\lsup LP)^\circ$.
We have $t \in N_{\wh G}(B_P) \cap N_{\wh G}(T_P)=T_P$.
By the construction of $\rho^{\Sems}$, we have $t\in U$.
So $t\in T_P\cap U=\{1\}$.
\end{proof}

\mysection{Maximally unramified tori}~

\mysubsection{Convention}
All homomorphisms of algebraic $F$-groups
$f:G\to H$
are meant to be defined over $F^s$,
unless it is said explicitly that it is defined over $F$.
For ease of notation,
we write $G$ for $G(F^s)$
if it is clear from the context that we are treating
it as a set.

Recall that we fixed an $F$-pinning $(B, T, \{X_\alpha\})$
of $G$.
Write $\Omega$ for the Weyl group $N_{G}(T)/T:=N_{G(F^s)}(T(F^s))/T(F^s)$,
which is canonically identified with
$N_{\wh G}(\wh T)/\wh T$.

\vspace{3mm}

Following \cite[Section 5.1]{Kal19a},
it is conceptually easier to think about framed maximal tori
(see Definition \ref{def:fmtori} and \ref{def:fmtori2})
instead of maximal tori.

Let $j: S \to G$ be a framed maximal tori.
By replacing $j$ by a conjugate $g j g^{-1}$, $g\in G(F^s)$,
we may assume $j(S)\subset T$.
Thus there is a canonical bijection

\bigequiv{$G(F^s)$-conjugacy classes of framed maximal tori of $G$}
{$\Omega$-conjugacy classes of pairs $(S, j)$
where $S$ is a $F$-torus and $j:S_{F^s}\to T_{F^s}$
is an isomorphism of tori}

Similarly, there is a bijection on the dual side

\bigequiv{$\wh G$-conjugacy classes of framed maximal tori of $\lsupp LG$}
{$\Omega$-conjugacy classes of pairs $(S, \wh \j)$
where $S$ is a $F$-torus and $\wh \j: \wh S \to \wh T$
is an isomorphism of tori}

By sending $(S, j)\mapsto (S, (\wh j)^{-1})$
where $\wh j:\wh T\to \wh S$ is the functorial dual torus map,
we get a duality between the geometric conjugacy classes
of maximal tori of $G$ and that of $\lsupp LG$.

Recall the following fact (\cite[Fact 3.4.1]{Kal19a}):

\mysubsection{Fact}
\label{fact:max-unr}
Let $S\subset G$ be a maximal torus and let $S_s\subset S$
be the maximal unramified torus (both embeddings are defined over $F$).
The following statements are equivalent.
\begin{itemize}
\item $S_s$ is of maximal dimension among the unramified subtori of $G$.
\item $S_s$ is not properly contained in an unramified subtorus of $G$.
\item $S$ is the centralizer of $S_s$ in $G$.
\item The action of $I_F$ on $R(S,G)$ preserves a set of positive roots.
\end{itemize}

\mysubsection{Definition}
\label{def:max-unr-tori}
Let $j:S\to G$ be a framed maximal torus.
The inertia $I_F$ acts on $X^*(S)$ by $\theta \cdot \alpha:= \theta \circ \alpha \circ \theta^{-1}$,
and thus acts on $X^*(j(S))$ by transport
($\theta\cdot \alpha:=(\theta\cdot (\alpha\circ j))\circ j^{-1}$).
Note that in general the $I_F$-action 
on $X^*(j(S))$ does not need to preserve the set of absolute roots
$R(j(S), G)\subset X^*(j(S))$.

We say $(S, j)$ is {\it maximally unramified}
if $I_F$ preserves $R(j(S), G)$
    and a set of positive roots thereof.

\mysubsection{Lemma}
Let $S\subset G$ be a maximal $F$-torus.
Then the tautological embedding $j:S \hookrightarrow G$
defines a maximally unramified framed maximal torus
if and only if $S$ satisfies
the equivalent statements in Fact \ref{fact:max-unr}.

\begin{proof}
Since $j$ is the tautological embedding of an $F$-torus,
$I_F$ preserves $R(j(S), G)$.
The Lemma now follows from the last bullet point of Fact \ref{fact:max-unr}.
\end{proof}

\mysubsection{Lemma}
$(S, j)$ is a maximally unramified framed torus of $G$
if and only if a geometric conjugate of $(S, j)$
    is maximally unramified.

\begin{proof}
Write $\Int(g)$ for $x\mapsto g x g^{-1}$, $g\in G(F^s)$.
We have 
$$
R(\Int(g)\circ j(S), G)
=
\{\alpha\circ\Int(g)^{-1}|\alpha\in R(j(S), G)\}.
$$
So it is clear $I_F$ preserves $R(\Int(g)\circ j(S), G)$
if and only if it preserves $R(j(S), G)$.
If $I_F$ preserves a set of positive roots
$R^+\subset R(j(S), G)$
than $\Int(g)^{-1}(R^+)$ is a set of positive roots
of $R(\Int(g)\circ j(S), G)$ that $I_F$ preserves.
\end{proof}

\mysubsection{Definition}
\label{def:max-unr-tori-dual}
Let $(S, \wh \j)$, $\wh \j:\wh S\to \wh G$
be a framed maximal torus of $\lsup LG$.
The $I_F$-action on $X^*(\wh S)$
transports to $X^*(\wh \j (\wh S))$.
We say $(S, \wh \j)$ is maximally unramified
if $I_F$ preserves $R(\wh \j (\wh S), \wh G)$
and a set of positive roots thereof.

\mysubsection{Lemma}
\label{lem:dually-max-unr}
Let $(S, j)$ be a framed maximal torus of $G$
whose geometric conjugacy class is $\Gal(F^s/F)$-stable.
Then $(S, j)$ is maximally unramified if and only if
its dual $(S, \wh \j)$, $\wh \j:\wh S\to \wh G$
is a maximally unramified framed maximal torus of $\lsup LG$.

\begin{proof}
By Theorem \ref{thm:kot},
we can assume $S\subset G$ is an $F$-subtorus.
Let $R^+\subset R(S, G)$ be a set of positive roots
preserved by $I_F$.
$R^+$ defines a Borel $B\supset S_{F^s}$ of $G_{F^s}$.
Since $I_F$ permutes the root groups $U_\alpha$ ($\alpha\in R^+$) of $B$,
$B$ is $I_F$-stable.
As a consequence, $I_F$ also permutes
the positive coroots of $B$,
which translates to the maximal unramifiedness of $(S, \wh \j)$.
\end{proof}

We summarize the duality in the following proposition.

\mysubsection{Proposition}
\label{prop:dual-max-unr}
There is a canonical bijection:

\bigequiv
{$\Gal_F$-stable $G(F^s)$-conjugacy classes
of maximally unramified framed maximal torus of $G$}
{$\Gal_F$-stable conjugacy classes
of maximally unramified framed maximal torus of $\lsup LG$}

\vspace{3mm}

Recall that a {\it twisted Levi} of $G$ is
a subgroup $M\subset G$ defined over $F$
which becomes a Levi subgroup after bash change to $F^s$.
We have the following fact.

\mysubsection{Langlands-Shelstad $L$-embedding} (\cite[Lemma 5.2.6]{Kal19a})
Let $M\subset G$ be a tame twisted Levi,
and let $\wh M_{\bC}\to \wh G_{\bC}$ be the natural inclusion.
Write $\lsupw LM$ and $\lsupw LG$ for the Weil form of the $L$-group.
There exists an extension of $\wh M\to\wh G$
to a tame $L$-embedding ${\lsupw LM}_{\bC} \to {\lsupw LG}_{\bC}$.

The conjugacy class of the tame $L$-embedding ${\lsupw LM}_{\bC} \to {\lsupw LG}_{\bC}$ depends on the choice of the so-called {\it $\chi$-data}.
We remark that it is important to use the Weil form.
An $L$-embedding $\lsup LM_{\bC}\to\lsup LG_{\bC}$
of the absolute Galois form does not exist in general
(see \cite{Kal21}).

When $M\subset G$ be a maximally unramified torus (Definition \ref{def:max-unr-tori}),
there is a canonical choice of $\chi$-data (the unramified $\chi$-data).
We establish the integral version
in the following theorem.

\mysubsection{Theorem}
\label{thm:LS-maxunr}
Let $S\subset G$ be a maximal $F$-torus
such that the tautological embedding $S\hookrightarrow$
defines a maximally unramified framed maximal torus.

For each choice of Borel $B_S\supset S$
of $G$ defined over $F^s$,
there exists an element $g\in G(F^s)$
such that $(B_S, S) = g(B, T)g^{-1}$
(recall that we fixed a pinning $(B, T, \{X_\alpha\})$ of $G$).
The conjugation-by-$g$ map $\Int(g): T\to S$
(which is only defined over $F^s$)
has a dual map $\wh{\Int(g)}:\wh S\to \wh T$.

Let $A$ be a finite ring.
Feeding the unramified $\chi$-data into
the recipe in \cite[Section 2.6]{LS87}, we get
an $L$-embedding
$\lsup Lj:\lsup LS_A \to \lsup LG_A$ 
extending the natural embedding
$\wh S\xrightarrow{\wh{\Int(g)}}\wh T\subset \wh G$.

Moreover, $\lsup Lj$ is uniquely determined up to
$\wh T$-conjugacy by the choice of
\begin{itemize}
\item the $\Gal_F$-pinning $(\wh B, \wh T, \{Y_\alpha\})$
of $\wh G$, and
\item the Borel $B_S$.
\end{itemize}

\begin{proof}
By Lemma \ref{lem:Lgrp-1},
it suffices to produce an $L$-embedding
for the Weil form of $L$-groups.
We refer the reader to \cite[Section 2.2]{Kal19a}
for the notion of ramified/unramified symmetric/asymmetric roots,
    and the definition of $\chi$-data.
Since $I_F$ stabilizes a set of positive roots (Definition \ref{def:max-unr-tori}),
there is no symmetric ramified root (the existence of symmetric ramified roots implies there exists a ramified quadratic extension $E/F$ such that the $\Gal(E/F)$-orbit of a positive root $\alpha$ is $\pm\alpha$).
Consequently, we can choose the unramified $\chi$-data
$\{\chi_\alpha\}$
such that $\chi_\alpha=1$ if $\alpha$ is asymmetric
and $\chi_\alpha$ is the unramified quadratic character
if otherwise.
In particular, all these characters $\chi_\alpha$
are valued in $\{\pm 1\}$.
We remark that unramified $\chi$-data are minimally ramified in the sense of \cite[Definition 4.6.1]{Kal19a}.

So it suffices to show the Langlands-Shelstad $L$-embedding in \cite[Section 2]{LS87}
for unramified $\chi$-data
    is defined over $A$.
The formula in \cite[Section 2.6]{LS87} for the $L$-embedding reads
$$
\xi(w)=r_p(w)n(\omega_T(\sigma))\times w
$$
(see loc. cit. for the notations).
Note that $n(\omega_T(\sigma))$ is a product of the image of
root vectors under the exponential map
(see \cite[Section 2.1]{LS87})
and the exponential map is well-defined for
split group schemes over an arbitrary base scheme
\cite[Theorem 4.1.4]{Crd11}.
The map $r_p(w)$ has explicit formula in \cite[Section 2.5]{LS87}
and it is valued in $\{\pm 1\}\otimes X^*(T)$
because all the characters $\chi_\alpha$ are valued in $\{\pm 1\}$.
Since $\{\pm 1\}\subset A^\times$ for an arbitrary finite ring $A$,
we have $\{\pm 1\}\otimes X^*(T)\subset A^\times\otimes X^*(T)=\wh T(A)$.
In particular, $\xi$ is defined over $A$.
The ``moreover'' part is the remark right above \cite[Paragraph 2.6.1]{LS87}.
\end{proof}

For lack of a better terminology,
we call $L$-embeddings $\lsup LS\to \lsup LG$
that appear in
the theorem above {\it canonical $L$-embeddings}.

Let $E$ be a finite extension of $\Qp$ with ring of integers $\cO$
and uniformizer $\varpi$.
Taking the inverse limit
of $\lsup LS(\cO/\varpi^n)\to \lsup LG(\cO/\varpi^n)$,
we get an $L$-embedding
$\lsup LS(\cO)\to \lsup LG(\cO)$.
Taking the union of
$\lsup LS(\cO)\to \lsup LG(\cO)$
as $E$ varies, we 
get an $L$-embedding
$\lsup LS(\bZp)\to \lsup LG(\bZp)$.
As we have remarked before,
we don't think an $L$-embedding
$\lsup LS(\bQp)\to \lsup LG(\bQp)$
exists in general (\cite{Kal21}).

\mysection{Langlands-Shelstad factorization of $L$-parameters}

\mysubsection{Theorem}
\label{thm:LS}
For each semisimple $L$-parameter $\rho: \Gal_F \to \lsup LG(\bFp)$,
there exists a maximally unramified tame torus $S$ of $G$ defined over $F$
and a canonical $L$-homomorphism $\lsup Lj:\lsupp LS\to \lsupp LG$
such that $\rho$ factors through $\lsup Lj$.

\vspace{3mm}


As a byproduct of the construction,
we can ensure
the ramification index of $S$ is equal to that of $G$
and that there exists a $\rho(I_F)$-stable Borel $\wh B\subset \wh G$
that contains $\wh S$.

\begin{proof}
Write $\sigma\in \Gal_F$ for a Frobenius element,
and write $\tau\in \Gal_F$ for an element generating the tame quotient of the inertia.

By Corollary \ref{cor:metacyclic} and Theorem \ref{thm:ab-qss},
there exists a maximal torus $\wh S\subset \wh G$
which is stable under $\rho(\sigma)$ and $\rho(\tau)$,
and a $\rho(\tau)$-stable Borel $\wh B$ of $\wh G$ containing $\wh S$.
Thus the torus $\wh S$, equipped with Galois action
$\gamma \cdot x:= \rho(\gamma) x \rho(\gamma)^{-1}$
($x\in S$, $\gamma\in \Gal_F$),
together with the tautological embedding
$\wh S\to \wh G$
is maximally unramified in the sense of Definition \ref{def:max-unr-tori-dual}.
By Proposition \ref{prop:dual-max-unr} and Theorem \ref{thm:kot},
$S$ exists. 
By Theorem \ref{thm:LS-maxunr}, $\lsup Lj$ exists.
Write $\varphi: \Gal_F\to \lsup LS(\bFp)$
the map $\gamma\mapsto 1\rtimes \gamma$.
Write $\rho':=\lsup Lj \circ \varphi$.
Note that for each $\gamma\in\Gal_F$,
$\rho(\gamma)\rho'(\gamma)^{-1}$ acts trivially on $\wh S$,
and thus $\rho(\gamma)\rho'(\gamma)^{-1}\in Z_{\wh G}(\wh S)=\wh S$.
So $\Img \rho \subset \lsup Lj(\lsup LS(\bFp))$.
Since $\lsup Lj$ is an embedding, we have a factorization as desired.

Write $e$ for the ramification index of $G$,
we have $\rho(\tau)^e(1\rtimes \tau^{-e})\in \wh G \cap N_{\lsupp LG}(\wh B) \cap N_{\lsupp LG}(\wh S)=\wh S$,
which implies the ramification index of $S$ is at most $e$.
\end{proof}

\mychapter{The Deligne-Lusztig map}~

\etocsettocdepth{2}
\localtableofcontents
\vspace{3mm}
\label{sec:DL}

The uniqueness of canonical $L$-embeddings
(see Theorem \ref{thm:LS-maxunr})
allows us to attach
characters of $F$-tori of $G$
to semisimple $L$-parameters for $G$.

We start with the following basic lemma,
which claims that
homomorphisms of tori
which respect $F$-points
are necessarily defined over $F$.

\mysubsection{Lemma}
\label{lem:unique-F-str}
Let $T$ and $S$ be two $F$-tori.
Let $f: T_{F^s}\to S_{F^s}$ be an isomorphism of $F^s$-tori
such that $f(T(F))=S(F)$.
Then $f$ is an $F$-isomorphism.

\begin{proof}
Let $L/F$ be a finite extension splitting both $T$ and $S$.
We may assume $f$ is an $L$-isomorphism between $T_L$ and $S_L$.
Write
\begin{align*}
j_T:& T \to \Res_{L/F} T_L\\
j_S:& S \to \Res_{L/F} S_L
\end{align*}
for the adjunction maps relating Weil restriction and base change.
On the level of functor of points, the map
$$
j_T(F): T(F) \to (\Res_{L/F}T_L)(F) = T(L\otimes_F F) = T(L)
$$
is the map induced by the inclusion $F\subset L$ (\cite[Proposition A.5.7]{CGP15}).
By assumption, we have
$\Res_{L/F}f(j_T(T(F))) = j_S(S(F))$.
By \cite[Proposition A.5.7]{CGP15},
$j_T$ is a closed immersion defined over $F$.
It is clear that $\Res_{L/F} f$ is an $F$-isomorphism.
Write $j(T)$ for the subtorus $\Res_{L/F} f(j_T(T)) \subset \Res_{L/F}S_L$,
and write $j(S)$ for the subtorus $j_S(S) \subset \Res_{L/F}S_L$.
Note that both $j(T)$ and $j(S)$ are $F$-subtorus of $\Res_{L/F}S_L$.
We have shown $j(T)(F) = j(S)(F)$.
By \cite[Corollary 13.3.10]{Sp98}, since $F$ is an infinite field,
$j(T)(F)$ is Zariski dense in $j(T)$;
and
$j(S)(F)$ is Zariski dense in $j(S)$.
Thus $j(S)=j(T)$.
Now $T \cong j(T) =j(S) \cong S$ as an $F$-torus.
\end{proof}

\mysection{Stable conjugacy}

Fix an $F$-pinning $(B, T, \{X_\alpha\})$
of $G$
and its dual pinning $(\wh B, \wh T, \{Y_\alpha\})$
once for all.

\mysubsection{Definition}
A {\it Deligne-Lusztig datum}
is a pair $(S, \chi)$
consisting of a maximally unramified maximal $F$-torus
$S\subset G$
and a character $\chi: S(F)\to \bFp^\times$.

Write $\Int(g)$
for the conjugation by $g$ map $x\mapsto g x g^{-1}$.

Two Deligne-Lusztig data (DL data for short)
$(S_1, \chi_1)$ and $(S_2, \chi_2)$
are said to be stably conjugate
if there exists an element $g\in G(F^s)$
such that $S_2(F) = g S_1(F) g^{-1}$
and $\chi_2 = \Int(g)_*\chi_1$.
By Lemma \ref{lem:unique-F-str},
$\Int(g)$ is an $F$-isomorphism.

A {\it based Deligne-Lusztig datum}
is a tuple $(S, \chi, B_{S})$
where $(S, \chi)$ is a Deligne-Lusztig datum
and $B_{S}\subset G_{F^s}$ is a Borel defined over $F^s$
containing $S$.

\mysubsection{The Deligne-Lusztig map}
By Theorem \ref{thm:LS-maxunr},
given a Borel pair $(S, B_S)$,
there is a canonical $L$-embedding
$\lsup Lj: \lsupp LS \to \lsupp LG$,
unique up to $\wh T$-conjugacy.
Given a DL datum $(S, \chi)$,
by the Local Langlands Correspondence for tori,
we can attach to it an $L$-homomorphism
$\rho_\chi: \Gal_F\to \lsup LS(\bFp)$,
well-defined up to $\wh S$-conjugacy.
Define the {\it Deligne-Lusztig map}
$$
\DL: (S, \chi, B_S) \mapsto \lsup Lj \circ \rho_\chi,
$$
which is well-defined up to $\wh T$-conjugacy.

\mysubsection{Lemma}
Let $(S, \chi)$ be a DL datum
and let $(S, \chi, B_S)$ and $(S, \chi, B_S')$
be two enhancements of $(S, \chi)$.
Then $\DL(S, \chi, B_S)$
and $\DL(S, \chi, B_S')$
are $\wh G$-conjugate to each other.

\begin{proof}
By \cite[Lemma 2.6.A]{LS87},
the $\wh G$-conjugacy class of
$\lsup Lj: \lsupp LS \to \lsupp LG$
does not depend on the choice of $B_S$.
See the proof of Theorem \ref{thm:LS-maxunr}
for characteristic $p$ issues.
\end{proof}

As a consequence of the lemma above,
for each DL datum $(S, \chi)$,
$\DL(S, \chi) := \DL(S, \chi, B_S)$
is well-defined up to $\wh G$-conjugacy.

\mysubsection{Lemma}
\label{lem:DL-1}
Let $(S, \chi)$ and $(S', \chi')$
be two stably conjugate DL data.
Then $\DL(S, \chi)$
and $\DL(S', \chi')$
are $\wh G$-conjugate.

\begin{proof}
Let $(S, \chi, B_S)$
be an enhancement of $(S, \chi)$.
Let $g\in G(F^s)$ be an element such that
$g S(F) g^{-1} = S'(F)$
and $\chi' = \Int(g)_*\chi$.
Set $B_S':= g B_S g^{-1}$.
By Lemma \ref{lem:unique-F-str},
$\Int(g): S\to S'$ is an $F$-isomorphism
and thus induces
an isomorphism $\lsup L\Int(g): \lsup LS'\xrightarrow{\cong} \lsup LS$
of $L$-groups.

By the functoriality of the LLC for tori,
$\rho_{\chi} = \lsup L \Int(g) \circ \rho_{\chi'}$.
Thus $\DL(S, \chi, B_S) = \DL(S', \chi', B_S')$.
The lemma now follows from the definition
of $\DL(S, \chi)$.
\end{proof}

\mysubsection{Proposition}
\label{prop:DL-1}
The map $\DL$
induces a surjective map
from the set of stable conjugacy classes of 
DL data $(S, \chi)$ to the set of semisimple $L$-parameters $\Gal_F\to \lsup LG(\bFp)$.

\begin{proof}
The well-definedness follows from Lemma \ref{lem:DL-1}
and the surjectivity follows from Theorem \ref{thm:LS}.
\end{proof}

\mysubsection{Failure of injectivity}~

The map $\DL$ in
Proposition \ref{prop:DL-1}
is not injective in general.

Consider the spherical $L$-parameter
$\rho:\Gal_{F}\to \GL_2$
for $\GL_2$
which is trivial on the inertial $I_F$
and sends a Frobenius element to
$
\begin{bmatrix}
0 & 1 \\
-1 & 0
\end{bmatrix}
$.
On the one hand,
$\rho$ factors through the $L$-group
of an elliptic maximal torus $S_e$ of $\GL_2$.
On the other hand,
since the matrix $\begin{bmatrix}
0 & 1 \\
-1 & 0
\end{bmatrix}$
is diagonalizable,
$\rho$ also factors through the $L$-group
of a split maximal torus $S_s$ of $\GL_2$.
Therefore, two DL data
$(S_e, \chi_e)$ and $(S_s, \chi_s)$
are both sent to $\rho$ under $\DL$.
However, since $S_e$ and $S_s$
are not $F$-isomorphic,
by Lemma \ref{lem:unique-F-str},
$(S_e, \chi_e)$ and $(S_s, \chi_s)$
are not stably conjugate.

Denote by $\uF$ the strict henselization of $F$.
Note that both $S_s$ and $S_e$
are maximally unramified maximal torus of $\GL_2$
(indeed, they are both unramified),
and they both become split after base change to $\uF$.
Since any two split maximal tori are rationally conjugate 
to each other, there exists
an element $g\in G(\uF)$
such that $g S_s(\uF) g^{-1} = S_e(\uF)$.
This observation suggests that 
the inertial version of $\DL$
has a chance to be a bijection.

\mysection{Characters of tori over finite fields}

\mysubsection{Lemma}
\label{lem:ca-tori-1}
Let $X$ be a finite free abelian group equipped
with a finite order automorphism $\pi$.
There exists isomorphisms
\begin{align*}
\frac{X}{(q-\pi)X} &\cong (X \otimes \QZp)^{q\pi^{-1}}\\
\frac{X^\vee}{(q-\pi)X^\vee} &\cong \Hom((X \otimes \QZp)^{q\pi^{-1}}, \QZp)
\end{align*}
where $X^\vee=\Hom(X,\bZ)$.

\begin{proof}
It is \cite[Proposition 3.2.2]{Ca93}
and \cite[Proposition 3.2.3]{Ca93}.
\end{proof}

\label{rem:fin-char}
Indeed, for each $X$ as in the lemma above,
there exists a (unique up to isomorphism) torus $\uT$ defined over $\bF_q$
such that $X = X_*(\uT)$ (Proposition \ref{prop:F-tori}),
and there is a isomorphism
$(X \otimes \QZp)^{q \pi^{-1}} 
\cong (X \otimes \bFp^\times)^{\pi^{-1}\otimes \Frob_q} = \uT(\bF_q)$.
By the lemma above, there are short exact sequences
\begin{equation}
\label{eq:fin-tori}
0 \to X_*(\uT) \xrightarrow{q-\pi} X_*(\uT) \xrightarrow{\Xi}
\uT(\bF_q) \to 0,\text{~and}
\end{equation}
\begin{equation}
\label{eq:fin-tori2}
0 \to X^*(\uT) \xrightarrow{q-\pi} X^*(\uT) \xrightarrow{\Xi}
\Hom(\uT(\bF_q), \bFp^\times) \to 0.
\end{equation}

The short exact sequence (\ref{eq:fin-tori})
enables us to ``forget the $\bF_q$-structure'' on $\uT$
when working with characters of $\uT(\bF_q)$.
Given a continuous homomorphism
$$
\chi: \uT(\bF_q) \to \bFp^\times,
$$
the composition $\chi \circ \Xi$
is a character of the free abelian group $X_*(\uT)$,
which is independent of the $\bF_q$-structure on $\uT$.
Moreover the composite $\chi \circ \Xi$ uniquely determines
the original character $\chi$, once the $\bF_q$-structure
on $\uT$ is specified.

We write $\uT_{\pi}(\bF_q)=\uT(\bF_q)$
to emphasize the $\bF_q$-points are taken
with respect to $\pi$.

\mysubsection{Example}
The $1$-dimensional torus $\uGm$ over $\bFp$
admits two $\bF_q$-structures $\pi$ and $\pi'$
that splits over $\bF_{q^2}$:
$\pi$ corresponds to the split $1$-dim $\bF_q$-torus
and $\pi'$ corresponds to the (unique up to isomorphism) nonsplit $1$-dim $\bF_q$-torus.
Let $$\chi:\uGm_{\pi}(\bF_q) \to \bFp^\times$$
and $$\chi':\uGm_{\pi'}(\bF_q)\to \bFp^\times$$
be trivial characters $x\mapsto 1$.
Even though $\uGm_{\pi}(\bF_q)$ and $\uGm_{\pi'}(\bF_q)$
are distinct subsets of $\uGm(\bFp)$,
we have $\chi\circ\Xi = \chi' \circ \Xi' = 1$.

\mysubsection{Definition}
\label{def:fin-equiv}
Let $\uT$ be a torus over $\bFp$
equipped with two (possibly different) $\bF_q$-structures
$\pi$ and $\pi'$.
A character $\chi:\uT_{\pi}(\bF_q) \to \bFp^\times$
and 
a character $\chi':\uT_{\pi'}(\bF_q) \to \bFp^\times$
are said to be equivalent
if they define the same character of $X_*(\uT)$.
We write $\chi \cong \chi'$
if they are equivalent.

\mysubsection{Definition}
\label{def:chi-equiv}
Let $\uS_1$, $\uS_2$ be two tori defined over a field $\bF_q$.
Let $\uchi_1: \uS_1(\bF_q)\to \bFp^\times$
and $\uchi_2: \uS_2(\bF_q)\to \bFp^\times$
be characters.

A $\bFp$-isomorphism $f:(\uS_1)_{\bFp}\to (\uS_2)_{\bFp}$
is said to define an equivalence of $(\uS_1, \uchi_1)$ and $(\uS_2, \uchi_2)$
if 
 $f_*(\chi_1)\cong\chi_2$ in the sense of Definition \ref{def:fin-equiv};
we write $(\uS_1, \uchi_1) \cong_f (\uS_2, \uchi_2)$.

\mysubsection{Lemma}
\label{lem:ca-tori-0}
Let $X$ be a finite free abelian group.
We have canonical isomorphisms
\begin{align*}
\Hom_{\cts}(X \otimes \invlim{E/F\text{~unramified}}\kappa_E^\times,~\bFp^\times)
&=
\dirlim{E/F\text{~unramified}}\Hom_{\cts}(X \otimes \kappa_E^\times,~\bFp^\times)\\
&=
\dirlim{E/F\text{~unramified}}\Hom_{\cts}(X \otimes \cO_E^\times,~\bFp^\times).
\end{align*}
Here transition maps are norm maps.

\begin{proof}
It follows from the Stone duality that
in the category of profinite groups we have
$$
\Hom_{\cts}(X \otimes \invlim{E/F\text{~unramified}}\kappa_E^\times,~\bF_q^\times)
=
\dirlim{E/F\text{~unramified}}\Hom_{\cts}(X \otimes \kappa_E^\times,~\bF_q^\times).
$$
Since $X \otimes \invlim{E/F\text{~unramified}}\kappa_E^\times$
is a profinite group
and is thus a compact object,
we have
\begin{align*}
\Hom_{\cts}(X \otimes \invlim{E/F\text{~unramified}}\kappa_E^\times,~\bFp^\times)
&=
\dirlim{q}\Hom_{\cts}(X \otimes \invlim{E/F\text{~unramified}}\kappa_E^\times,~\bF_q^\times)\\
&=
\dirlim{q}\dirlim{E/F\text{~unramified}}\Hom_{\cts}(X \otimes \kappa_E^\times,~\bF_q^\times)\\
&=
\dirlim{E/F\text{~unramified}}\dirlim{q}\Hom_{\cts}(X \otimes \kappa_E^\times,~\bF_q^\times)\\
&=\dirlim{E/F\text{~unramified}}\Hom_{\cts}(X \otimes \kappa_E^\times,~\bFp^\times)\qedhere
\end{align*}
\end{proof}

\mysubsection{Definition}
Write $\kNm F$ for $$\invlim{E/F\text{~unramified}}\kappa_E^\times,$$
and
write $\ONm F$ for $$\invlim{E/F\text{~unramified}}\cO_E^\times.$$

\mysubsection{Corollary}
\label{cor:ca-tori-1}
Let $\uT$ be a $\kappa_F$-torus.

(1)
For each unramified extension $E/F$,
the map
$$
\Hom(\uT(\kappa_F), \bFp^\times)
\xrightarrow{\chi\mapsto \chi\circ \Nm_{\kappa_E/\kappa_F}}
\Hom(\uT(\kappa_E), \bFp^\times)
$$
is an inclusion.

(2)
We can identify
$\Hom(\uT(\kappa_F), \bFp^\times)$
as a subgroup of
$
\Hom(X_*(\uT) \otimes \kNm F, \bFp^\times)
$.
Under this identification, we have
$$
\Hom(\uT(\kappa_E), \bFp^\times)
=\Hom(X_*(\uT) \otimes \kNm F, \bFp^\times)^{\Gal_{\kappa_F}}.
$$

\begin{proof}
(1)
Write $\wh \uT$ for the dual torus of $\uT$.
By Lemma \ref{lem:ca-tori-1},
we have $\Hom(\uT(\kappa_E), \bFp^\times) \cong \wh \uT(\kappa_E)$
and $\Hom(\uT(\kappa_F), \bFp^\times) \cong \wh \uT(\kappa_F)$.
It is clear that $\uT(\kappa_F)=\uT(\kappa_E)^{\Gal(\kappa_E/\kappa_F)}$.

(2)
It follows from part (1) and Lemma \ref{lem:ca-tori-0}.
Note that for all $E$ sufficiently large (splitting $\uT$),
$\uT(\kappa_E) = X_*(uT)\otimes \kappa_E^\times$.
\end{proof}

\mysubsection{Corollary}
\label{cor:ca-tori-2}
Let $T$ be an unramified $F$-torus,
and write $T(F)^0$ for the Iwahori subgroup of $T(F)$.
We have
$$
\Hom(T(F)^0, \bFp^\times) = \Hom(X_*(T)\otimes \ONm F, \bFp^\times)^{\Gal_F}.
$$

\begin{proof}
It is equivalent to Corollary \ref{cor:ca-tori-1}.
\end{proof}

\mysubsection{Remark}
If we choose $T=\Gm$, then
Corollary \ref{cor:ca-tori-2} together with
the Hochschild-Serre spectral sequence
immediately implies that the Artin repository map
$\Hom(F^{\times}, \bFp^\times) \to \Hom(\Gal_F, \bFp^\times)$ is an isomorphism.
The profinite group $\ONm F$ can be identified with the abelianized inertia
of $\Gal_F$.

\mysection{Inertial refinement of the mod $p$ LLC for tame tori}~
\label{sec:ILLC}

In this subsection, we
fix an $F$-tori $T$, with splitting field $L$.
The local Langlands correspondence for tori holds for
all divisible coefficients.

\mysubsection{Theorem}(\cite{Ch20}, \cite[7.5]{Yu09})
\label{thm:llct}
If $F$ is a local field,
then there exists an isomorphism $$\beta_T:H^1_{\cts}(W_{L/F}, X^*(T)\otimes D)\cong \Hom_{\cts}(T(F), D)$$
for any divisible abelian topological group $D$.
Moreover, $\beta_T$ is additive functorial in $T$
(in the sense that $\beta_T$ is an additive natural transformation between additive functors).

\vspace{3mm}
We will now fix a torsion divisible group $D$
equipped with discrete topology.

\mysubsection{Lemma}
\label{lem:llct-0}
We have $H^1_{\cts}(W_{L/F}, X^*(T)\otimes D)=H^1_{\cts}(\Gal_F, X^*(T)\otimes D)$.

\begin{proof}
It is clear since $D$-valued $L$-parameters have finite image.
\end{proof}

We also need to understand how the LLC for tori behaves 
under base change.

\mysubsection{Proposition}
\label{prop:llct-0}
Let $E\subset F^s$ be a finite extension of $F$,
and write $\Nm_{E/F}:E \to F$ for the norm map.
There exists a commutative diagram
$$
\xymatrix{
H^1_{\cts}(\Gal_F, X^*(T)\otimes D)
\ar[d]^{\rho:\mapsto \rho|_{\Gal_E}}
\ar[r] &
\Hom_{\cts}(T(E), D)\ar[d]^{\chi:\mapsto \chi \circ \Nm_{E/F}}
\\
H^1_{\cts}(\Gal_E, X^*(T)\otimes D)
\ar[r] &
\Hom_{\cts}(T(F), D)
}
$$
\begin{proof}
Let $T':=\Res_{E/F}T$.
Let $f: T'\to T$ be the norm morphism.
By unravelling definitions,
$f$ induces the restriction map $\rho \mapsto \rho|_{\Gal_E}$
for $L$-parameters,
and induces the norm map $\chi\mapsto \chi\circ \Nm_{E/F}$
for characters of tori.
The proposition is a special case of the functoriality of $\beta_T$.
\end{proof}

\mysubsection{Definition}
\label{def:llct-0}
Let $\rho\in H^1_{\cts}(\Gal_F, X^*(T)\otimes D)$
be an $L$-parameter.
The {\it inertial type} of $\rho$
is defined to be the image of $\rho$ in
$\dirlim{E/F\text{~finite unramified}}
\Hom_{\cts}(T(E), D)$,
and is denoted by $\beta_{T,I}(\rho)$.

\vspace{3mm}
We summarize basic results on integral models of tori
as follows.

\mysubsection{Proposition}
\label{prop:tori-bd}
(1) There exists a unique maximal compact subgroup
$T(F)^1\subset T(F)$.

(2) We have
$$
T(F)^1 = \{x\in T(F)| |\chi(x)|_p=0 \text{~for all~}\chi\in X^*(T)\}.
$$

(3)
Let $\Tft{T}$ be the {\it ft-N\'eron model} of $T$.
We have
$$
T(F)^1 = \Tft{T}(\cO_F).
$$

(4)
Let $\To{T}$ be the {\it connected N\'eron model} of $T$,
and let $T(F)^0$ be the Iwahori subgroup of $T(F)$.
We have
$$
T(F)^0 = \To{T}(\cO_F).
$$

(5)
Let $T'\to T$ be a morphism of $F$-tori. Then
$T'(F)\to T(F)$ maps $T'(F)^1$ to $T(F)^1$,
and maps $T'(F)^0$ to $T(F)^0$,

(6)
There exists a finite unramified extension
$E/F$ such that
$T(E)/T(E)^0 = T(\breve{F})/T(\breve{F})^0$.

\begin{proof}
(1) and (2): See \cite[Proposition 2.5.8]{KP22}.

(3): It is \cite[Proposition B.7.2]{KP22}.

(3): It is \cite[Proposition B.8.7]{KP22}.

(5): It is \cite[Proposition 2.5.9]{KP22}
and \cite[Proposition 2.5.19]{KP22}.

(6): It follows from \cite[Corollary 11.1.6]{KP22}.
\end{proof}

\mysubsection{Lemma}
\label{lem:tori-llc-1}
We have
$$
\dirlim{E/F\text{~finite unramified}}
\Hom_{\cts}(T(E), D)
=
\dirlim{E/F\text{~finite unramified}}
\Hom_{\cts}(T(E)^0, D).
$$

\begin{proof}
By part (6) of Proposition \ref{prop:tori-bd},
there exists a finite unramified extension $E/F$
such that $T(E)/T(E)^0=T(\breve{F})/T(\breve{F})^0$.

Let $F_1/F$ be a finite unramified extension containing $E$,
and let $E_1/F_1$ be a finite unramified extension.
It suffices to show for $\chi_1, \chi_2:T(F_1)\to D$
such that $\chi_1|_{T(F_1)^0}=\chi_2|_{T(F_1)^0}$,
$\chi_1\circ \Nm_{E_1/F_1}=\chi_2 \circ \Nm_{E_1/F_1}$.

Write $\wt \chi_i = \chi_i\circ \Nm_{E/F}$, $i=1,2$.
Suppose $\chi_1|_{T(F_1)^0}=\chi_2|_{T(F_1)^0}$.
Let $x\in T(E_1)$ be an arbitrary element.
Since $T(E_1)/T(E_1)^0=T(F_1)/T(F_1)^0$,
there exists $y\in T(F_1)$ such that
$\frac{x}{y}\in T(E_1)^0$.

Part (5) of Proposition \ref{prop:tori-bd} applied to
$\Res_{E_1/F_1}(T_{F_1}) \xrightarrow{\Nm_{E_1/F_1}} T_{F_1}$
shows that 
$$
\wt \chi_1|_{T(E_1)^0}=\wt \chi_2|_{T(E_1)^0}.
$$
In particular,
$$
\wt\chi_1(x/y)
=\wt\chi_2(x/y).
$$
Equivalently
$$
\frac{\wt\chi_1(x)}{\wt \chi_2(x)}
= (\frac{\chi_1(y)}{\chi_2(y)})^{[E_1:F_1]}.
$$
Since $D$ is assumed to be a torsion divisible group,
there exists an integer $n$ such that $(\frac{\chi_1(y)}{\chi_2(y)})^n=1$.
Choose $E_1$ to be the unramified extension of $F_1$ of degree $n$ inside $F^s$,
and we have $\wt \chi_1(x) = \wt \chi_2(x)$.
Since the group $T(F_1)$ is a finitely generated abelian group,
we can choose $E_1/F_1$ such that
$\wt \chi_1(x) = \wt \chi_2(x)$ for all $x\in T(E_1)$.
\end{proof}

\mysubsection{Corollary}
\label{cor:llct-0}
If $D=\bFp^\times$, we have
$$
\dirlim{E/F\text{~finite unramified}}
\Hom_{\cts}(T(E), D)
=
\dirlim{E/F\text{~finite unramified}}
\Hom(\To{T}(\kappa_E), D).
$$
where $\kappa_E$ is the residue field of $E$.

\begin{proof}
Combine Lemma \ref{lem:tori-llc-1}
and part (4) of Proposition \ref{prop:tori-bd};
also note that 
$\Hom(\To{T}(\kappa_E), \bFp^\times)=\Hom_{\cts}(\To{T}(\cO_E), \bFp^\times)$.
\end{proof}

\mysubsection{Lemma}
\label{lem:llct-1}
Assume $D=\bFp^\times$.
For each $\rho\in H^1_{\cts}(\Gal_F, \wh T(\bFp))$,
$$\beta_{T,I}(\rho)\in \Img(\Hom(\To{T}(\kappa_F), D) \longrightarrow \dirlim{E/F\text{~finite unramified}}
\Hom(\To{T}(\kappa_E), D)).$$

\begin{proof}
We define a $\Gal_{\kappa_F}$-action
on $\dirlim{E/F\text{~finite unramified}}
\Hom(\To{T}(\kappa_E), D)$
as follows,
for $\sigma\in \Gal_{\kappa_F}$ and $\chi: \Hom(\To{T}(\kappa_E), D))$, set $\sigma\cdot \chi:= \chi\circ \sigma^{-1}$.

It is clear that $\beta_{T,I}(\rho)$ is a $\Gal_{\kappa_F}$-fixed point.
So it remains to show that 
\begin{align*}
&\Img(\Hom(\To{T}(\kappa_F), D) \to \dirlim{E/F\text{~finite unramified}}\Hom(\To{T}(\kappa_E), D)) \\=& 
(\dirlim{E/F\text{~finite unramified}}
\Hom(\To{T}(\kappa_E), D))^{\Gal_{\kappa_F}}.
\end{align*}
The map above is clearly a well-defined injection.
Let $\chi\in \Hom(\To{T}(\kappa_E), D)$ be a $\Gal_{\kappa_F}$-invariant character.
Write $\chi_0\in \Hom(\To{T}(\kappa_F), D)$
for the character $\frac{1}{[E:F]}\chi|_{\To{T}(\kappa_F)}$
    (note that $D$ is a divisible group).
It is clear that $\chi = \chi_0 \circ \Nm_{E/F}$.
\end{proof}

\mysubsection{Proposition}
\label{prop:llct-1}
(1)
There is a commutative diagram
$$
\xymatrix{
\Hom_{\cts}(T(F), \bFp^\times) \ar[r]^{\beta_T^{-1}} \ar[d] &
H^1_{\cts}(\Gal_F, \wh T(\bFp))\ar[d] \\
\Hom_{\cts}(T(F)^0, \bFp^\times)
\ar[r]^{^{\beta_{T,I}^{-1}}}&
H^1_{\cts}(I_F, \wh T(\bFp))
}
$$
where both vertical maps are restriction maps.

(2)
Write $T_0\subset T$ for the maximal unramified subtorus.
All arrows in the following commutative diagram
$$
\xymatrix{
\Hom_{\cts}(T(F)^0, \bFp^\times) \ar[r]_{\beta_{T,I}^{-1}}^*[@]{\cong} \ar[d]^*[@]{\cong} &
H^1_{\cts}(I_F, \wh T(\bFp))^{\Frob}\ar[d]^*[@]{\cong} \\
\Hom_{\cts}(T_0(F)^0, \bFp^\times)
\ar[r]_{^{\beta_{T,I}^{-1}}}^*[@]{\cong}&
H^1_{\cts}(I_F, \wh T_0(\bFp))^{\Frob}
}
$$
are isomorphisms.
Here $\Frob\in \Gal_F$ is a Frobenius element
(that is, a topological generator of $\Gal_F$ modulo $I_F$).

\begin{proof}
(1)
It suffices to show
the composition
$$
\Hom_{\cts}(T(F), \bFp) 
\to \dirlim{E/F\text{~finite unramified}}\Hom(\To{T}(\kappa_E), D))
\xrightarrow{\cong} H^1_{\cts}(I_F, \wh T(\bFp))
$$
factors through $\Hom_{\cts}(T(F)^0, \bFp)$,
which is exactly the content of the Lemma \ref{lem:llct-1}.

(2)
It suffices to show the vertical morphisms and the bottom horizontal morphism are isomorphisms.

The left vertical map is an isomorphism by \cite[Corollary B.7.12]{KP22}.
The right vertical map is a special case of Proposition \ref{prop:QF}
(which will be elaborated in \ref{par:QF}).

So it remains to analyze the bottom horizontal map.
By local class field theory (\cite[6.11]{Iw86}) and the fact that
abelianization commutes with colimits,
there is a canonical $\Frob$-equivariant isomorphism
$$
\kNm{F}=\invlim{E/F\text{~finite unramified}}\kappa_E^\times \cong I_{F}/P_{F}
$$
where $P_F$ is the wild inertia.
Since all $L$-parameters we consider are tamely ramified, we have
$$
H^1_{\cts}(I_F, \wh T_0(\bFp))
=
H^1_{\cts}(I_F/P_F, \wh T_0(\bFp))
=
H^1_{\cts}(\kNm{F}, \wh T_0(\bFp)).
$$
Since $T_0$ is an unramified torus,
$I_F$ acts on $\wh T_0(\bFp)$ trivially, and thus
we have $\Frob$-equivariant isomorphisms
\begin{align*}
H^1_{\cts}(\kNm{F}, \wh T_0(\bFp))
&=\Hom_{\cts}(\kNm{F}, \wh T_0(\bFp))\\
&=\Hom_{\cts}(\kNm{F}, \bFp^\times\otimes X^*(T_0))\\
&=\Hom_{\cts}(X_*(T_0)\otimes\kNm{F} , \bFp^\times).
\end{align*}
Write $\To{T}$ for the (connected) N\'eron model of $T^0$,
we have
\begin{align*}
\Hom_{\cts}(T_0(F)^0, \bFp^\times)
&=
\Hom_{\cts}(\To{T}(\cO_F), \bFp^\times)\\
&=
\Hom_{\cts}(X_*(T_0)\otimes\kNm{F} , \bFp^\times)
\end{align*}
where the last step is Corollary \ref{cor:ca-tori-1}.
\end{proof}

\mysection{Inertial refinement of the Deligne-Lusztig map}~
\label{}

\mysubsection{The restriction map}
Denote by $\SC{G}$ the set of
stable conjugacy classes
of Deligne-Lusztig data $(S, \chi)$ for $G$.

Let $E/F$ be a finite unramified extension.
There exists a map
$$
\SC{G} \to \SC{G_E}
$$
sending $(S, \chi)$ to $(S_E, \chi\circ\Nm_{E/F})$
where $\Nm_{E/F}:S(E)\to S(F)$ is the norm map.

Define $$\SC{I_F, G}:=\Img(\SC{G} \to \dirlim{\text{$E/F$ unramified}}~\SC{G_E}).$$
Also denote by
$
\TI{G}
\subset H^1(I_F, \wh G(\bFp))_{\text{tame}}
$
the set of tame inertial Langlands parameters
$I_F\to \lsup LG(\bFp)$
that can be extended to $\Gal_F$.

With the new notations introduced,
Proposition \ref{prop:DL-1} can be rephrased as
the image of $\DL:\SC{G}\to H^1(\Gal_F, \wh G(\bFp))$
consists precisely the semisimple $L$-parameters.

\mysubsection{Lemma}
\label{lem:TI-0}
$\DL: \SC{I_F, G}\to \TI{G}$ is surjective.
So is 
$\dirlim{E/F\text{~unramified}}~\SC{G_E}\to H^1(I_F, \wh G(\bFp))_{\text{tame}}$.

\begin{proof}
By Lemma \ref{lem:Frob-ss-2},
a tame inertial type can be extended to 
a semisimple $L$-parameter.
The lemma now follows from Proposition \ref{prop:DL-1}.
\end{proof}

To show $\DL: \SC{I_F, G}\to \TI{G}$ is actually
a bijection,
we will need to give a combinatorial description of both
sets.

\mysubsection{Parahorics}
For each maximal $F$-split torus $S_s$,
we can attach an affine apartment $\cA(S_s)$ of
the Bruhat-Tits building $\cB(G)$.
For each vertex $x\in \cA(S_s)$,
a smooth group scheme $\cG^\circ_x$
with generic fiber $G$ and connected special fiber
such that $\cG_x^{\circ}(\cO)\subset \Fix(x)$.
There exists a closed $\cO$-split
torus $\cS_s\subset \cG_x^0$
with generic fiber $S_s$;
the special fiber $\bar\cS_s$ of $\cS_s$
is a maximal $\kappa_F$-split torus
in the special fiber of $\cG_x^0$
(\cite[4.1.20]{KP22}).

\mysubsection{Notation}
Recall that we fixed an $F$-pinning $(B, T, \{X_\alpha\})$
of $G$.
We can choose $T$ so that it is maximally unramified
and maximally split since $G$ is quasi-split.
Write $T_s\subset T$ for the maximal $F$-split subtorus
and write $T_0\subset T$ for the maximal unramified subtorus.
The Chevalley valuation associated to the
pinning $(B, T, \{X_\alpha\})$
determines a superspecial vertex (\cite[Definition 3.4.8]{Kal19a}) in the apartment $\cA(T_s)$.

Let $\cT_0$ be the maximal unramified $\cO$-torus
of $\cG^\circ_x$ with generic fiber $T_0$.
Write $\uT$ for the special fiber of $\cT_0$.
Also write $\uG$ for the reductive quotient of
(the special fiber of) $\cG^\circ_x$.

\mysubsection{Definition}
\label{def:iDLd}
An {\it inertial Deligne-Lusztig datum}
is a pair $(\uS, \uchi)$
where $\uS$ is a maximal $\kappa_F$-torus of $\uG$
and $\uchi:\uS(\kappa_F) \to \bFp^\times$
is a character.

A {\it based inertial Deligne-Lusztig datum}
is a pair $(\uS, \uchi, \uB_{\uS})$
where $(\uS, \uchi)$ is an inertial Deligne-Lusztig datum and $\uB_{\uS}\subset \uG$
is a Borel defined over $\bFp$ containing $\uS$.

Two inertial Deligne-Lusztig data
$(\uS, \uchi)$ and $(\uS', \uchi')$
are said to be
{\it geometrically conjugate}
if there exists an element $g\in \uG(\bFp)$
    such that 
$(\uS, \uchi)\cong_{\Int(g)}(\uS', \uchi')$
in the sense of Definition \ref{def:chi-equiv}.

Denote by $\SC{\uG}$
the set of geometric conjugacy classes of inertial
Deligne-Lusztig data.

\mysubsection{Proposition}
There is a natural bijection
$$
\dirlim{E/F\text{~unramified}}~\SC{\uG_{\kappa_E}} \to \dirlim{E/F\text{~unramified}}~\SC{G_E}.
$$

\begin{proof}
Let $(\uS, \uchi)\in \SC{\uG}$
be an inertial Deligne-Lusztig datum.
By \cite[Proposition 8.2.1 (1)]{KP22},
there exists a closed unramified $\cO$-torus $\cS_0$
of $\cG^0_x$ whose special fiber is $\uS$.
The generic fiber $S_0$ of $\cS_0$
is a maximal unramified torus of $G$
and the centralizer $S$ of $S_0$ in $G$
is a maximally unramified maximal $F$-torus.
Write $\chi_0^0$ for the inflated character
$$
\chi_0^0:S_0(F)^0\to \uS(\kappa_F) \xrightarrow{\chi_0} \bFp^\times.
$$
By part (2) of Proposition \ref{prop:llct-1},
the character $\chi_0^0$
can be extend uniquely to a character
$\chi^0: S(F)^0\to \bFp^\times$.
Since $\bFp^\times$ is a divisible group,
$\Hom(-, \bFp^\times)$ is exact,
the character $\chi^0$
can be extended to $\chi: S(F)\to \bFp^\times$.
The pair $(S, \chi)$ is a Deligne-Lusztig datum.

We first show the well-definedness.
Suppose $(\uS, \uchi)$ and $(\uS', \uchi')$
are geometrically conjugate.
By \cite[Proposition 8.2.1 (5)]{KP22},
there exists an element $g\in \cG^\circ_x(\cO_\uF)\subset G(\uF)$
such that $g S_0 g^{-1} = S_0'$
and $\Int(g)_*\uchi\cong \uchi'$.
Since $S$ (resp. $S'$) is the centralizer of
$S_0$ (resp. $S_0'$),
we have $g S g^{-1} = S'$.
Let $E/F$ be an unramified extension
such that $g\in G(E)$ and both
$\uS_{\kappa_E}, \uS_{\kappa_E}'$
are split.
Then we have $\Int(g)_*\uchi= \uchi'$
and $g S(E)g^{-1} = S(E)'$.
By Corollary \ref{cor:llct-0},
there exists a finite unramified extension $E'/E$
such that $\Int(g)_*\uchi\circ\Nm_{E'/E} = \uchi'\circ\Nm_{E'/E}$.
As a consequence, $(S, \chi)$ and $(S', \chi')$
are stably conjugate after restricting to $E'$.
Corollary \ref{cor:llct-0} also ensures
the map is injective.
Suppose $(S, \chi)$ and $(S', \chi')$
are stably conjugate after restricting to some $E$.
Since the apartment of $S$ and of $S'$
both contains the vertex $x$,
by \cite[Lemma 3.4.12]{Kal19a}, there exists
an element $g\in \cG^\circ_x(\cO_{\uF})$
such that $g S g^{-1} = S'$.
Corollary \ref{cor:llct-0} implies
$(\uS, \uchi)$ and $(\uS', \uchi')$
are geometrically conjugate (under the reduction $\bar g$ of $g$) after restricting to some finite unramified extension $E'/E$.

Next, we show the surjectivity of the map.
Let $S$, $\uS$ and $S_0$ be as in the first paragraph
of the proof, and
let $(S', \chi')$ be an arbitrary Deligne-Lusztig datum.
Since both $S$ and $S'$ are both maximally unramified
(i.e. they become split after base change to $\uF$),
there exists an element $g\in G(\uF)$
such that $g S' g^{-1} =S$.
We may assume $g\in G(E)$ for some finite unramified extension $E/F$.
Write $\chi_1$ for $\Int(g)_*\chi'$.
The restriction to $S_0(E)^0$
if $\chi_1$ factors through
$\uchi_1:\uS(\kappa_E)\to \bFp^\times$.
The inertia Deligne-Lusztig datum $(\uS_{\kappa_E}, \uchi_1)$
is mapped to the equivalence class of $(S', \chi')$.
\end{proof}

\mysubsection{Root system for the reductive quotient
of a superspecial parahoric}
The relative root system $\Phi_{\uF}(T_0, G)$ of $G$
with respect to $T_0$
is not reduced in general.
The absolute root system $\Phi_{\bFp}(\uT, \uG)$ of $\uG$
is a reduced modification
of the possibly non-reduced root system
$\Phi_{\uF}(T_0, G)$
(see \cite[Proposition 8.4.8]{KP22}).
In particular, $\Phi_{\uF}(T_0, G)$
and $\Phi_{\bFp}(\uT, \uG)$ have the same Weyl group
(see \cite[Remark 3.4]{Ha18}),
and their Weyl group
can be identified with
the $I_F$-invariant subgroup
(which we denote by $\Omega^\theta$)
of the absolute Weyl group
$\Omega:=N_{T(F^s)}(G(F^s))/T(F^s)$
(\cite[Lemma 4.2]{Ha15}).

\mysubsection{Theorem}
\label{thm:DL-SC}
(1)
Let $E/F$ be a finite unramified extension.
The restriction functor
$$
\SC{\uG} \to \SC{\uG_{\kappa_E}}
$$
sending $(\uS, \uchi)$ 
to $(\uS_{\kappa_E}, \uchi \circ \Nm_{\kappa_E/\kappa_F})$
is injective.

(2)
The set $\SC{\uG}$
is in natural bijection with
$$
((X^*(\uT) \otimes \QZp)/ \Omega^\theta)^{\Frob}
\cong
((X^*(\uT) \otimes \bFp^\times)/ \Omega^\theta)^{\Frob}
$$
where 
$\Frob$ is the Frobenius map
on $\uT$ corresponding to the standard $\kappa_F$-rational structure.

(3)
The set $\dirlim{E/F\text{~unramified}}~\SC{\uG_{\kappa_E}}$
is in natural bijection with
$$
(X^*(\uT) \otimes \QZp)/ \Omega^\theta\cong
(X^*(\uT) \otimes \bFp^\times)/ \Omega^\theta.
$$

\begin{proof}
(1)
It is \cite[Proposition 5.4]{DL76}.

(2, 3)
It is \cite[Proposition 5.7]{DL76}.
\end{proof}

\mysubsection{Summary}
\label{par:TI-0}
We have established the following diagram:
$$
\xymatrix{
    \SC{\uG} \ar@{^{(}->}[d] \ar@{^{(}->}[r]&
    \SC{I_F, G} \ar@{^{(}->}[d] \ar@{->>}[r]&
    \TI{G}\ar@{^{(}->}[d]
    \\
\dirlim{E/F\text{~unramified}}~\SC{\uG_{\kappa_E}}
\ar@{=}[d]
\ar@{=}[r]&
\dirlim{E/F\text{~unramified}}~\SC{G_E}
\ar@{->>}[r]&
H^1(I_F, \wh G(\bFp))_{\text{tame}}
&
\\
(X^*(\uT) \otimes \QZp)/ \Omega^\theta
}
$$

\mysection{Tame types and semisimple $L$-parameters}~

Let $q:=p^f$ be a $p$-power integer.
Denote by $\Frob_q:x\mapsto x^q$ the relative $q$-Frobenius map.

In the rest of this subsection, the dual groups $\wh G$ and $\lsupp LG$ 
are always defined over $\bFp$.
Write $L$ for the splitting field of $G$.
We regard $\Gal(L/F)$ as a subgroup of $\lsupp LG=\wh G\rtimes \Gal(L/F)$.
Fix a Frobenius element $\sigma\in \Gal_F$
and denote its image in $\lsupp LG$ and $\Gal(L/F)$ by $\bar\sigma$;
also fix a generator $\theta$ of the tame inertia of $\Gal_F$
and denote its image in $\lsupp LG$ and $\Gal(L/F)$ by $\bar\theta$.
Denote by $e$ the ramification index of $\Gal(L/F)$.

Define the {\it twisted $q$-Frobenius map}
$$F_\theta^\sigma: \wh G\to \wh G, \hspace{5mm}g\mapsto \bar\sigma^{-1}(\prod_{i=0}^{q-1} \bar \theta^i g \bar \theta^{-i})\bar\sigma.$$

\mysubsection{Definition}
A semisimple $\bar\theta$-twisted conjugacy class $[s]_{\bar\theta}$ of $\wh G$
is said to be {\it $F^\sigma_\theta$-stable}
if $s\in [s]_{\bar\theta}$ implies $F^\sigma_\theta(s)\in [s]_{\bar \theta}$.

\mysubsection{Lemma}
\label{lem:Frob-ss-1}
A semisimple $\bar\theta$-twisted conjugacy class $[s]_{\bar\theta}$ is
$F^\sigma_\theta$-stable
if and only if for each representative $s$ of $[s]_{\bar\theta}$,
there exists a tamely ramified $L$-parameter
$\rho:\Gal_F\to \lsupp LG$
such that $\rho(\theta)=s\bar \theta$.

\begin{proof}
Unravel the definitions.
\end{proof}

\mysubsection{Corollary}
\label{cor:Frob-ss-1}
A semisimple $\bar\theta$-twisted conjugacy class $[s]_{\bar\theta}$ is
$F^\sigma_\theta$-stable
if and only if for some representative $s$ of $[s]_{\bar\theta}$,
there exists a semisimple $L$-parameter
$\rho:\Gal_F\to \lsupp LG$
such that $\rho(\theta)=s\bar \theta$.

\begin{proof}
Combine Lemma \ref{lem:Frob-ss-1} and Lemma \ref{lem:Frob-ss-2}.
\end{proof}

\mysubsection{Proposition}
\label{par:QF}
(1) 
The map $F^\sigma_\theta$ is a $\theta$-twisted Frobenius endomorphism
in the sense of \ref{def:QF}.

(2)
There is a one-to-one correspondence between
$\bar\theta$-twisted semisimple conjugacy classes in $\wh G$
and $\wh G$-conjugacy classes of $L$-parameters $I_F \to \lsupp LG$,
given by $[s]_{\bar\theta} \mapsto (\theta\mapsto s \bar \theta)$.

(3)
There is a one-to-one correspondence between
$F^\sigma_\theta$-stable
$\bar\theta$-twisted semisimple conjugacy classes in $\wh G$
and $\wh G$-conjugacy classes of $L$-parameters $I_F \to \lsupp LG$
that can be extended to a semisimple $L$-parameter
$\rho:\Gal_F\to \lsupp LG$.

\begin{proof}
(1), (2): Unravel the definitions.

(3): It is Corollary \ref{cor:Frob-ss-1}.
\end{proof}

\mysubsection{Corollary}
\label{cor:DL}
In the context of Paragraph \ref{par:TI-0},
there are natural bijections
$$
\SC{\uG} \cong \SC{I_F, G} \cong \TI{G}.
$$

\begin{proof}
By Paragraph \ref{par:TI-0},
it suffices to show $\SC{\uG}\to \TI{G}$
is bijective.
By Theorem \ref{thm:DL-SC}, we have
$\SC{\uG}\cong ((X^*(\uT) \otimes \bFp^\times)/ \Omega^{\theta})^{\Frob}$.
By Proposition \ref{prop:QF}, we have
$\TI{G}\cong (X_*(\wh T)_{\theta,\Tf}\otimes \bFp^\times/\Omega^{\theta})^{\varphi\otimes \Frob_q}$.
It remains to show $X_*(\wh T)_{\theta,\Tf}$ and $X^*(\uT)=X^*(T_0)$
are canonically identified.
Note that $X_*(\wh T)_{\theta,\Tf}$
is by definition 
the maximal unramified finite free quotient
of $X_*(\wh T)$.
Since $X_*(T_0)$ is the maximal unramified
subgroup of $X_*(T)$,
by duality, $X^*(T_0)$
is the maximal unramified finite free quotient
of $X^*(T)$.
Since $X^*(T)=X_*(\wh T)$,
$X_*(\wh T)_{\theta,\Tf}$ and $X^*(\uT)$
are unramified.
\end{proof}

\mychapter{Digression: de Rham lifts of semisimple mod $p$ $L$-parameters of regular Hodge type}~

\etocsettocdepth{2}
\localtableofcontents
\vspace{3mm}

To study the Emerton-Gee stacks,
it is important to construct regular de Rham
lifts of all mod $p$ $L$-parameters.

For semisimple $L$-parameters,
such lifts can be easily constructed via
the Langlands-Shelstad factorization.
Recall that if $S\subset G$ is a maximally
unramified $F$-torus,
then there exists an $L$-embedding
$\lsup Lj: \lsup LS(\bZp) \to \lsup LG(\bZp)$
(see the remarks after Theorem \ref{thm:LS-maxunr}).

\mysection{Basic facts about de Rham $L$-parameters}~

Let $\Lambda\supset \bZ_p$
be a discrete valuation ring.

\mysubsection{Definition}
An $L$-parameter $\rho:\Gal_F \to \lsup LG(\Lambda)$
is said to be semistable (resp. crystalline)
if for some closed embedding of algebraic groups
$\lsup LG\hookrightarrow\GL_d$
the composite $\Gal_F\to \GL_d(\Lambda)$
is semistable (resp. crystalline).

An $L$-parameter $\rho:\Gal_F \to \lsup LG(\Lambda)$
is said to be potentially semistable (resp. crystalline)
if there exists a finite extension $E/F$
    such that $\rho|_{\Gal_E}$
is semistable (resp. crystalline).

\mysubsection{Lemma}
\label{lem:ss-inertia}
An $L$-parameter $\rho:\Gal_F\to \lsup LG(\Lambda)$
is semistable (resp. crystalline)
if and only if its restriction to inertia
$\rho|_{I_F}:I_F\to \lsup LG(\Lambda)$
is semistable (resp. crystalline).
Here we regard $I_F$ as the absolute Galois
group of the strict henselization $\uF$ of $F$.

\begin{proof}
It is \cite[Proposition 9.3.1]{BC08}.
\end{proof}


\mysubsection{Hodge-Tate theory}
\label{par:HT}
The main reference is \cite[Section 5]{L22} and \cite[Section 2.8]{BG19}.
Write $E$ for $\Lambda[1/p]$.
Let $\rho: \Gal_F\to \lsup LG(E)$ be a potentially
semistable $L$-parameter.
Then we can associate to $\rho$
a tuple of cocharacter
$\HT^\iota(\rho)$ of $\wh G_{\bC}$
for each embedding $\iota:E\hookrightarrow \bC$.
For each embedding $f:\lsup LG\to \lsup LH$,
we have $\HT^\iota(f\circ\rho) = f\circ \HT^\iota(\rho)$.

We can regard $\HT_{\bC}(\rho):=\boxtimes_{\iota:E\hookrightarrow \bC}\HT^\iota(\rho)$
as a cocharacter of $\prod_{\iota}\wh G_{\bC}\cong \Res_{\bC\otimes_{\Qp} E/\bC}\wh G$.
By Tannakian formalism, a cocharacter is equivalent to an exact tensor grading.
Write $\dR_{\bC}(\rho)$ for the canonical exact tensor filtration
associated to $\HT_{\bC}(\rho)$.
When $\rho$ is potentially semistable (say $\rho|_{\Gal_{F'}}$ is semistable),
then $\rho$ corresponds to a filtered $(\varphi, N, \Gal(F'/F))$-module $D_0$ with $\lsup LG$-structure.
Write $F'_0\subset F'$ for the maximal subfield unramified over $\Qp$.
Note that $D:=D_0\otimes_{F'_0}F'$ is an $\lsup LG$-torsor over $F'\otimes_{\Qp}E$
equipped with a $\Gal(F'/F)$-stable exact tensor filtration $\dR_{F'}(\rho)$
which recovers $\dR_{\bC}(\rho)$ after base change to $\bC\otimes_{\Qp}E$.
On the other hand, $D$ descends to an $\lsup LG$-torsor
over $D_F$ over $F\otimes_{\Qp}E$
and $\dR_{F'}(\rho)$ descends to an exact tensor filtration $\dR_F(\rho)$
on $D_F$.
Following \cite{BG19}, we will call $\dR_F(\rho)$ the {\it Hodge type} of $\rho$.
The potentially semistable $L$-parameter $\rho$ is said to be
{\it of regular Hodge type} or {\it Hodge regular}
if the stabilizer of $\dR_F(\rho)$ is a Borel subgroup of (a form of) $\Res_{F\otimes_{\Qp}E/E}\lsup LG$.
Since the property of being a Borel subgroup descends along field extensions,
we see $\rho$ is Hodge regular if and only if
each $\HT^\iota(\rho)$ is a regular cocharacter 
in the cocharacter lattice $X_*(\wh T)$ (here $\wh T$ is a maximal torus of $\wh G$ containing the image of $\HT^\iota(\rho)$), that is, $\HT^\iota(\rho)$
is not killed by any root of $\wh G$ with respect to $\wh T$.

\mysection{The $p$-adic Hodge theoretic refinement of the LLC for tori}~
\label{subsec:Hodge-LLC}

The results in this subsection is standard.
The main reference is \cite{Se89}.

Let $T$ be an $F$-torus which splits over $L$.
Fix an algebraic closure $\bQp$
and
let $E\subset \bQp$ be a finite extension of $L$.

Since $\bZp^\times$ is a divisible abelian group,
by \ref{thm:llct}, there exists a functorial isomorphism
$$\beta_T:H^1_{\cts}(W_{L/F}, X^*(T)\otimes \bZp^\times)\cong \Hom_{\cts}(T(F), \bZp^\times).$$
Since we are working with integral coefficients, the Galois form of $L$-parameters
and the Weil form of $L$-parameters are equivalent,
so $H^1_{\cts}(W_{L/F}, X^*(T)\otimes \bZp^\times)=H^1_{\cts}(\Gal_F, \wh T(\bZp))$
(see \ref{lem:Lgrp-1}).

\mysubsection{Locally algebraic characters}
A character $\chi\in \Hom_{\cts}(T(F), E^\times)$
is said to be {\it algebraic}
if there exists an algebraic character
$\psi\in \Hom(\Res_{F/\Qp}T, \Res_{E/\Qp}\Gm)$
such that $\chi = \psi(\Qp): T(F)\to \Gm(E)$.

A character $\chi\in \Hom_{\cts}(T(F), E^\times)$
is said to be {\it locally algebraic}
if $\chi$ coincides with some algebraic character
    in an open neighborhood of $1$.

\noindent
{\bf Facts} (See, for example, \cite[Appendix B]{C11})
$\chi$ is locally algebraic if and only if $\beta_T^{-1}(\chi)$ is Hodge-Tate,
if and only if $\beta_T^{-1}(\chi)$
is potentially semistable,
and if and only if $\beta_T^{-1}(\chi)$
is potentially crystalline.

\vspace{3mm}

\noindent
{\bf Notation}
Write $\Hom_{\cts}(T(F), -)_{\la}$
for the locally algebraic subgroup
of $\Hom_{\cts}(T(F), -)$,
and write
$H^1_{\cts}(\Gal_F, X^*(T)\otimes -)_{\HT}$
for the Hodge-Tate subset of
$H^1_{\cts}(\Gal_F, X^*(T)\otimes -)$.

\mysubsection{Lemma}
The composition
$$
\Hom_{\cts}(T(F), \bZp^\times)_{\la}
\xrightarrow{\beta_T^{-1}}
H^1_{\cts}(\Gal_F, X^*(T)\otimes \bZp^\times)_{\HT}
\xrightarrow{\HT}\prod_{\iota:\bQp\hookrightarrow \bC}X_*(T)
$$
is a group homomorphism.

\begin{proof}
Since the formation of co-labelled Hodge-Tate cocharacter 
is insensitive to restriction of Galois groups,
the general case is reduced to the split tori case,
which is clear.
\end{proof}

\mysubsection{Definition}
Denote by $$\fH^\iota:\Hom_{\cts}(T(F), \bZp^\times)_{\la} \xrightarrow{\HT\circ\beta_T^{-1}}\prod_{\iota:\bQp\hookrightarrow \bC}X_*(T)\xrightarrow{\text{$\iota$-th component}} X_*(T)$$
the group homomorphism
attaching to a character of $T(F)$
its $\iota$-colabelled Hodge-Tate cocharacter.

\mysubsection{Lubin-Tate Galois characters}
\label{par:LT}
Write $\iota_0:E\to \bC$ 
for the distinguished embedding
(recall that $E$ is a subfield of the fixed algebraic closure $\bQp$).

The identity map $(\Gm)_E \xrightarrow{=}(\Gm)_E$
induces a $E^\times$-valued character
$$
\chi_{\LT}: \Gm(E) \to E^\times=\Gm(E)
$$
of $\Gm(E)$.
Under the Local Langlands for split tori (= Local Class Field Theory),
the associated $L$-parameter
$\rho_{\LT}: \Gal_E \to E^\times$
is the so-called {\it Lubin-Tate Galois character}.

The Hodge-Tate weights for Lubin-Tate characters
are computed by Serre.
We have
$\HT^\iota(\rho_{\LT}) = 0$
if $\iota\ne \iota_0$
    and 
$\HT^{\iota_0}(\rho_{\LT}) = -1$
(see, for example, \cite[Appendix A.3]{L22}).
Here we naturally label the cocharacters of $\Gm$ by integers.

The Lubin-Tate character satisfies $\rho_{\LT}(I_E)\subset \cO_E^\times$.
If we choose a uniformizer $\varpi_E$ of $E$,
and let $\chi_{\varpi}:\Gm(E)\to \cO_E^\times$
be the character which sends $\varpi_E$ to $1$
and agrees with $\chi_{LT}$ on $\cO_E^\times$.
Then $\chi_{\varpi}$ corresponds to an integral Galois character
$\rho_{\varpi}:\Gal_E\to \cO_E^\times$,
which is usually called {\it the fundamental character}.


\mysubsection{Construction of crystalline Galois characters}~

Let $x\in X_*(\wh T)=X^*(T)\cong \bZ^{\oplus \dim_{F} T}$ be an arbitrary cocharacter.

Since $T$ splits over $L$, we have
$T(L)=(L^\times)^{\times \dim_FT}$.
Write
$\wt\chi_x: T(L) \to \cO_L^\times$
be the locally algebraic character which
is the composition of $x$ and the fundamental character.
Set $\chi_x:=\wt \chi_x|_{T(F)}$.
Note that $\Gal(L/F)$ acts on $X_*(T)$.

\mysubsection{Lemma}
\label{lem:tori-lift}
Assume $L/F$ is a Galois extension.
We have $\fH^{\iota_0\circ \theta}(\chi_x) = -\theta^{-1}\cdot x$
for all $\theta\in \Gal(L/F)$.

\begin{proof}
Since the formation of co-labelled Hodge-Tate cocharacter
is insensitive to restriction of field,
it suffices to compute the co-labelled Hodge-Tate cocharacters
of the composite
$$
\Res_{L/F}T(F) \xrightarrow{\Nm_{L/F}}
T(F) \xrightarrow{\chi_x} L^\times.
$$
Note that
$$
\chi_x \circ \Nm_{L/F}
= \prod_{\sigma\in \Gal(L/F)}
\wt \chi_x\circ\sigma.
$$
We have
\begin{align*}
\fH^{\iota_0\circ \theta}(\chi_x \circ \Nm_{L/F})
&= \sum_{\sigma\in \Gal(L/F)}
\fH^{\iota_0\circ \theta}(\wt \chi_x\circ\sigma)\\
&= \sum_{\sigma\in \Gal(L/F)}
\fH^{\iota_0\circ \theta}(\sigma\circ\wt \chi_{\sigma^{-1}\cdot x})\\
&= \sum_{\theta\in \Gal(L/F)}
\fH^{\iota_0\circ \theta}(\sigma\circ\wt \chi_{\sigma^{-1}\cdot x})\\
&= \sum_{\theta\in \Gal(L/F)}
\fH^{\iota_0\circ (\theta\sigma^{-1})}(\wt \chi_{\sigma^{-1}\cdot x})\\
&=
\sum_{\theta\in \Gal(L/F)}
\begin{cases}
0 & \theta\sigma^{-1} \ne 1\\
-\sigma^{-1} \cdot x & \theta\sigma^{-1}=1
\end{cases}\\
&= -\theta^{-1}\cdot x.
\end{align*}
The fourth equality follows from \cite[Corollary 4, Appendix A.2]{L22},
and the fifth equality follows from
the Lubin-Tate Galois character computation
(\ref{par:LT}).
\end{proof}

\mysection{Existence of de Rham lifts of prescribed Hodge types}~

\mysubsection{Theorem}
\label{thm:dR-lift}
Let $G$ be a quasi-split tame group over $F$,
and let $\bar\rho: \Gal_F\to \lsup LG(\bFp)$
be a semisimple mod $p$ $L$-parameter.

There exists a maximal $F$-torus $S$ of $G$
such that for
each cocharacter $x\in X_*(\wh S)\subset X_*(\wh G)$,
there exists a potentially crystalline lift 
$\rho: \Gal_F\to \lsup LG(E)$ of $\bar\rho$
such that $\HT^{\iota_0\circ \theta}(\rho)=-\theta^{-1}\cdot x$ for all $\theta\in \Gal(E/F)$.
Here $E$ is a sufficiently large extension of $F$
containing $L$.

In particular, $\bar\rho$ admits a de Rham lift of regular Hodge type.

\begin{proof}
By Theorem \ref{thm:LS},
$\bar\rho$ factors through $\lsup LS(\bFp)$
for some maximally unramified torus
$S$ of $G$.
By possibly enlarging $L$,
we assume $L$ is a splitting field of $S$.

By the LLC for tori, $\bar\rho$
corresponds to a character
$\chi_{\bar\rho}: S(F)\to \bFp^\times$,
which admits a Teichm\"uller lift
$[\chi_{\bar\rho}]: S(F)\to W(\bFp)^\times$.
Applying the LLC for tori once again,
$[\chi_{\bar\rho}]: S(F)\to W(\bFp)^\times\subset \bZp^\times$
corresponds to a finite image
$L$-parameter $[\bar\rho]: \Gal_F \to \lsup LS(\bZp)$.
The lift $[\bar\rho]$ is potentially
crystalline of trivial Hodge type.

It remains to modify $[\bar\rho]$
so that it has the desired Hodge type.
By Lemma \ref{lem:tori-lift},
there exists a locally algebraic character
$\chi_x: S(F)\to \cO_L^\times$
such that $\fH^{\iota_0\circ \theta}(\chi_x) = -\theta^{-1}\cdot x$
for all $\theta\in \Gal(L/F)$.
Write $\bar\chi_x: S(F)\to \bFp^\times$
for the reduction mod $p$ of $\chi_x$.
The product
$$
\chi:=\chi_x [\bar\chi_x^{-1}\chi_{\bar\rho}]
:S(F)\to (LW(\bFp))^{\times}
$$
is a locally algebraic character.
By the LLC for tori,
$\chi$ corresponds to
an $L$-parameter
$\rho:\Gal_F\to \lsup LS(\cO_E)$
lifting $\bar\rho$.
Here $E$ is a sufficiently large coefficient field.

To show $\rho$ can be made Hodge regular,
by the discussion in \ref{par:HT},
it suffices to ensure
$\theta^{-1}\cdot x$ is a regular cocharacter
for all $\theta$.
Irregular cocharacters
lie on the wall of Weyl chambers.
The $\Gal(L/F)$-orbit of irregular cocharacters is contained in a finite union of hyperplanes of $X_*(\wh S)\otimes_{\bZ}\bR$.
Since it is impossible to cover all integral points of $\bR^{d}$ using a finite number of hyperplanes, there is a choice of $x$ which makes $\rho$ Hodge regular.
\end{proof}

\mychapter{Parahoric Serre weights}~

\etocsettocdepth{2}
\localtableofcontents

\mysection{Algebraic representations}~

\mysubsection{Simple $\uG$-modules}~
Let $\uG$ be a connected reductive group
over $\bFp$.

Let $(\uB, \uT)$ be a Borel pair of $\uG$.
For each dominant character $\lambda\in X^*(\uT)$,
there exists a simple $\uG$-module $L(\lambda)$
of highest weight $\lambda$.
Any simple $\uG$-module is isomorphic to exactly one such
$L(\lambda)$ (\cite[II.2.4]{Jan03}).

Let $(\uB', \uT')$ be another Borel pair of $\uG$.
There exists an element $g\in\uG$ such that
$(\uB', \uT')=g(\uB, \uT)g^{-1}$.
Then a simple $\uG$-module $V$
has highest weight $\lambda$ with respect to $(\uB, \uT)$
if and only if it has highest weight $\lambda\circ\Int(g^{-1})$
with respect to $(\uB', \uT')$.

\mysubsection{Simple $\uG^\bfF$-modules}
Suppose the derived subgroup of $\uG$
is simply-connected.
Equip $\uG$ with a $\bF_{p^r}$-structure.
Since all reductive groups over a finite field is quasi-split,
specifying a $\bF_{p^r}$-structure amounts to
specifying a finite order automorphism $\pi$ of the based root datum of $\uG$.
Such a $\pi$ induces a Frobenius map $\bfF:\uG\to \uG$.
Fix a $\bfF$-stable Borel pair $(\uB, \uT)$.
Simple $\uG^\bfF$-modules arise from restrictions of simple $\uG$-modules
(\cite[Appendix A]{Her09}).
Write 
$$X_r(\uT):=\{\lambda\in X^*(\uT)|0\le \langle \lambda, \alpha^\vee\rangle < p^r,~\forall \alpha\in \Delta(\uB, \uT)\}$$
and
$$
X^0(\uT):=\{\lambda\in X^*(\uT)|\langle \lambda, \alpha^\vee\rangle = 0,~ \forall \alpha\in R(\uB, \uT)\}.
$$
Then we have a bijection (\cite[Proposition A.1.3]{Her09})
\begin{align*}
\frac{X_r(\uT)}{(p^r-\pi)X^0(\uT)}  &\xrightarrow{\cong}
\{\text{Irreducible representations of $\uG(\bF_{p^r})=\uG^\bfF$}\}\\
\lambda& \mapsto L(\lambda)|_{\uG^\bfF}
\end{align*}

\mysubsection{Based inertial Deligne-Lusztig data and simple $\uG^{\bfF}$-modules}~
Recall that a based inertial Deligne-Lusztig datum is a tuple
$(\uS, \uchi, \uB_\uS)$ (Definition \ref{def:iDLd}).
Recall 
that there exists a short exact sequence (Equation (\ref{eq:fin-tori2}))
$$
0 \to X^*(\uS) \xrightarrow{p^r-\pi} X^*(\uS) \xrightarrow{\Xi}
\Hom(\uS^\bfF, \bFp^\times) \to 0.
$$

\mysubsection{Lemma}
\label{lem:serre-1}
Let $(\uS, \uchi, \uB_\uS)$ be a based inertial Deligne-Lusztig datum.
If $\uB_{\uS}$ is $\bfF^m$-stable, then
the map
$$
X_{rm}(\uS) \to \frac{X^*(\uS)}{(p^{rm}-\pi^m)X^*(\uS)}
$$
is surjective.

\begin{proof}
It is harmless to assume $m=1$.
Let $\{\alpha\}=\Delta:=\Delta(\uB_\uS, \uS)$ be the set of simple roots with respect to $\uB_\uS$.
Write $\{\alpha^\vee\}_{\alpha\in \Delta}$ for the set of coroots,
and write $\{\omega_\alpha\}_{\alpha\in \Delta}\subset X^*(\uS)$
for the fundamental weights; they form a basis of the weight lattice
that is dual to $\{\alpha^\vee\}$.

We can replace $\uG$ by its derived subgroup and thus assume $\uG$ is semisimple.
Since $\uG$ is a simply-connected semisimple group,
the character lattice of $\uG$ coincides with the weight lattice.
We say two characters $\lambda_1$ and $\lambda_2\in X^*(\uS)$
are equivalent if they have the same image in
$\frac{X^*(\uS)}{(p^{r}-\pi)X^*(\uS)}$.
Let $\lambda=\sum n_\alpha \omega_\alpha\in X^*(\uS)$,
where $n_\alpha$ are integers.

\vspace{3mm}

{\bf Claim} $\lambda$ is equivalent to $\lambda'=\sum n_\alpha' \omega_\alpha$
where all $n_\alpha' \ge 0$.
\begin{proof}
We have $\mu:=(p^r-1)\sum \omega_\alpha = \sum (p^r-\pi)\omega_\alpha\in (p^r-\pi)X_*(\uS)$.
For a sufficiently large integer $N$,
$\lambda'=N\mu + \lambda$ is a dominant character.
\end{proof}

Now we assume $\lambda$ is dominant.
Among all dominant characters equivalent to $\lambda$,
we assume $|\lambda|:=\sum n_\alpha\in \bZ_{\ge 0}$ is the smallest.

Next we show that $\lambda\in X_{p^r}(\uS)$.
Assume some $n_\beta \ge p^r$.
Write $n_\beta = s + p^r t$, $s, t\in \bZ_{\ge 0}$.
We have
$$
\sum n_\alpha\omega_\alpha \sim \sum_{\alpha\ne \beta, \pi \beta} n_\alpha \omega_\alpha
+ s \omega_\beta + (n_{\pi \beta}+ t)\omega_{\pi\beta}=:\lambda'.
$$
We have
$$
|\lambda| - |\lambda'| = n_\alpha -s -t =(p^r-1)t>0
$$
which contradicts the assumption that $|\lambda|$ is minimal among
all equivalent dominant characters.
\end{proof}

\mysubsection{Definition}
A based inertial Deligne-Lusztig datum
$(\uS, \uchi, \uB_{\uS})$
is said to be of {\it niveau $m$} if
$\uB_\uS$ is $\bfF^m$-stable.

\mysubsection{Proposition}
\label{prop:serre-1}
There exists a natural surjective map
from the set of equivalence classes
of irreducible representations
of $\uG^{\bfF^m}$
to the set of geometric conjugacy classes
of niveau-$m$ based inertial Deligne-Lusztig data.
This map does not depend on any choices.

\begin{proof}
An irreducible representation of $\uG^{\bfF^m}$
is the restriction of an irreducible
algebraic representation of highest weight $\lambda\in X_{rm}(\uT)$
where $(\uT, \uB)$ is a $\bfF^m$-stable Borel pair.
The proposition follows from Lemma \ref{lem:serre-1}.
\end{proof}

\mysection{Deligne-Lusztig representations}~

$\uG$ is a connected reductive group over $\bFp$ with simply-connected derived subgroup,
equipped with a $\bF_p$-structure
$\bfF:\uG\to\uG$.
Fix a $\bfF$-stable pinning $(\uB, \uT, \{u_\alpha\})$ of $\uG$,
which exists because all reductive groups over finite fields
are quasi-split.

\mysubsection{Herzig's presentation of Deligne-Lusztig datum}
\label{subsec:herzig-pre}
In \cite{Her09}, a tame type is described by a pair
$(w, \mu)\in \Omega(\uG, \uT) \times X^*(\uT)$.

Let $(\uS, \uchi, \uB_\uS)$ be a based inertial Deligne-Lusztig datum.
There exists an element $g\in \uG$ such that
$(\uB, \uT)=g(\uB_\uS, \uS)g^{-1}$.
Two different choices of $g$ differ by left translation by
an element of $\uT$, so we get a well-defined
identification of $X^*(\uT)\cong X^*(\uS)$.
By Equation (\ref{eq:fin-tori2}),
$\uchi\in X^*(\uS)/(p-\pi)X^*(\uS)$,
and we let $\mu$ be the element of $X^*(\uT)$ lifting $\uchi$.
Let $\Delta:=\Delta(\uB, \uT)\cong \Delta(\uB_\uS, \uS)$ be the set of simple roots.
Since $\uS$ is $\bfF$-stable,
$\bfF$ acts on $X^*(\uS)$;
however, since $\uB_\uS$ is not $\bfF$-stable in general,
$\Delta(\uB_\uS, \uS)$ is not $\bfF$-stable in general.
There exists an element $w\in \Omega(\uG, \uS)\cong \Omega(\uG, \uT)$ such that
$\bfF(\Delta(\uB_\uS, \uS)) = w \Delta(\uB_\uS, \uS)$.

The based inertial Deligne-Lusztig datum $(\uS, \uchi, \uB_\uS)$
is presented by the pair $(w, \mu)$
by many authors.

\mysubsection{Regular $p$-restricted weights}
\label{subsubsec:regular}
An element $\lambda\in X_1(\uT)$ is said to be {\it regular}
if $\langle \lambda, \alpha^\vee\rangle \in [0, p-1)$ for all $\alpha\in \Delta(\uB, \uT)$.

\mysubsection{Lemma}
\label{lem:serre-2}
The natural map in Proposition \ref{prop:serre-1}
is injective when restricted to
irreducible representations
of $\uG^\bfF$
of regular $p$-restricted highest weight.

\begin{proof}
We keep notations in the proof of Proposition \ref{prop:serre-1}.
Let $\lambda = n_\alpha \omega_\alpha$
and $\lambda'=n_\alpha'\omega_\alpha'$
such that
$\lambda-\lambda'\in (p-\pi)X^*(\uT)$.
Say
$$
\lambda-\lambda' = \sum c_\alpha(p\omega_\alpha-\omega_{\pi \alpha})
=\sum (p c_\alpha-c_{\pi^{-1}\alpha})\omega_\alpha.
$$
We have
$$
-p < p c_\alpha - c_{\pi^{-1}\alpha}
< p
$$
for all $\alpha\in \Delta(\uB, \uT)$.
If $c_\alpha$ is the largest integer
for various $\alpha$,
we have
$p >p c_\alpha - c_{\pi^{-1}\alpha}\ge (p-1)c_\alpha$
and thus $c_\alpha\le 1$.
Similarly,
$-p < (p-1)c_\alpha$
and thus $c_\alpha\in [-1,1]$.
We have
$p c_\alpha - c_{\pi^{-1}\alpha}\in \{p-1, 0, -p+1\}$.
In particular, if $n_\alpha\ne n_\alpha'$,
then $\{n_\alpha, n_\alpha'\}=\{0, p-1\}$.
\end{proof}

\mysubsection{Definition}
A based inertial Deligne-Lusztig datum
of niveau $1$ is said to be {\it regular}
if it has a Herzig presentation
$(w, \mu)$
where $\mu$ is a regular $p$-restricted
weight.

An irreducible representation of
$\uG^\bfF$ is said to be {\it regular}
if it is the restriction of
a simple $\uG$-module
$L(\mu)$ to $\uG^\bfF$
where $\mu$ is a regular $p$-restricted
weight.

\mysubsection{Proposition}
\label{prop:parahoric-based}
There is a natural bijection between
equivalence classes of
regular irreducible representations
of $\uG^\bfF$
and geometric conjugacy classes
of regular based inertial Deligne-Lusztig data of niveau $1$.

\begin{proof}
Combine Proposition \ref{prop:serre-1}
and Lemma \ref{lem:serre-2}.
\end{proof}

\mysubsection{Mod $p$ twisting element}
Following \cite{GHS}, we denote by
$\eta\in X^*(\uT)$
an element that is Frobenius-stable and 
$\langle \eta, \alpha^\vee\rangle = 1$
for all $\alpha\in \Delta(\uB, \uT)$.

\mysubsection{Herzig's reflection operator $\clR$}
We define an involution
operator on the set of
isomorphism classes of regular irreducible
representations of $\uG^\bfF$
by
$$
\clR: L(\mu) \mapsto L(w_0\cdot (\mu - p\eta))
$$
where $w_0$ is the longest Weyl group element and $\eta$ is a mod $p$ twisting element.

\mysubsection{Deligne-Lusztig induction}
Let $(\uS, \uchi)$ be an inertial Deligne-Lusztig datum.
Write
$\bar V(\uS, \uchi)$
for the reduction mod $p$
of the Deligne-Lusztig induction
$\varepsilon_{\uG}\varepsilon_{\uS}R^{\uchi}_{\uS}$
(see \cite[Section 4.1]{Her09} for unfamiliar notations.)

\mysection{A generalization of Herzig's recipe}~
\label{subsec:Serre}

Let $G$ be a quasi-split tamely ramified reductive group over $\Qp$.
Fix a $\Gal_{\Qp}$-stable pinning
$(B, T, \{X_\alpha\})$ of $G$
such that $T$ is maximally unramified,
which determines
a superspecial parahoric $\cG^\circ$
of $G$.
Write $\uG$ for the reductive quotient
of $\cG^\circ$.
The pinning of $G$
determines a Frobenius-stable pinning
$(\uB, \uT, \{u_\alpha\})$
of $\uG$.

\mysubsection{Definition}
A {\it Serre weight} for $G$
is an irreducible $\bFp$-representation
of $\cG$,
where $\cG$ is the maximally bounded subgroup of $G$
containing $\cG^\circ$.

A {\it parahoric Serre weight} for $G$
is an irreducible $\bFp$-representation of $\uG(\bF_p)$.

\mysubsection{Assumption}
Assume both
$G$ and $\uG$ admit a local twisting element,
and assume $\uG$ has a simply-connected derived subgroup.
More precisely,
there exists an element $\eta_{\Qp}\in X^*(T)^{\Gal_{\Qp}}$
such that $\langle \eta_{\Qp}, \alpha^\vee\rangle=1$
for all $\alpha\in \Delta(B, T)$;
and there exists an element $\eta_{\bF_p}\in X^*(\uT)^{\Gal_{\bF_p}}$
such that $\langle \eta_{\bF_p}, \alpha^\vee\rangle=1$
for all $\alpha\in \Delta(\uB, \uT)$.

\vspace{3mm}
Write $T_0\subset T$ for the maximally unramified subtorus.
The identification
$X^*(T_0)\cong X^*(\uT)$
allows us to define a reduction map
$$
X^*(T) \xrightarrow{\text{restriction}}
X^*(T_0)\to X^*(\uT).
$$
By abuse of notation,
we denote by $\eta_{\Qp}$
the image of $\eta_{\Qp}$
in $X^*(\uT)$.

\mysubsection{The speculative recipe}
Let $\tau:I_{\Qp}\to \lsup LG(\bFp)$
be a tame inertial $L$-parameter.
Define
$$
W^?(\tau):= \clR(\JH(\bar V(\DL^{-1}(\tau)) \otimes W(w_0(\eta_{\bF_p}-\eta_{\Qp})))).
$$
Here $\DL$ is the Deligne-Lusztig map
(see Section \ref{sec:DL}),
$w_0$ is the longest Weyl group element,
$W(-)$ denotes the Weyl module,
and $\JH(-)$ denotes
the set of Jordan-H\"older factors.

\begin{appendix}

\myappsection{Maximal tori of quasi-split groups}
\label{app:tori}
In this appendix, we clarify the relation between maximal tori of a quasi-split group and that of its Langlands dual group in the natural generality.

Let $F$ be a perfect field with a fixed separable closure $F^s$, and
let $G$ be a quasi-split group over $F$.
Let $k$ be an arbitrary algebraically closed field,
and write $\lsup LG$ for the Langlands dual group of $G$ defined over $k$.

We start with recalling the classification of tori.

\myappsubsection{Proposition}
\label{prop:F-tori}
The set of $F$-isomorphism classes of $F$-tori of dimension $n$ is in natural bijection
with group homomorphisms $\Gal_F \to \GL_n(\bZ) = \Aut((\Gm^{\oplus n})_k)$ having finite image, up to $\GL_n(\bZ)$-conjugacy.

\begin{proof}
By \cite[Lemma 7.1.1]{Crd11}, the $F$-isomorphism classes
of $F$-tori are classified by the cohomology set
$H^1(\Gal_F, \Aut(\Gm^{\oplus n}))=H^1(\Gal_F, \GL_n(\bZ))$.
The proposition follows from the fact that the Galois action
on $\GL_n(\bZ)$ is trivial.
\end{proof}

\myappsubsection{Definition}
\label{def:fmtori}
A {\it framed maximal torus} of $G$ is a torus defined over $F$ of dimension equal to the rank of $G$, together with an embedding $j: S\to G$ defined over $F^s$.

The Galois group $\Gal_F$ acts on the collection of framed maximal tori by $j \mapsto \sigma \circ j \circ \sigma^{-1}=:j^\sigma$, $\sigma\in \Gal_F$.

A {\it conjugacy class of framed maximal tori}
is a set $\{g \circ j \circ g^{-1}=:\Int(g) \circ j| g\in G(F^s)\}$
where $j$ is a framed maximal tori.

\myappsubsection{Theorem} (Kottwitz)
\label{thm:kot}
Let $J$ be a conjugacy class of framed maximal tori
which is stable under $\Gal_F$-action.
Then there exists some $j\in J$ which is an algebraic group homomorphism defined over $F$.

\begin{proof}
It is \cite[Corollary 2.2]{Kot82}.
\end{proof}

\myappsubsection{Definition}
\label{def:fmtori2}
On the dual side, we define a {\it framed maximal torus} of $\lsup LG$
to be
an $k$-torus $\wh S$ of dimension $\rk G$ together with
\begin{itemize}
\item 
a Galois action
$\psi:\Gal_F\to \Aut(\wh S)$ and
\item
an embedding $\hat \j: \wh S \to \wh G$.
\end{itemize}
Note that the embedding $\hat \j$ is arbitrary and does not have to respect Galois actions.
The Galois group $\Gal_F$ acts on the collection of framed maximal tori by $(\psi, \hat \j) \mapsto (\psi, \sigma \circ \hat \j \circ \psi(\sigma)^{-1})$, $\sigma\in \Gal_F$.

\myappsubsection{Proposition}
\label{prop:kel}
There is a natural one-to-one correspondence $(S, j)\mapsto (\psi_S, \hat \j)$ between $\Gal_F$-stable
conjugacy classes of framed maximal tori of $G$,
and $\Gal_F$-stable
conjugacy classes of framed maximal tori of $\lsup LG$ where
the map $S\mapsto \psi_S$ is the dual map of the map defined in Proposition \ref{prop:F-tori}.

\begin{proof}
See \cite[Section 5.1]{Kal19a} for the construction.
Note that loc. cit. requires that $\operatorname{char} k=0$,
which is not necessary because the equivalence is validated by
a Galois cohomology computation with $\Omega(T,G)$-coefficients,
where $\Omega(T,G)$ is the absolute Weyl group and does not depend
on $k$.
\end{proof}
%
\end{appendix}

\addcontentsline{toc}{section}{References}
\printbibliography
\end{document}